\documentclass[smallextended]{svjour3}      
\smartqed  
\usepackage{lineno}
\modulolinenumbers[5]

\usepackage{mathptmx}      

\usepackage[colorlinks=true]{hyperref}
\usepackage{xparse}
\usepackage{amsmath}
\usepackage{amsfonts}

\newcommand{\bxi}{\boldsymbol{\xi}}

\newcommand{\dVblank}{\, {\rm d}V}
\newcommand{\dvblank}{\, {\rm d}v}
\newcommand{\dX}{\, {\rm d}\mbf{X}}

\newcommand{\eps}{\varepsilon}

\NewDocumentCommand \dV{ o }{%
    \IfNoValueTF{#1}{\dVblank}%
    {
        \dV_{#1}
    }
}

\NewDocumentCommand \dv{ o }{%
    \IfNoValueTF{#1}{\dvblank}%
    {
        \dv_{#1}
    }
}

\NewDocumentCommand \mc{ m }{%
    \mathcal{#1}%
}

\NewDocumentCommand \an{ m }{%
    \langle {#1} \rangle%
}

\NewDocumentCommand \mbf{ m }{%
    \mathbf{#1}
}

\NewDocumentCommand \bbar{ m }{%
    \mbf{\bar{#1}}
}

\NewDocumentCommand \mbt{ m }{%
    \mbf{\tilde{#1}}
}

\NewDocumentCommand \bs{ m }{%
    \boldsymbol{#1}
}

\NewDocumentCommand \stateOne{ m }{
    \underline{\mbf{#1}}%
}

\NewDocumentCommand \stateTwo{ m o }{%
    \IfNoValueTF{#2}{\stateOne{#1}}
    {
    \stateOne{#1}\an{#2}%
    }
}

\NewDocumentCommand \s{ m o o }{%
    \IfNoValueTF{#3}{\stateTwo{#1}[#2]}%
    {
    \underline{\mbf{#1}}(#3)\an{#2}%
    }
}

\NewDocumentCommand \stateOneI{ m m }{
    \underline{#1}_{#2}%
}

\NewDocumentCommand \stateTwoI{ m m o }{%
    \IfNoValueTF{#3}{\stateOneI{#1}{#2}}
    {
    \stateOneI{#1}{#2}\an{#3}%
    }
}

\NewDocumentCommand \sI{ m m o o }{%
    \IfNoValueTF{#4}{\stateTwoI{#1}{#2}[#3]}%
    {
        \stateOneI{#1}{#2}(\mbf{#4})\an{#3}
    }
}

\NewDocumentCommand \snsI{ m m o }{%
    \IfNoValueTF{#3}{\underline{\mbf{#1}}_{#2}}%
    {
      \underline{\mbf{#1}}_{#2}(#3)%
    }
}

\NewDocumentCommand \stateOneN{ m }{
    \underline{#1}%
}

\NewDocumentCommand \stateTwoN{ m o }{%
    \IfNoValueTF{#2}{\stateOneN{#1}}
    {
    \stateOneN{#1}\an{#2}%
    }
}

\NewDocumentCommand \sN{ m o o }{%
    \IfNoValueTF{#3}{\stateTwoN{#1}[#2]}%
    {
    \underline{#1}(#3)\an{#2}%
    }
}

\NewDocumentCommand\presuper{ m m }{%
  {\mathop{}%
   \mathopen{\vphantom{#2}}^{#1}%
   \kern-1\scriptspace%
   #2}
}

\NewDocumentCommand\presuperz{ m m m }{%
  {\mathop{}%
   \mathopen{\vphantom{#2}}^{#1}%
   \kern-3\scriptspace%
   #2\kern-8\scriptspace_{#3}%
}
}

\usepackage{amssymb}
\usepackage[numbers]{natbib}

\usepackage{pgfplots}
\usetikzlibrary{arrows}
\usetikzlibrary{calc}

\usepackage[english]{babel}
\usepackage[utf8x]{inputenc}
\usepackage[T1]{fontenc}

\usepackage{amsmath}
\usepackage{graphicx}
\usepackage{subfig}
\usepackage{bm} 

\makeatletter
\let\cl@chapter\undefined
\makeatletter

\usepackage[nameinlink]{cleveref}
\crefname{equation}{Eq.}{Eqs.}
\crefname{figure}{Fig.}{Figs.}
\crefname{section}{Section}{Sections}
\crefname{remark}{Remark}{Remarks}
\crefname{example}{Example}{Examples}
\crefname{appendix}{}{}

\hypersetup{allcolors=blue}

\newcommand{\uvec}[1]{\, \boldsymbol{\hat{\textbf{#1}}}}

\makeatletter
\newcommand{\oset}[3][0ex]{%
  \mathrel{\mathop{#3}\limits^{
    \vbox to#1{\kern-2\ex@
    \hbox{$\scriptstyle#2$}\vss}}}}
\makeatother

\begin{document}

\setlength{\abovedisplayskip}{6pt}
\setlength{\belowdisplayskip}{6pt}

\title{A unified, stable and accurate meshfree framework for peridynamic correspondence modeling. Part I: core methods}

\author{Masoud Behzadinasab \and
        Nathaniel Trask \and
        Yuri Bazilevs
}


\institute{M. Behzadinasab \at
              School of Engineering, Brown University, 184 Hope St., Providence, RI 02912, USA \\
              \email{masoud\_behzadinasab@brown.edu}           
           \and
           Nathaniel Trask \at
              Center for Computing Research, Sandia National Laboratories, Albuquerque, NM, USA
           \and
           Yuri Bazilevs \at
              School of Engineering, Brown University, 184 Hope St., Providence, RI 02912, USA\vspace{5pt}\\Sandia National Laboratories is a multi-mission laboratory managed and operated by National Technology and Engineering Solutions of Sandia, LLC., a wholly owned subsidiary of Honeywell International, Inc., for the U.S. Department of Energy’s National Nuclear Security Administration under contract DE-NA0003525. This paper describes objective technical results and analysis. Any subjective views or opinions that might be expressed in the paper do not necessarily represent the views of the U.S. Department of Energy or the United States Government.
}

\date{Received: date / Accepted: date}

\maketitle

\begin{abstract}
The overarching goal of this work is to develop an accurate, robust, and stable methodology for finite deformation modeling using strong-form peridynamics (PD) and the correspondence modeling framework. We adopt recently developed methods that make use of higher-order corrections to improve the computation of integrals in the correspondence formulation. A unified approach is presented that incorporates the reproducing kernel (RK) and generalized moving least square (GMLS) approximations in PD to obtain higher-order non-local gradients. We show, however, that the improved quadrature rule does not suffice to handle instability issues that have proven problematic for the correspondence-model based PD. In Part I of this paper, a bond-associative, higher-order core formulation is developed that naturally provides stability without introducing artificial stabilization parameters. Numerical examples are provided to study the convergence of RK-PD, GMLS-PD, and their bond-associated versions to a local counterpart, as the degree of non-locality (i.e., the horizon) approaches zero. Problems from linear elastostatics are utilized to verify the accuracy and stability of our approach. It is shown that the bond-associative approach improves the robustness of RK-PD and GMLS-PD formulations, which is essential for practical applications. The higher-order, bond-associated model can obtain second-order convergence for smooth problems and first-order convergence for problems involving field discontinuities, such as curvilinear free surfaces. In Part II of this paper we use our unified PD framework to: (a) study wave propagation phenomena, which have proven problematic for the state-based correspondence PD framework; (b) propose a new methodology to enforce natural boundary conditions in correspondence PD formulations, which should be particularly appealing to coupled problems. Our results indicate that bond-associative formulations accompanied by higher-order gradient correction provide the key ingredients to obtain the necessary accuracy, stability, and robustness characteristics needed for engineering-scale simulations.
\end{abstract}


\section{Introduction}
\label{sec:intro}

The classical (local) theory of solid mechanics involves the calculation of the divergence of the stress tensor to evaluate the internal state of the material. This approach involves differentiation of displacements, which becomes ill-defined with field discontinuities (e.g. cracks). The peridynamic theory (PD) \citep{silling2000reformulation,silling2007peridynamic} was developed as a non-local reformulation of the classical theory, replacing the divergence of the stress at a material point $\mbf{X}$ with an integral of the action of $\mbf{X}$ with other material points within a finite distance $\delta$, called the horizon. Peridynamics has been appealing to the computational mechanics community, given its natural capabilities in handling material failure without any complicated numerical treatment. The original PD model \citep{silling2000reformulation} is called the {\em bond-based} version as a bond force depends only upon its own deformation. Bond-based materials are limited to Poisson's ratio 1/3 and 1/4 in two and three dimensions, respectively. To address the limitation, the {\em state-based} PD framework \citep{silling2007peridynamic} was presented to incorporate the collective deformation of a neighborhood in deriving the bond actions.

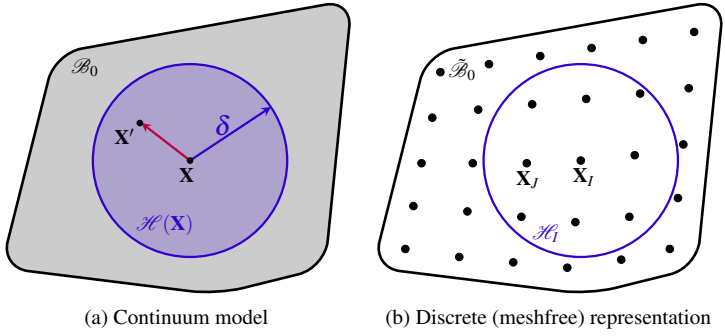
\begin{figure*}[!ht]
  \centering
  \subfloat[][Continuum model]{
\definecolor{wwqqzz}{rgb}{0.2,0.,0.8}
\resizebox {0.4\textwidth} {!} {
\begin{tikzpicture}[line cap=round,line join=round,>=triangle 45,x=1.0cm,y=1.0cm]

\path (0,0.5) coordinate (A) (4,0) coordinate (B) (6,0.5) coordinate (C)(6.5,5.5) coordinate (D) (1,4.5) coordinate (E);
\draw[rounded corners=5mm, line width=1.4pt,fill=black,fill opacity=0.2] (A) -- (B) -- (C) -- (D) -- (E) -- cycle;
\draw (1.2,4.5) node[anchor=north west] {\large $\mathcal{B}_0$};

\draw [line width=1.2pt,color=wwqqzz,fill=wwqqzz,fill opacity=0.2] (3.5,2.5) circle (1.8cm);
\draw (3.2,2.45) node[anchor=north west] {\large $\mathbf{X}$};
\draw [color=wwqqzz](2.4,1.6) node[anchor=north west] {\large $\mathcal{H}(\mathbf{X}) $};

\draw [->,>=stealth,line width=1.2pt,color=wwqqzz] (3.5,2.5) -- (5,3.5);
\draw [color=wwqqzz](3.85,3.45) node[anchor=north west] {\Large $\delta$};

\draw [->,>=stealth,line width=1.2pt,color=purple] (3.5,2.5) -- (2.6,3.2);

\draw (2,3.2) node[anchor=north west] {\large $\mathbf{X'}$};

\begin{scriptsize}
\draw [fill=black] (3.5,2.5) circle (1.5pt);
\draw [fill=black] (2.57,3.2) circle (1.5pt);
\end{scriptsize}

\end{tikzpicture}
}
}
  \subfloat[][Discrete (meshfree) representation]{
\definecolor{wwqqzz}{rgb}{0.2,0.,0.8}
\resizebox {0.4\textwidth} {!} {
\begin{tikzpicture}[line cap=round,line join=round,>=triangle 45,x=1.0cm,y=1.0cm]

\path (0,0.5) coordinate (A) (4,0) coordinate (B) (6,0.5) coordinate (C)(6.5,5.5) coordinate (D) (1,4.5) coordinate (E);
\draw[rounded corners=5mm, line width=1.4pt,fill=white,fill opacity=0.2] (A) -- (B) -- (C) -- (D) -- (E) -- cycle;
\draw (1.3,4.5) node[anchor=north west] {\large $\tilde{\mathcal{B}}_0$};

\draw [line width=1.2pt,color=wwqqzz,fill opacity=0.2] (3.85,2.5) circle (1.8cm);
\draw (3.6,2.45) node[anchor=north west] {\large $\mathbf{X}_I$};
\draw [color=wwqqzz](2.85,1.4) node[anchor=north west] {\large $\mathcal{H}_I$};


\draw (2.55,2.4) node[anchor=north west] {\large $\mathbf{X}_J$};

\begin{scriptsize}
\draw [fill=black] (0.6,0.85) circle (2pt);
\draw [fill=black] (1.6,0.7) circle (2pt);
\draw [fill=black] (2.6,0.6) circle (2pt);
\draw [fill=black] (3.6,0.5) circle (2pt);
\draw [fill=black] (4.6,0.65) circle (2pt);
\draw [fill=black] (5.5,0.85) circle (2pt);
\draw [fill=black] (5.65,1.8) circle (2pt);
\draw [fill=black] (5.75,2.8) circle (2pt);
\draw [fill=black] (5.85,3.85) circle (2pt);
\draw [fill=black] (5.95,4.9) circle (2pt);
\draw [fill=black] (5,4.75) circle (2pt);
\draw [fill=black] (4.05,4.6) circle (2pt);
\draw [fill=black] (3.1,4.45) circle (2pt);
\draw [fill=black] (2.15,4.3) circle (2pt);
\draw [fill=black] (1.25,4.15) circle (2pt);
\draw [fill=black] (1.1,3.3) circle (2pt);
\draw [fill=black] (0.9,2.45) circle (2pt);
\draw [fill=black] (0.75,1.65) circle (2pt);
\draw [fill=black] (1.75,1.55) circle (2pt);
\draw [fill=black] (2.75,1.45) circle (2pt);
\draw [fill=black] (3.75,1.35) circle (2pt);
\draw [fill=black] (4.75,1.45) circle (2pt);
\draw [fill=black] (4.85,2.6) circle (2pt);
\draw [fill=black] (4.95,3.75) circle (2pt);
\draw [fill=black] (3.95,3.65) circle (2pt);
\draw [fill=black] (2.95,3.55) circle (2pt);
\draw [fill=black] (1.95,3.45) circle (2pt);
\draw [fill=black] (1.85,2.45) circle (2pt);
\draw [fill=black] (2.85,2.45) circle (2pt);
\draw [fill=black] (3.85,2.5) circle (2pt);
\end{scriptsize}

\end{tikzpicture}
}
}
  \caption{Schematic of a peridynamic body in a reference configuration with horizon $\delta$.}
  \label{fig:PD} 
\end{figure*}

While Galerkin formulations can be developed for peridynamics (see, e.g. \citep{chen2011continuous,madenci2018weak}), the weak-form approach involves a six-dimensional (i.e., double) integral \citep{littlewood2015roadmap,bobaru2016handbook}, entailing considerable geometric complication and computational cost. Strong-form strategies are appealing, therefore, for practical applications. Meshfree techniques are commonly employed to solve the PD strong-form version using a nodal discretization of the PD equation of motion, in which the PD nodes act as both collocation and quadrature points \citep{silling2005meshfree}. \cref{fig:PD} shows a PD body in the continuum and discrete forms.

A special class of peridynamic materials, known as {\em correspondence models} \citep{silling2007peridynamic}, has been developed to incorporate classical constitutive equations (outcome of decades of experimental and theoretical research by the mechanics community) within the PD framework. In the correspondence formulation, a kinematic variable (e.g., the deformation gradient $\mbf{F}$) is computed using integration, rather than differentiation, and then used to evaluate its energy conjugate (e.g., the first Piola--Kirchhoff stress) through local material models. It was shown in \citep{bessa2014meshfree,hillman2020generalized} that there is a close connection between the PD correspondence model and meshfree discretizations of local theories, such as the reproducing kernel (RK) particle method \citep{liu1995reproducing,chen1996reproducing}. \citet{hillman2020generalized} showed the equivalence between the PD differential operator \citep{madenci2016peridynamic} and the RK implicit gradient operator \citep{chen2004implicit}, which both can be employed to obtain higher-order meshfree gradients, such as a higher-order $\mbf{F}$. Reproducing kernel enhanced approaches have been proposed to improve meshfree integration in peridynamic models \citep{pasetto2018reproducing,hillman2020generalized,leng2020asymptotically}.

Despite offering convenience in utilizing classical constitutive relations, PD correspondence materials have been found problematic in practice as they involve instability issues under non-uniform deformations. The unstable behavior appears as zero-energy mode oscillations in meshless peridynamic simulations and has been largely attributed to the cancellation of the non-uniform parts of deformation, which is allowed within the integration (averaging) technique \citep{breitenfeld2014non,tupek2014extended,silling2017stability,chowdhury2019modified,behzadinasab2020stability}. The instabilities can lead to large errors in practical applications \citep{kramer2019sandia}. To address the issue, stabilization techniques have been proposed, such as (1) addition of penalty terms to deviations from homogeneous deformations \citep{littlewood2010simulation,silling2017stability}, which involves artificial (non-physical) parameters, and (2) splitting of the integration domain into multiple sub-regions \citep{chowdhury2019modified}, which is unsettling as a continuum theory. A third group of methods called {\em bond-associative modeling} have have been developed \citep{breitzman2018bond,chen2018bond,behzadinasab2020semi}, which are continuum-based approaches and do not involve any numerical stabilization. The original correspondence model includes assigning the same kinematic variable ($\mbf{F}$) to all the bonds connected to a node, ignoring the non-uniform deformations in the neighborhood, which is the main root of the instability. The bond-associative modeling, on the other hand, increases the influence of each individual bond deformation on its own kinematic variable to take into effect inhomogeneous deformations. The bond-associated models can be seen as mixed bond-based, state-based formulation \citep[Section 5.4]{behzadinasab2020peridynamic}, incorporating strengths of both approaches and offering natural stability and robustness to the resulting formulation. 

In Part I of this two-part paper, we focus on a unified meshfree approach that considers non-local integration as a regularization of classical continuum mechanics to obtain higher-order local derivatives. The two methods of RK-PD \citep{hillman2020generalized} and generalized moving least squares (GMLS) \citep{trask2019asymptotically} serve our goal here. The RK and GMLS approximations were shown to have a close connection under a uniform-grid discretization \citep{leng2019asymptotically}. We consider convergence of these non-local models to a physically relevant local counterpart as the horizon approaches zero. This work disregards inherently non-local physical interactions, in which a finite horizon is attributed to a physically non-local phenomena (see e.g., \citep{bobaru2009convergence}). 

We analyze the robustness of the higher-order approaches of RK-PD and GMLS-PD. As we will demonstrate, these methods are unable to naturally provide stability; therefore, a bond-associated version of is proposed to naturally handle the instability issue without adding {\em tuning} parameters. We emphasize that without a stable methodology at hand, the applicability of RK-PD and GMLS-PD to practical simulations is questionable.

The remainder of Part I of this paper is as follows. A unified definition of the models is provided in \cref{sec:formulation}, after recalling first the PD correspondence model in \cref{subsec:PD}, RK-PD is summarized in \cref{subsec:RK-PD}, and a GMLS-PD version is developed in \cref{subsec:GMLS-PD}. We then present a bond-associative modeling approach in \cref{subsec:BA}. Numerical examples are provided in \cref{sec:numerics} to highlight the instability issue of the base models and study the robustness of the bond-associated versions. \cref{sec:conclusions} provides some concluding remarks. 

In Part II of this work, a framework for modeling natural boundary conditions in the strong-form peridynamics using the unified meshfree formulation is proposed, and numerical examples are provided. Dynamic wave propagation in the unified correspondence framework is also examined.

In what follows, all vectors are column vectors unless otherwise stated. Bold symbols indicate tensors of rank 1 (i.e., vectors) or higher.

\section{Model Formulation}
\label{sec:formulation}

The classical continuum theory of solid mechanics, in Lagrangian form, states that a static equilibrium state satisfies
\begin{equation}
  \nabla \cdot \mbf{P}(\mbf{X}) + \mbf{b}(\mbf{X}) = 0 , 
  \label{eqn:classical}
\end{equation}
where $\mbf{X}$ is a material point in the reference domain, $\mbf{u}$ is the displacement field, and $\mbf{b}$ is an applied body force density. $\nabla$ is the gradient operator, thus $\nabla \cdot$ is the divergence operator in the first term. $\mbf{P}$ is the first Piola--Kirchhoff stress, which is the energy conjugate of the deformation gradient, a kinematic variable based on gradient of displacement, i.e., 
$$\mbf{F} = \mbf{I} + \nabla \mbf{u} , $$
where $\mbf{I}$ is the Identity Matrix. This property is evaluated using integral formulations in peridynamics, which is covered in this section.

In the remainder of this section, we first recall correspondence material modeling in peridynamics, followed by a summary of the RK-PD formulation. We then develop a GMLS-PD model by applying the GMLS operator to the correspondence PD theory. Next, a bond-associated version of the two theories is presented, which provides a necessary tool for stability. 

\subsection{Peridynamic correspondence materials}
\label{subsec:PD}

As part of the state-based peridynamics, \citet{silling2007peridynamic} proposed a non-local deformation gradient as
\begin{gather}
  \bbar{F}(\mbf{X}) = \mbf{I} + \nabla \mbf{u}(\mbf{X}) , \notag \\
  \nabla \mbf{u}(\mbf{X}) = \int_{\mc{H}(\mbf{X})} \alpha(\mbf{X'}-\mbf{X}) \, \left[\mbf{u}(\mbf{X'})-\mbf{u}(\mbf{X})\right] \, \left[\mbf{X'}-\mbf{X}\right]^\intercal \, \bbar{K}^{-1}(\mbf{X}) \dX' ,
  \label{eqn:base-grad}
\end{gather}
where $\mc{H}(\mbf{X})$ is the {\em family} set (neighborhood) of $\mbf{X}$ defined in the reference body $\mc{B}_0$ as 
$$ \mc{H}(\mbf{X}) = \left\{ \mbf{X'} \ | \ \mbf{X'} \in \mc{B}_0, \, 0 < | \mbf{X'} - \mbf{X} | \leq \delta \right\} . $$
$\alpha(\mbf{X'}-\mbf{X})$ is a weighting function, called the {\em influence function}, which depends on the relative distances between material points with respect to the horizon. $\bbar{K}$ is called the shape tensor and defined by 
$$ \bbar{K}(\mbf{X}) = \int_{\mc{H}(X)} \alpha(\mbf{X'}-\mbf{X}) \, \left[\mbf{X'}-\mbf{X}\right] \, \left[\mbf{X'}-\mbf{X}\right]^\intercal \dX' .$$

In the state-based PD theory, {\em state} is a mathematical operator that maps a bond to a scalar, vector, or tensor-valued function. For example, the reference position vector state $\s{X}$ is defined as
$$ \s{X}(\mbf{X})\an{\mbf{X'}-\mbf{X}} = \mbf{X'} - \mbf{X} , $$
which maps the bond $(\mbf{X})\an{\mbf{X'}-\mbf{X}}$ to its reference image, the vector $[\mbf{X'} - \mbf{X}]$ originating from $\mbf{X}$.

\subsection{Reproducing kernel peridynamics}
\label{subsec:RK-PD}

The RK-PD formulation is based a non-local gradient operator, which uses the reproducing kernel shape functions to provide $n^\text{th}$-order accuracy in approximating local gradients. This non-local operator is constructed using a Taylor expansion procedure, based on a complete set of monomials, and incorporates the effects of higher-order terms in evaluating local gradients \citep[Section 4]{hillman2020generalized}. In this model, the deformation gradient is approximated by
\begin{align}
  \bbar{F}(\mbf{X}) 
   &= \mbf{I} + \int_{\mc{H}(\mbf{X})} \left[\mbf{u}(\mbf{X'})-\mbf{u}(\mbf{X})\right] \, \left[\s{\bs\Phi}(\mbf{X})\an{\mbf{X'}-\mbf{X}}\right]^\intercal \dX' ,
  \label{eqn:RK-F}
\end{align}
where $\s{\bs\Phi}$ is the RK-weight vector state 
$$\s{\bs\Phi}(\mbf{X})\an{\mbf{X'}-\mbf{X}} = \alpha(\mbf{X'}-\mbf{X}) \, \oset[.2ex]{\nabla}{\mbf{Q}}^\intercal \, \bbar{M}^{-1}_{[n]}(\mbf{X}) \, \mbf{Q}_{[n]}(\mbf{X'}-\mbf{X}) ,$$
in which $\mbf{Q}_{[n]}(\bxi)$ is a vector of the set of monomials $\{\bxi^{\, \beta}\}_{|\beta|=1}^{n}$. $\bbar{M}_{[n]}$ is called the moment matrix and defined as
$$\bbar{M}_{[n]}(\mbf{X}) = \int_{\mc{H}(\mbf{X})} \alpha(\mbf{X'}-\mbf{X}) \, \mbf{Q}_{[n]}(\mbf{X'}-\mbf{X}) \, \mbf{Q}_{[n]}^\intercal(\mbf{X'}-\mbf{X}) \dX' ,$$
and
$$\oset[.2ex]{\nabla}{\mbf{Q}} \; \equiv \mbf{Q}_{[n]}^{(\delta_{j1}, \, \delta_{j2}, \, \delta_{j3})} = [0 \ , \ \dots \ , \ 0 \ , \ \underset{jth \; {\rm entry}}{1} \ , \ 0 \ , \ \dots \ , \ 0]^\intercal.$$
$n$ is a free parameter, which determines the degree of accuracy of the non-local gradient operator. Note that the linear ($n=1$) RK-PD gradient operator is equivalent to the PD correspondence counterpart (cf. \cref{eqn:base-grad} with \cref{eqn:RK-F}).

\begin{example}
  \normalfont
  In three dimensions, for n=2 (quadratic accuracy) the involved $\mbf{Q}$-terms are 
  $$\mbf{Q}_{[2]}(\bxi) = [\xi_1, \, \xi_2, \, \xi_3, \, \xi_1^2, \, \xi_2^2, \, \xi_3^2, \, \xi_1\xi_2, \, \xi_1\xi_3, \, \xi_2\xi_3]^\intercal ,$$
  $$\oset[.2ex]{\nabla}{\mbf{Q}}_1 = [1, \, 0, \, 0, \, 0, \, 0, \, 0, \, 0, \, 0, \, 0]^\intercal ,$$
  $$\oset[.2ex]{\nabla}{\mbf{Q}}_2 = [0, \, 1, \, 0, \, 0, \, 0, \, 0, \, 0, \, 0, \, 0]^\intercal ,$$
  $$\oset[.2ex]{\nabla}{\mbf{Q}}_3 = [0, \, 0, \, 1, \, 0, \, 0, \, 0, \, 0, \, 0, \, 0]^\intercal ,$$
  $$\oset[.2ex]{\nabla}{\mbf{Q}} = [\oset[.2ex]{\nabla}{\mbf{Q}}_1, \oset[.2ex]{\nabla}{\mbf{Q}}_2, \oset[.2ex]{\nabla}{\mbf{Q}}_3] .$$
\end{example}

In the RK-PD model, the internal force is evaluated by approximating the divergence of the stress
\begin{align}
  \nabla \cdot \mbf{P}(\mbf{X}) 
  &\approx \nabla \cdot \bbar{P}(\mbf{X}) \notag \\
  &= \int_{\mc{H}(\mbf{X})} \left[\bbar{P}(\mbf{X'}) - \bbar{P}(\mbf{X})\right] \, \s{\bs\Phi}(\mbf{X})\an{\mbf{X'}-\mbf{X}} \dX' ,
  \label{eqn:RK-div}
\end{align}
where $\bbar{P}(\mbf{X}) = \mbf{P}(\bbar{F}(\mbf{X}))$. 

In the discrete form, \cref{eqn:RK-F,eqn:RK-div} are evaluated using 
\begin{align}
  \mbt{F}_I
  &= \mbf{I} + \nabla_h \mbf{u} \notag \\
  &= \mbf{I} + \sum_{J\in\mc{H}_I} \left[\mbf{u}_J - \mbf{u}_I\right] \, \bs\Phi_{IJ}^\intercal \, V_J ,
  \label{eqn:RK-F-discrete}
\end{align}
and
\begin{equation}
  \left(\nabla_h \cdot \mbt{P}\right)_I = \sum_{J\in\mc{H}_I} \left[\mbt{P}_J - \mbt{P}_I\right] \, \bs\Phi_{IJ} \, V_J ,
  \label{eqn:RK-div-discrete}
\end{equation}
where $I$ and $J$ refer to discrete nodes, and the family of $I$ is defined as 
$$ \mc{H}_I = \left\{ J \ | \ J \in \tilde{\mc{B}}_0, \, 0 < | \mbf{X}_J - \mbf{X}_I | \leq \delta \right\} . $$
The discrete RK-weight vector state is given as
$$\bs\Phi_{IJ} = \alpha_{IJ} \, \oset[.2ex]{\nabla}{\mbf{Q}}^\intercal \, \mbt{M}^{-1}_I \, \mbf{Q}_{IJ} ,$$
$$\mbt{M}_I = \sum_{J\in\mc{H}_I} \alpha_{IJ} \, \mbf{Q}_{IJ} \, \mbf{Q}_{IJ}^\intercal \, V_J .$$
Note that the subscript $[n]$ is implied in $\mbt{M}$ and $\mbf{Q}$ (dropped for notational simplicity).

\subsection{Generalized moving least square peridynamics}
\label{subsec:GMLS-PD}

To improve meshfree discretizations of peridynamic models and obtain higher-order strong-form formulations, \citet{trask2019asymptotically} developed a quadrature rule for non-local integro-differential equations, which has close ties to the approximation theory known as generalized moving least squares. It was shown that the quadrature rule can improve a PD elastic material model \citep{trask2019elastic}. Inspired by their framework, the GMLS algorithm is applied to \cref{eqn:base-grad} to construct a higher-order, non-local gradient operator at the discrete level:
\begin{equation}
  \left(\nabla_h \mbf{u}\right)_I = \sum_{J\in\mc{H}_I} \alpha_{IJ} \, \left[\mbf{u}_J-\mbf{u}_I\right] \, \left[\mbf{X}_J-\mbf{X}_I\right]^\intercal \bs\omega_{IJ} ,
  \label{eqn:GM-grad-0}
\end{equation}
in which $\bs\omega_{IJ}$ is a set of quadrature weights associated with the neighborhood of $I$ and are obtained by solving a local constrained optimization problem. While the original GMLS method was designed to improve the quadrature rule for the divergence operator, here we adopt it to improve integration for a non-local gradient operator at the discrete level. In this approach, we associate a quadrature weight per dimension to each bond and define it as a {\em diagonal} matrix (we define it as a matrix for the purpose of linear algebraic operations), i.e., in three dimensions, $\bs\omega_{IJ}$ is a diagonal matrix with size 3 by 3 for each bond. Note that it may still be possible to consider scalar weights, but a considerably larger size of neighbors is required then to ensure having enough unknowns for fulfilling the equality constraints in the optimization problem (for obtaining polynomial unisolvency).

For the purposes of this work, we choose the following as the influence function in \cref{eqn:GM-grad-0}:
\begin{equation}
  \alpha(\bxi) = \frac{1}{|\bxi|^2} ,
  \label{eqn:GM-inf}
\end{equation}
resulting in
\begin{equation}
  \left(\nabla_h \mbf{u}\right)_I = \sum_{J\in\mc{H}_I} \left[\mbf{u}_J-\mbf{u}_I\right] \, \frac{\left[\mbf{X}_J-\mbf{X}_I\right]^\intercal}{|\mbf{X}_J-\mbf{X}_I|^2} \bs\omega_{IJ} ,
  \label{eqn:GM-grad}
\end{equation}
and a unit-less quadrature weight (the following procedure should be straightforward to modify for other influence functions). In fact, $\left[\mbf{u}_J-\mbf{u}_I\right] \, \dfrac{\left[\mbf{X}_J-\mbf{X}_I\right]^\intercal}{|\mbf{X}_J-\mbf{X}_I|^2}$ can be seen as a bond-based gradient operator; thus, $\left(\nabla_h \mbf{u}\right)_I$ is a weighted average of the bond-based quantity over $\mc{H}_I$.

The quadrature weights are obtained by solving the following equality-constrained least squares problem 
\begin{align}
  & \underset{\bs\omega_{IJ}}{\rm argmin} \sum_{J\in\mc{H}_I} [ \bs\omega_{IJ}:\bs\omega_{IJ} ] \notag \\
  & {\rm such \ that,} ~~~~ \left(\nabla_h [p]\right)_I = \left(\nabla [p]\right)_I ~~~ \forall p \in \mbf{V}_h , 
  \label{eqn:optimization}
\end{align}
where $\mbf{V}_h$ denotes a Banach space of functions whose gradients should be computed exactly. While $\mbf{V}_h$ can be chosen for enforcing reproduction on any desired function, selecting it as the space of $n^\text{th}$-order polynomials would be beneficial for obtaining higher-order gradients. For the case of polynomials, it can be shown that $\Big[ \dfrac{(n+d)!}{n!d!}-1 \Big]$ number of non-colinear (in 2D) or non-coplanar (in 3D) bonds are required to satisfy the equality constraints \citep{liu1997moving}, where $d$ is the size of dimensions. For uniform discretizations, the horizon size must be chosen greater than $n \times h$ with $h$ being the nodal spacing.

\begin{example}
  \normalfont
  In order to obtain a 2nd-order deformation gradient, the non-local gradient operator must be able to exactly compute gradients of linear and quadratic functions. For establishing the weights $\bs\omega_{IJ}$, a local coordinate system is designed with its origin at $I$. The following equations must be used as constraints in the optimization problem, i.e. \cref{eqn:optimization}:
  \begin{itemize}
    \item $ p=x \ : \quad \left( \nabla [x] \right) |_{(0,0,0)} = 1 \uvec{i} + 0 \uvec{j} + 0 \uvec{k} ,$
    \item $ p=y \ : \quad \left( \nabla [y] \right) |_{(0,0,0)} = 0 \uvec{i} + 1 \uvec{j} + 0 \uvec{k} ,$
    \item $ p=z \ : \quad \left( \nabla [z] \right) |_{(0,0,0)} = 0 \uvec{i} + 0 \uvec{j} + 1 \uvec{k} ,$
    \item $ p=x^2 \ : \quad \left( \nabla [x^2] \right) |_{(0,0,0)} = 0 \uvec{i} + 0 \uvec{j} + 0 \uvec{k} ,$
    \item $ p=y^2 \ : \quad \left( \nabla [y^2] \right) |_{(0,0,0)} = 0 \uvec{i} + 0 \uvec{j} + 0 \uvec{k} ,$
    \item $ p=z^2 \ : \quad \left( \nabla [z^2] \right) |_{(0,0,0)} = 0 \uvec{i} + 0 \uvec{j} + 0 \uvec{k} ,$
    \item $ p=xy \ : \quad \left( \nabla [xy] \right) |_{(0,0,0)} = 0 \uvec{i} + 0 \uvec{j} + 0 \uvec{k} ,$
    \item $ p=xz \ : \quad \left( \nabla [xz] \right) |_{(0,0,0)} = 0 \uvec{i} + 0 \uvec{j} + 0 \uvec{k} ,$
    \item $ p=yz \ : \quad \left( \nabla [yz] \right) |_{(0,0,0)} = 0 \uvec{i} + 0 \uvec{j} + 0 \uvec{k} .$
  \end{itemize}
  For each coordinate component in the above system, $\left( \nabla_h [p] \right)_I = \left( \nabla [p] \right)_I$ , resulting in 27 total equations. Considering there are 3 weights per neighbor, at least 9 non-coplanar bonds are required to satisfy these constraints. Since the discrete gradient $\nabla_h$ is a linear operator, if it satisfies the above equations, it would satisfy all possible combinations of them as well; therefore, it will compute the gradient of any quadratic (and linear) function {\em exactly}.
\end{example}

Upon computing the quadrature weights, the discrete-level deformation gradient $\mbt{F}$ is computed using \cref{eqn:GM-grad}
\begin{align}
  \mbt{F}_I
  &= \mbf{I} + \sum_{J\in\mc{H}_I} \left[\mbf{u}_J - \mbf{u}_I\right] \, \frac{\left[\mbf{X}_J-\mbf{X}_I\right]^\intercal}{\left|\mbf{X}_J-\mbf{X}_I\right|^2} \bs\omega_{IJ} .
  \label{eqn:GM-F}
\end{align}

Since the divergence operator is the trace of the gradient, the discrete-level divergence of the stress can be computed using the same quadrature weights, i.e., 
\begin{align}
  \left(\nabla_h \cdot \mbt{P}\right)_I = \sum_{J\in\mc{H}_I} \left[\mbt{P}_J - \mbt{P}_I\right] \, \bs\omega_{IJ} \, \frac{\left[\mbf{X}_J-\mbf{X}_I\right]}{\left|\mbf{X}_J-\mbf{X}_I\right|^2} ,
  \label{eqn:GM-div}
\end{align}
noting that $\bs\omega_{IJ}^\intercal = \bs\omega_{IJ}$ is used. 

\begin{remark}
  \citet{leng2019asymptotically} shows connections between the RK and GMLS operators under uniform nodal discretization. It should however be noted that the RK shape functions are constructed on the basis of monomials and are more restricted compared to GMLS, which computes the quadrature weights completely based on an optimization framework. These two methods can generally result in different sets of weights, in particular for non-uniform nodal spacings; however, RK may always be expressed as a GMLS operator with particular choice of weighting. Our results in \cref{sec:numerics} show the different behavior between these two models. We will see that the particular choice of improved quadrature is less important than the stability gained via bond-associative modeling.
\end{remark}

\subsection{Bond-associative modeling}
\label{subsec:BA}

While the accuracy of the original PD gradient operator, i.e.~\cref{eqn:base-grad}, is enhanced by an improved integration rule through either the RK-PD or GMLS-PD formulation, they still suffer from instability issues (see \cref{sec:numerics}). A bond-associated formulation is developed here for RK-PD and GMLS-PD, which incorporates a naturally stabilizing technique without additional tunable parameters.

To employ the bond-associative corrections in the RK-PD and GMLS-PD formulations, we propose the following modification of the non-local divergence of stress 
\begin{align}
  \left(\nabla_h \cdot \mbt{P}\right)_I = \sum_{J\in\mc{H}_I} \left[\mbt{P}_{JI} - \mbt{P}_I\right] \, \bs\gamma_{IJ} ~ , \quad 
  \begin{cases}
    \bs\gamma_{IJ} = \bs\Phi_{IJ} \, V_J  ~~~ & \text{RK-PD}  \\
    \bs\gamma_{IJ} = \bs\omega_{IJ} \, \dfrac{\left[\mbf{X}_J-\mbf{X}_I\right]}{\left|\mbf{X}_J-\mbf{X}_I\right|^2}  ~~~ & \text{GMLS-PD} 
  \end{cases} ,
  \label{eqn:BA-div}
\end{align}
where $\mbt{P}_{JI}$ is called the {\em bond-associated} first Piola--Kirchhoff stress (a tensor state) and depends on the bond-associated deformation gradient (a tensor state), i.e., $\mbt{P}_{JI} = \mbf{P}(\mbt{F}_{JI})$, where $\mbt{F}_{JI}$, at the discrete level, is defined as
\begin{align}
  \mbt{F}_{JI} 
  &= \mbt{F}_J + \Delta \mbt{F}_{JI}^{\, nh} \notag \\
  &= \mbt{F}_J + \left[ \mbf{x}_J - \mbf{x}_I - \frac{\mbt{F}_I+\mbt{F}_J}{2}[\mbf{X}_J-\mbf{X}_I] \right] \frac{[\mbf{X}_J-\mbf{X}_I]^\intercal}{|\mbf{X}_J-\mbf{X}_I|^2} ,
  \label{eqn:BA-F}
\end{align}
where $\mbf{x}$ is the current position field, i.e., $\mbf{x} = \mbf{X} + \mbf{u}$. The new property $\Delta \mbt{F}_{JI}^{\, nh}$ is added to the bond-level deformation gradient, to account for non-uniform deformations. $\Delta \mbt{F}_{JI}^{\, nh} = \mbf{0}$ for the homogeneous case, i.e. if $\mbf{x} = \mbf{F} \, \mbf{X}$, $\mbt{F}_I = \mbt{F}_J = \mbf{F}$. Note that if the non-uniform part of deformation is neglected in this equation, i.e., $\Delta \mbt{F}_{IJ}^{\, nh} = \mbf{0}$, $\mbt{F}_{JI} = \mbf{F}_J$ and \cref{eqn:RK-div-discrete,eqn:GM-div} are recovered by \cref{eqn:BA-div}.

\begin{remark}
  Truncation error analysis of the bond-level deformation gradient shows that this property is limited to quadratic accuracy. Consider using the Taylor expansion process, applied on $I \rightarrow J$ and $J \rightarrow I$ subsequently as follows:
  \begin{subequations}
    \begin{align}
      & \mbf{x}_J = \mbf{x}_I + \mbf{F}_I \, [\mbf{X}_J-\mbf{X}_I] + \mbf{G}_I : \Big[ [\mbf{X}_J-\mbf{X}_I] [\mbf{X}_J-\mbf{X}_I]^\intercal \Big] + \mc{O}(h^3) , \label{eqn:Tij} \\
      & \mbf{x}_I = \mbf{x}_J + \mbf{F}_J \, [\mbf{X}_I-\mbf{X}_J] + \mbf{G}_J : \Big[ [\mbf{X}_I-\mbf{X}_J] [\mbf{X}_I-\mbf{X}_J]^\intercal \Big] + \mc{O}(h^3) , \label{eqn:Tji} 
    \end{align}
  \end{subequations}
  where $\mbf{G}$ is the 2nd-order spatial gradient of deformation, i.e. $\mbf{G} = \nabla \nabla \mbf{u} = \nabla \mbf{F}$, and $h= |\mbf{X}_J-\mbf{X}_I|$. Subtracting \cref{eqn:Tji} from \cref{eqn:Tij}, and dividing by 2, we obtain
  \begin{equation}
    \mbf{x}_J - \mbf{x}_I = \frac{\mbf{F}_I+\mbf{F}_J}{2} [\mbf{X}_J-\mbf{X}_I] + \frac{\mbf{G}_I-\mbf{G}_J}{2} : \Big[ [\mbf{X}_J-\mbf{X}_I] [\mbf{X}_J-\mbf{X}_I]^\intercal \Big] + \mc{O}(h^3) .
    \label{eqn:xij}
  \end{equation}
  Applying the Taylor series expansion, this time on $\mbf{G}$,
  \begin{equation*}
    \mbf{G}_I - \mbf{G}_J = \mc{O}(h) ,
  \end{equation*}
  thus,
  \begin{equation}
    [\mbf{G}_I - \mbf{G}_J] : \Big[ [\mbf{X}_J-\mbf{X}_I] [\mbf{X}_J-\mbf{X}_I]^\intercal \Big] = \mc{O}(h^3) .
    \label{eqn:Gij}
  \end{equation}
  Using \cref{eqn:xij,eqn:Gij},
  \begin{equation*}
    \mbf{x}_J - \mbf{x}_I - \frac{\mbf{F}_I+\mbf{F}_J}{2} [\mbf{X}_J-\mbf{X}_I] = \mc{O}(h^3) ,
  \end{equation*}
  \begin{equation*}
    \left[ \mbf{x}_J - \mbf{x}_I - \frac{\mbf{F}_I+\mbf{F}_J}{2} [\mbf{X}_J-\mbf{X}_I] \right] \frac{[\mbf{X}_J-\mbf{X}_I]^\intercal}{|\mbf{X}_J-\mbf{X}_I|^2} = \mc{O}(h^2) ,
  \end{equation*}
  and, at the discrete level,
  \begin{align}
    \mc{O}\left(\Delta \mbt{F}_{IJ}^{\, nh}\right) 
    &= \mc{O}\left( \, \left[ \mbf{x}_J - \mbf{x}_I - \frac{\mbt{F}_I+\mbt{F}_J}{2} [\mbf{X}_J-\mbf{X}_I] \right] \frac{[\mbf{X}_J-\mbf{X}_I]^\intercal}{|\mbf{X}_J-\mbf{X}_I|^2} \right) \notag \\
    &= \mc{O}(\mbt{F}_I) + \mc{O}(h^2) ,
    \label{eqn:Fnh-error}
  \end{align}
  noting that $\mc{O}(\mbt{F}_I) = \mc{O}(\mbt{F}_J) = \mc{O}(h^n)$. Using \cref{eqn:BA-F,eqn:Fnh-error},
  \begin{equation*}
    \mc{O}(\mbt{F}_{IJ}) = \mc{O}(h^p) \ , \quad p = \min(2,n).
  \end{equation*}
  Therefore, even if the accuracy of node-level deformation gradient is greater than quadratic ($n>2$), the bond-level property is of $\mc{O}(h^2)$. 
  \label{rmk:BA-F}
\end{remark}

\section{Numerical Examples}
\label{sec:numerics}

In this section, the reproducing kernel peridynamic method and its bond-associated version are denoted as RK-PD and BA-RK-PD, respectively. The generalized moving least square peridynamics and its bond-associated counterpart are called GMLS-PD and BA-GMLS-PD, respectively. The bond-associated model, without the higher-order corrections, is denoted as BA-PD, which is equivalent to the linear BA-RK-PD. The following numerical examples are provided for a threefold goal: (1) to test the convergence of the different variants of the correspondence model, (2) to demonstrate the unstable behavior of RK-PD and GMLS-PD and that it can be fixed through the bond-associative modeling approach, and (3) to study the robustness obtained by a combination of the bond-associative approach and higher-order corrections.

Linear, quadratic, and cubic models are considered in a set of two-dimensional (plane-strain) static problems. In each case, the horizon size $\delta$ is chosen according to the average nodal spacing and the order of the formulation ($n$) to ensure sufficient number of neighbors for obtaining polynomial unisolvency. Essential boundary conditions are enforced on a fictitious layer of size $\delta$, and the state of equilibrium is computed using a linear solver.

In the following examples, we consider the convergence of the above formulations to linear elastostatic solutions. Therefore, the horizon size and nodal spacing approach zero concurrently. Since the dependence of the non-local deformation gradients, i.e. \cref{eqn:RK-F-discrete,eqn:GM-F,eqn:BA-F}, on displacement is linear, the following constitutive model is constructed for consistency with linear elastostatics:
$$ \bs\eps = (\mbf{F}+\mbf{F}^T)/2 - \mbf{I} , $$
$$ \mbf{P} = \lambda \ {\rm tr}(\bs\eps) \, \mbf{I} + 2 \, \mu \, \bs\eps, $$
where $\lambda$ and $\mu$ are the Lam\'e parameters.

A cubic B-spline kernel (a smooth function) is utilized as the influence function for RK-PD (and its BA variant). This radial function is defined as:
\begin{align*}
  \alpha(\hat\xi) = 
  \begin{cases}
    \dfrac{2}{3} - 4\,\hat\xi^2 + 4\hat\xi^3  ~~~ & \text{for} ~~~ 0 < \hat\xi \leq \dfrac{1}{2} \\
    \dfrac{4}{3} - 4\,\hat\xi + 4\,\hat\xi^2 - \dfrac{4}{3}\hat\xi^3  ~~~ & \text{for} ~~~ \dfrac{1}{2} < \hat\xi \leq 1 \\
    0  ~~~ & \text{otherwise} 
  \end{cases}
  , 
\end{align*}
in which
$$\hat\xi \equiv \frac{|\mbf{X}_J-\mbf{X}_I|}{\delta} . $$
As noted previously, the influence function for GMLS-PD and BA-GMLS-PD is adopted as in \cref{eqn:GM-inf} and included in the corresponding formulations, i.e., \cref{eqn:GM-div,eqn:GM-F,eqn:BA-div}. 

To model the free surface behavior, fictitious nodes are placed in the free space region. In peridynamics language, all the bonds associated with the free-surface nodes are considered {\em broken}. In this framework, a broken bond carries a zero stress tensor that directly contributes to the force evaluation. A broken bond, however, does not engage in the kinematic variable computation; therefore, two sets of quadrature weights are required for integration. In other words, to obtain correct kinematic-gradient and stress-divergence operators, for material nodes with broken bonds (with free-surface neighbors here), two different set of quadrature weights are used: (1) a set of $\bs\Phi_{IJ}$ (for RK) and $\bs\omega_{IJ}$ (for GMLS) constructed using only the material points in the bulk of the body (including Dirichlet-bc nodes), for computing the kinematic variable (i.e., deformation gradient), and (2) a different set of $\bs\Phi_{IJ}$ and $\bs\omega_{IJ}$ obtained by considering the full neighborhood, for the internal force evaluation (i.e., divergence of stress). The fictitious free-surface nodes are only employed in the latter case. The broken bonds, associated with the free-surface nodes here, are given zero stress, i.e.
$$\mbt{P}_{JI} = \mbf{0} \  ~~~ \forall J \in \Omega_{\rm fs} ~ , ~ \forall I \in \mc{H}_J , $$
where $\Omega_{\rm fs}$ is the free-surface domain.

To clarify, the computed nodal kinematic variable must always pass a patch test, regardless of the boundary effects (e.g., free surface or damaged region). The divergence of the stress, however, should account for boundary conditions. Therefore, only neighbors from the bulk of material should contribute to the kinematic variable computation, while the full neighborhood should be involved in the internal force evaluation.


Error is calculated using the root-mean-square (RMS) norm
$$ ||e||_{2} = \sqrt{\frac{\sum_i^N e_i^2}{N}} ,$$
which is equivalent to the full $L_2$ norm for quasi-uniform pointsets \citep{wendland2004scattered}.

\subsection{Convergence to a manufactured solution}
\label{subsec:square}

The following 2D manufactured solution is considered over a square domain $[-1,1] \times [-1,1]$:
\begin{align*}
  & u_1(x,y) = A \sin(\pi x/2) \cos(\pi y/2) + B \exp(x) \exp(y) , \notag \\
  & u_2(x,y) = C \cos(\pi x/2) \sin(\pi y/2) + D \exp(x) \exp(y) ,
\end{align*}
with the corresponding body force $\mbf{b}$, which satisfies \cref{eqn:classical},
\begin{align*}
  b_1(x,y) = & - \frac{\pi^2}{4} \Big[ [A+C] \lambda + [3A+C] \mu \Big] \sin(\pi x/2) \cos(\pi y/2) \notag \\
      & + \Big[ [B+D] \lambda + [3B+D] \mu \Big] \exp(x) \exp(y) , \notag \\
  b_2(x,y) = & - \frac{\pi^2}{4} \Big[ [A+C] \lambda + [A+3C] \mu \Big] \cos(\pi x/2) \sin(\pi y/2) \notag \\
      & + \Big[ [B+D] \lambda + [B+3D] \mu \Big] \exp(x) \exp(y) .
\end{align*}
The constants are chosen as $A=0.2$, $B=-0.15$, $C=-0.15$, and $D=0.1$. A material with Young's modulus $E=100,000$ and Poisson's ratio $\nu=0.3$ is considered in this test.

\begin{figure*}[!ht]
  \centering
  \subfloat[][Base level (L0)]{\includegraphics[height=0.25\textwidth]{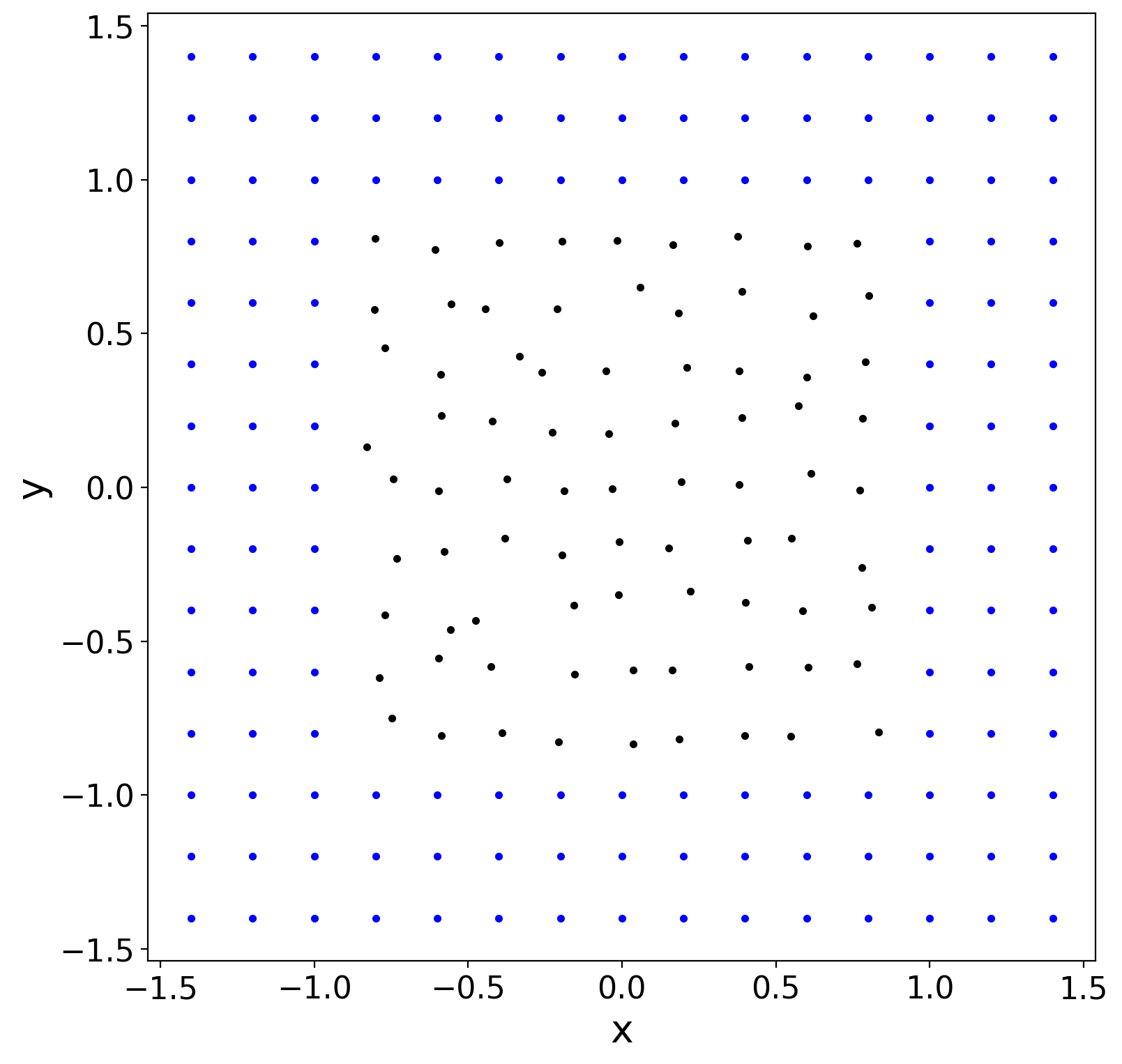}}
  \hspace*{0.1cm}
  \subfloat[][L1]{\includegraphics[height=0.25\textwidth,trim={2.5cm 0 0 0},clip]{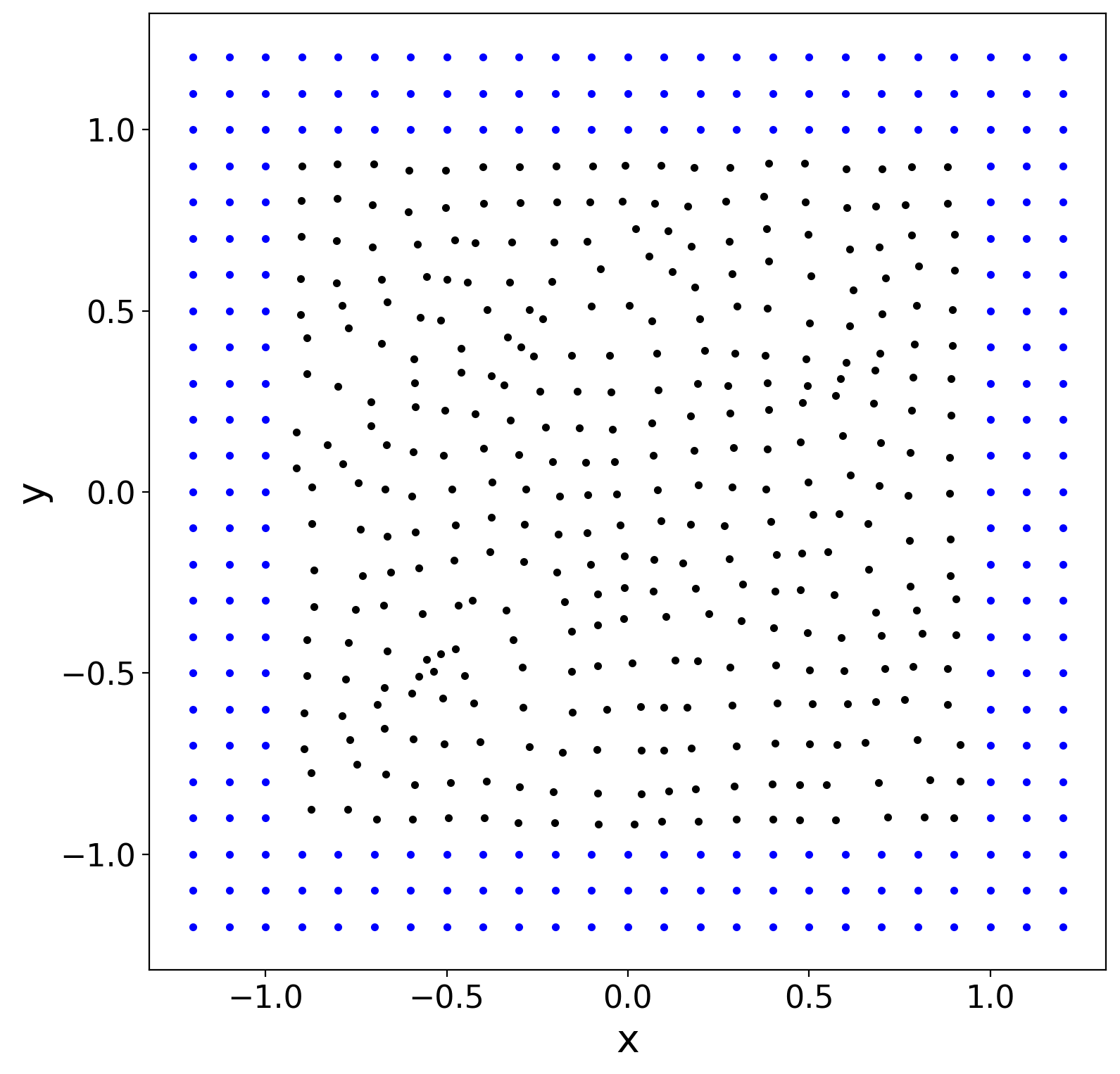}}
  \hspace*{0.1cm}
  \subfloat[][L2]{\includegraphics[height=0.25\textwidth,trim={2.5cm 0 0 0},clip]{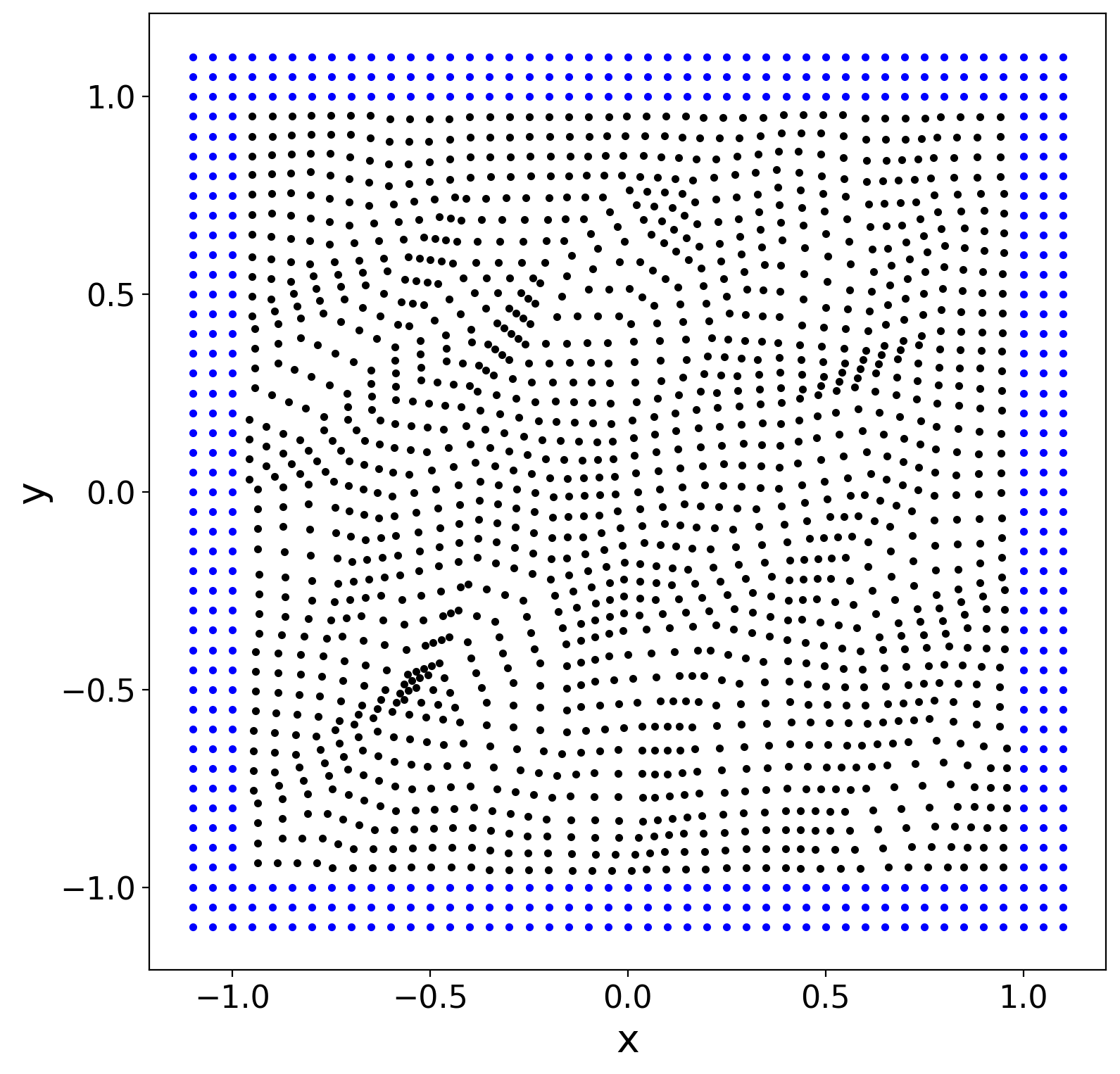}}
  \hspace*{0.1cm}
  \subfloat[][L3]{\includegraphics[height=0.25\textwidth,trim={2.5cm 0 0 0},clip]{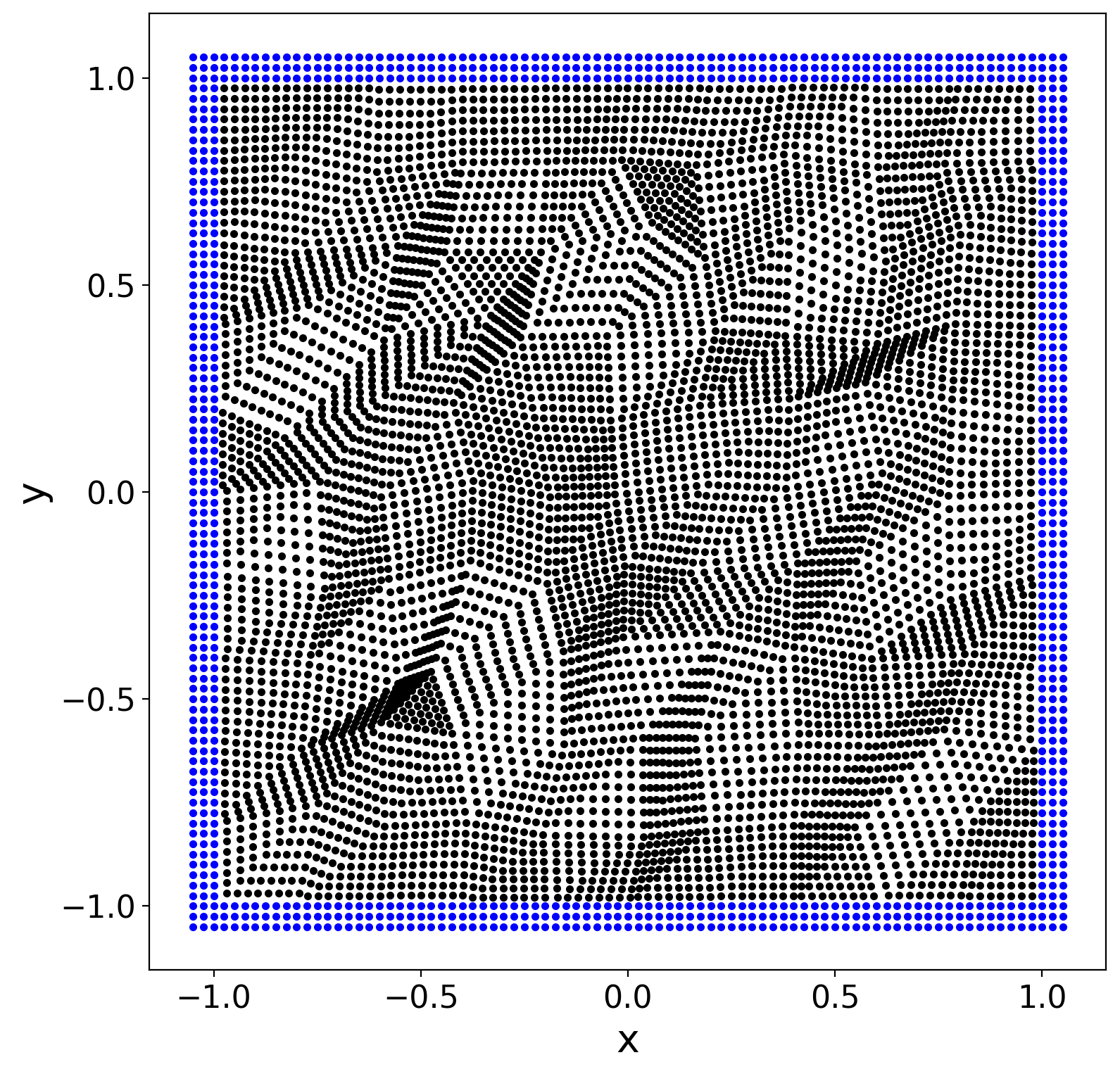}}
  \caption{Non-uniform nodal discretizations and refinements for the 2D manufactured problem. The blue color denotes the nodes used to prescribe displacement-controlled boundary conditions.}
  \label{fig:square-mesh}
\end{figure*}

This problem is solved on uniform and non-uniform discretizations, with average spacing $h=[0.2, \, 0.1, \, 0.05, \, 0.025]$. To obtain the non-uniform discretizations, the nodes in the base level ($h=0.2$) of the uniform case are perturbed with a random normal distribution with the standard deviation of 0.03 (15\% of the average nodal spacing), then refined in a systematic way, i.e., nodes are placed on midpoints during refinements. \cref{fig:square-mesh} shows the non-uniform nodal distributions for this problem. The horizon size for the linear, quadratic, and cubic models is chosen as $\delta = 2.5 h$, $\delta = 3.5 h$, and $\delta = 4.5 h$, respectively.

\begin{figure*}[!ht]
  \centering
  \subfloat[][Uniform discretization - Linear]{\includegraphics[height=0.39\textwidth]{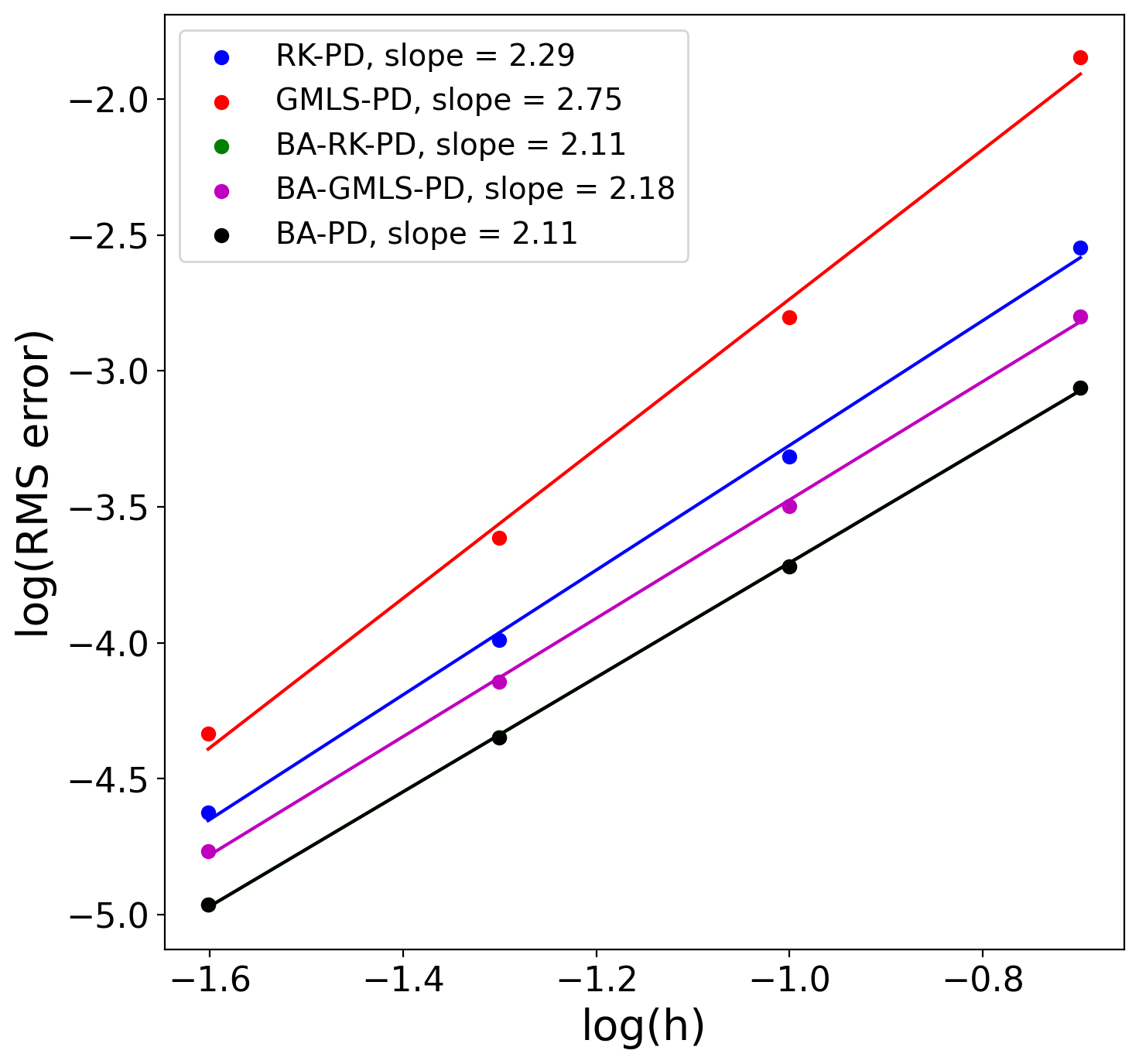}}
  \hspace*{1cm}
  \subfloat[][Non-uniform discretization - Linear]{\includegraphics[height=0.39\textwidth]{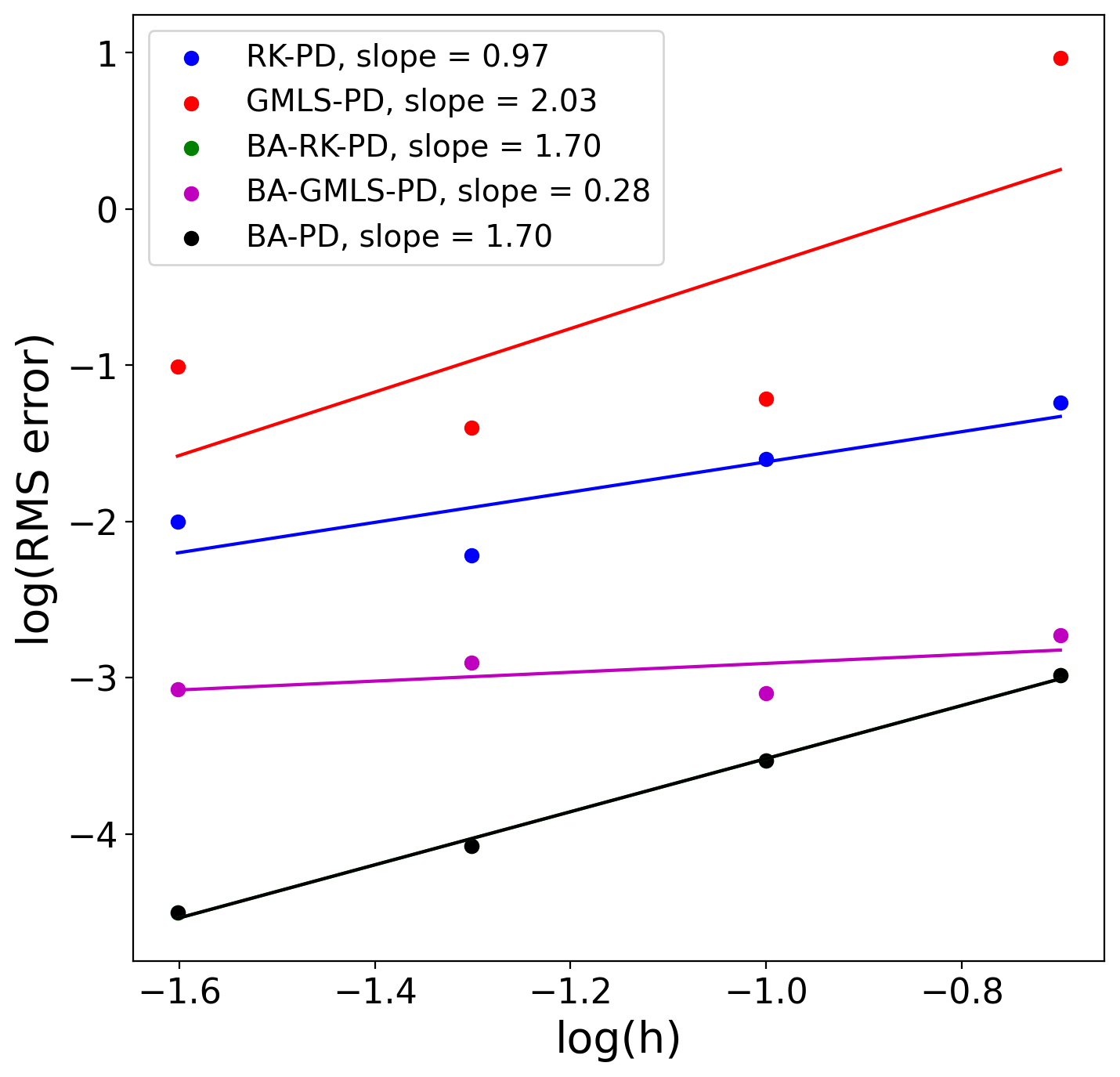}}

  \subfloat[][Uniform discretization - Quadratic]{\includegraphics[height=0.39\textwidth]{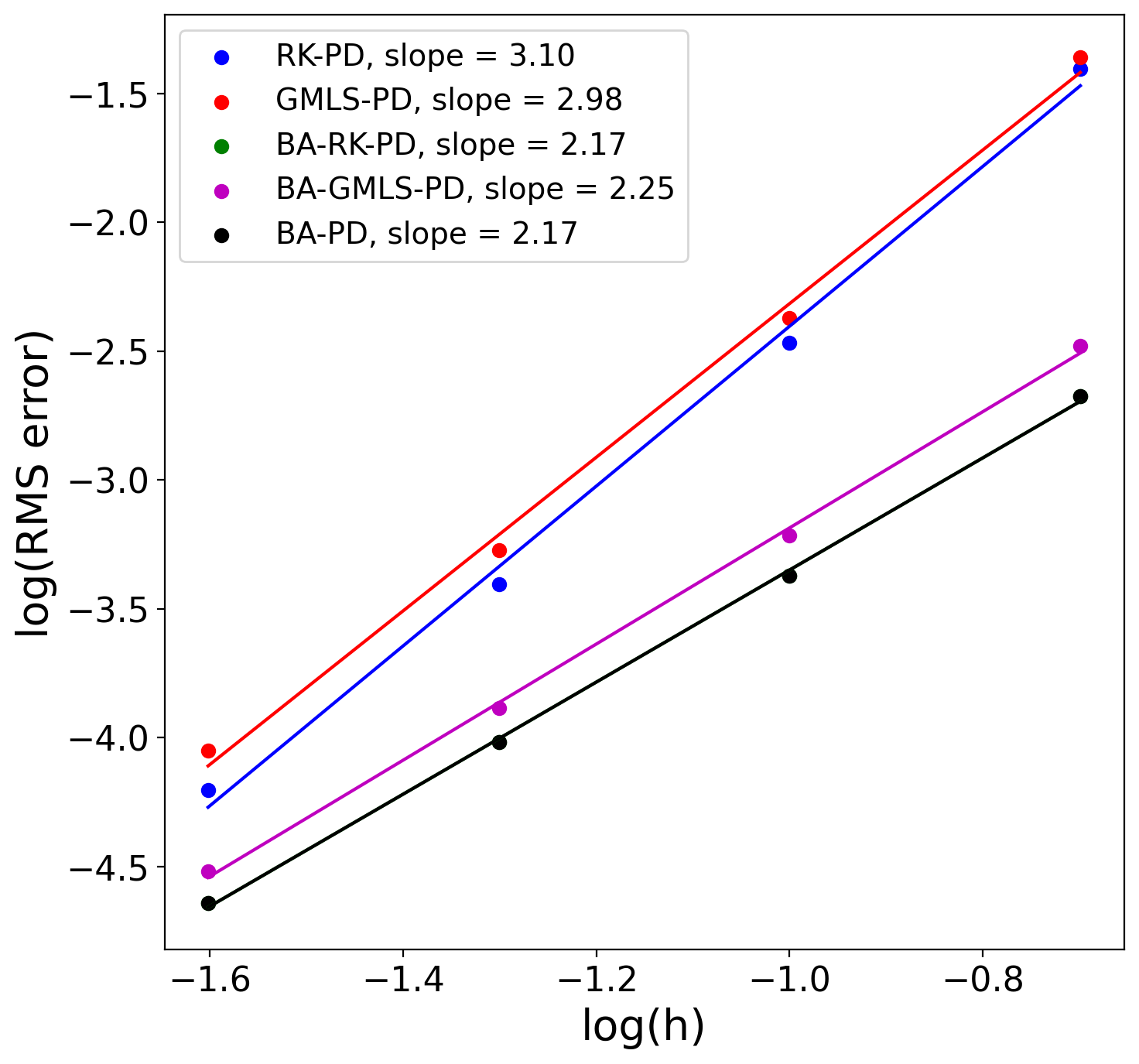}}
  \hspace*{1cm}
  \subfloat[][Non-uniform discretization - Quadratic]{\includegraphics[height=0.39\textwidth]{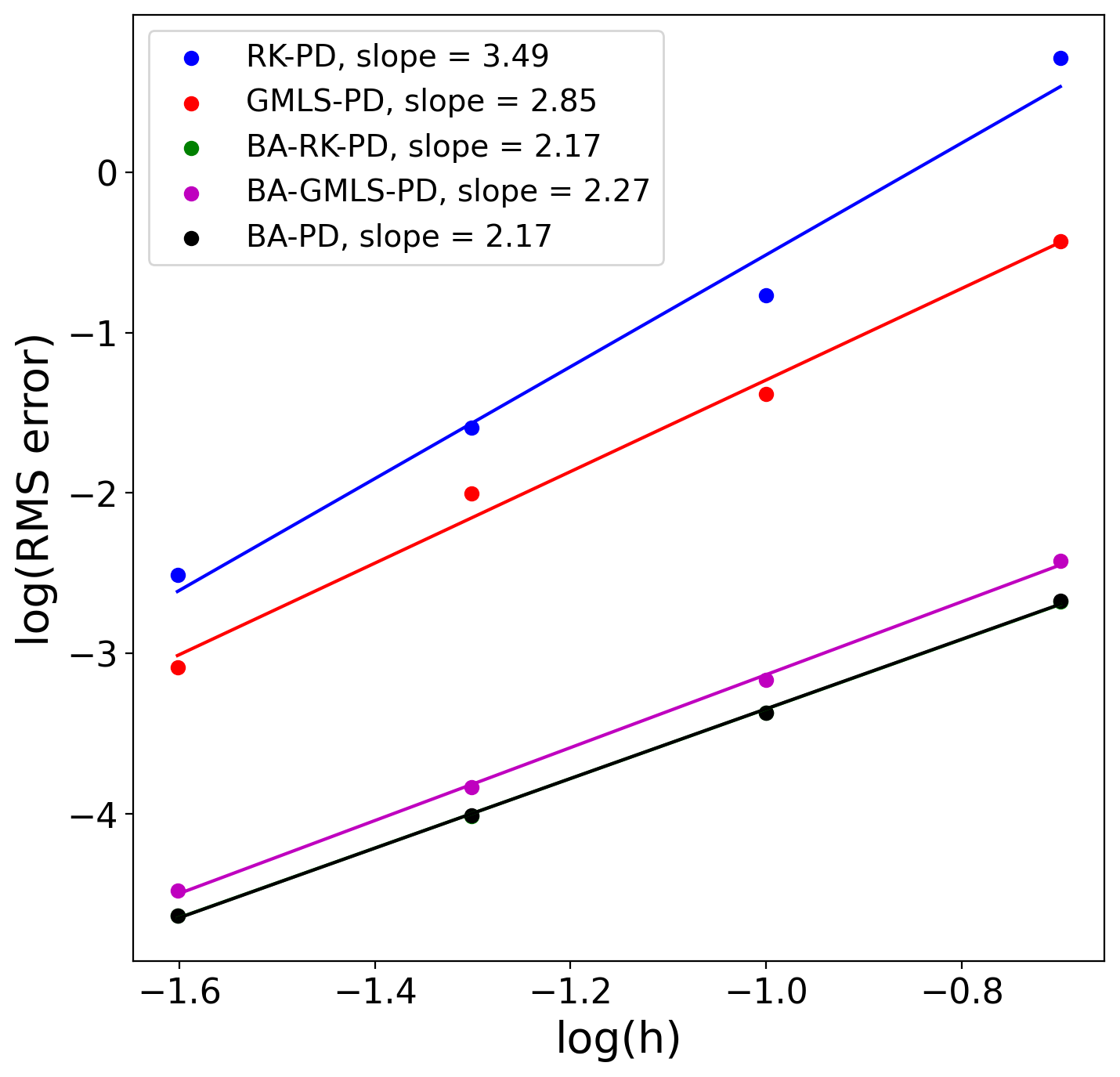}}

  \subfloat[][Uniform discretization - Cubic]{\includegraphics[height=0.39\textwidth]{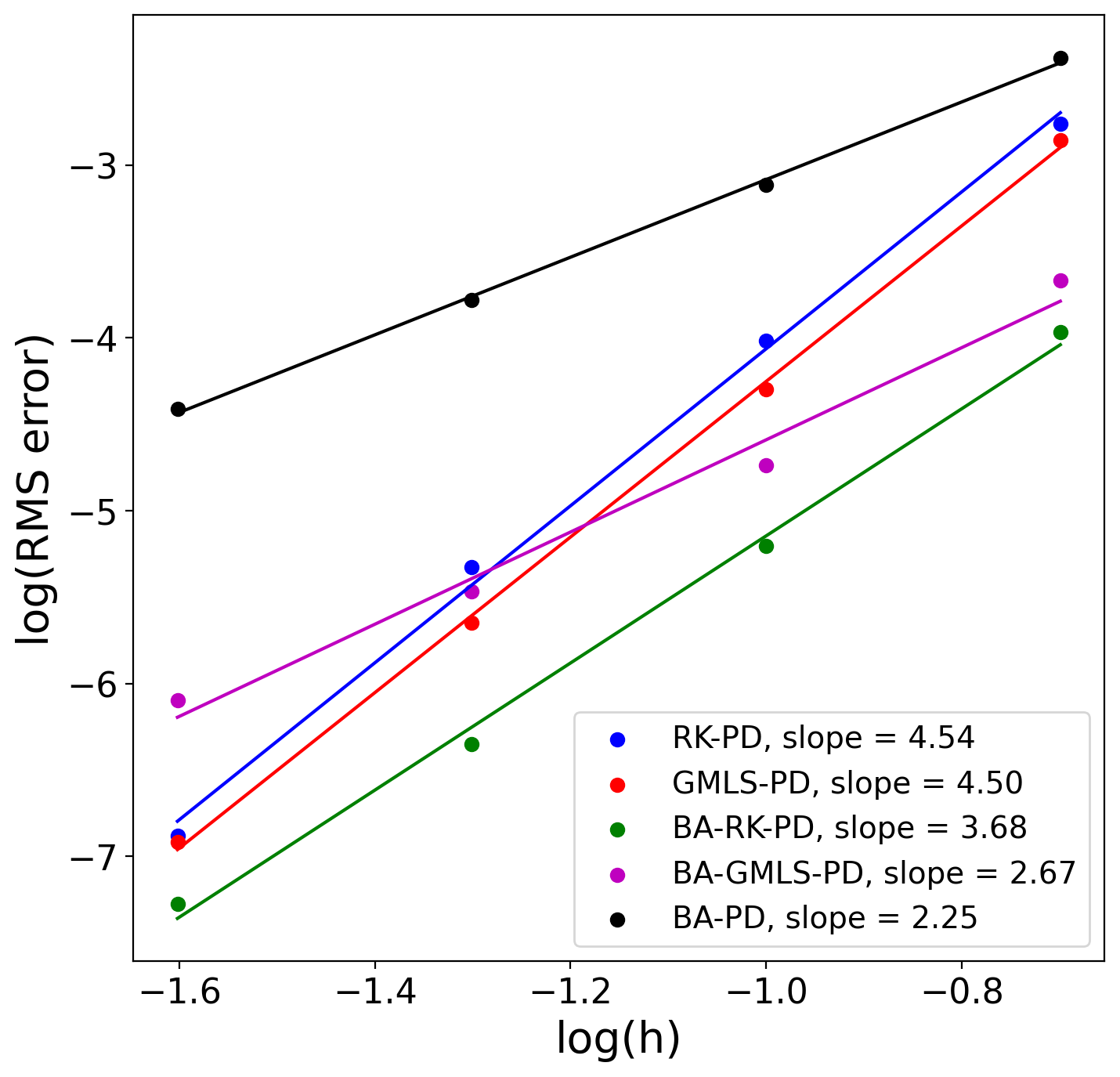}}
  \hspace*{1cm}
  \subfloat[][Non-uniform discretization - Cubic]{\includegraphics[height=0.39\textwidth]{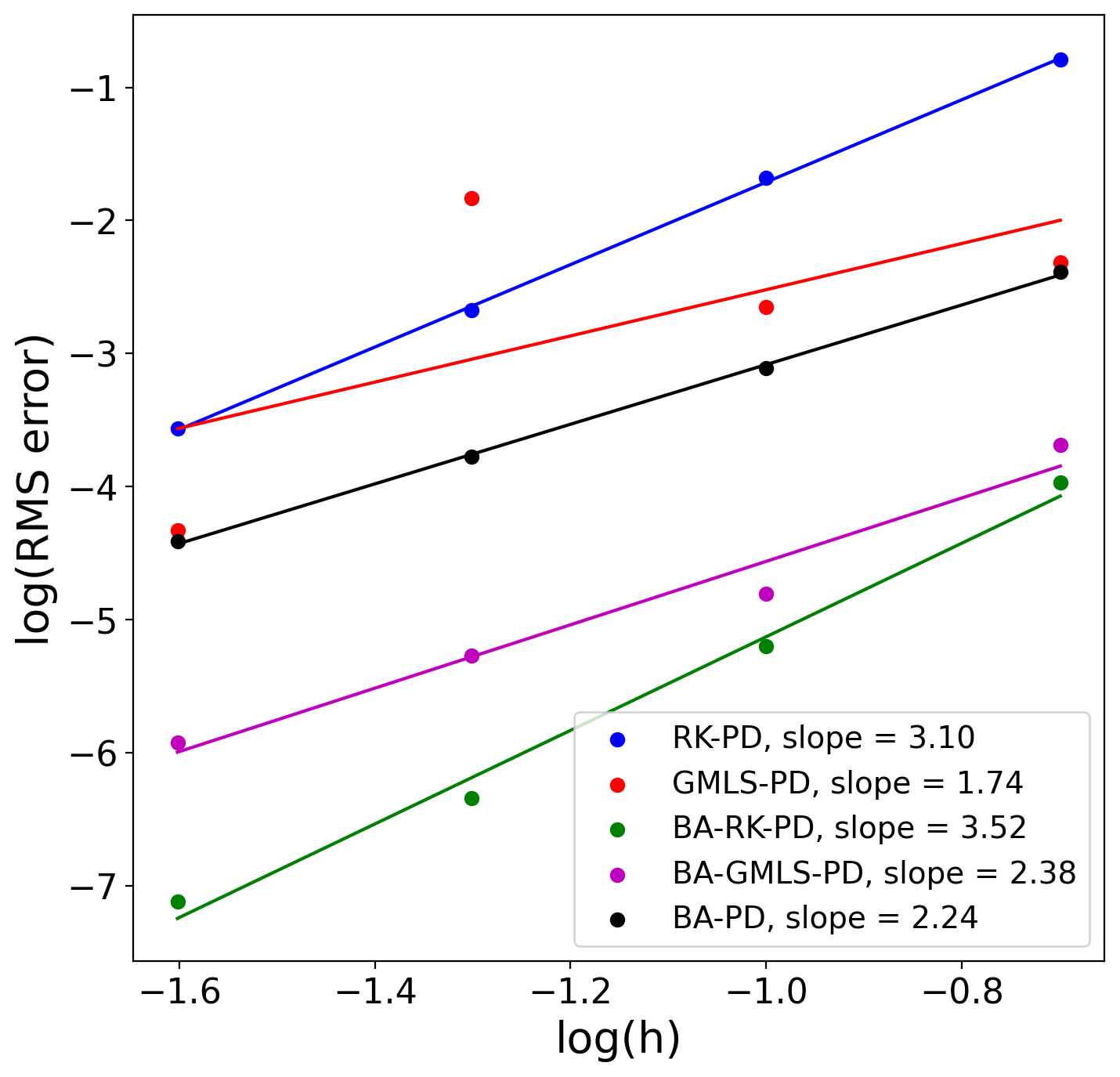}}
  \caption{Two-dimensional convergence test for RK-PD, GMLS-PD, BA-RK-PD, BA-GMLS-PD, and BA-PD on a manufactured solution. Linear, quadratic, and cubic formulations are tested on uniform and non-uniform discretizations. BA-PD and BA-RK-PD overlap in (a--d). Plots show the displacement RMS error values.}
  \label{fig:square-convergence}
\end{figure*}

\cref{fig:square-convergence} shows the RMS error in displacements. In this problem, generally speaking, RK-PD and GMLS-PD results are similar, with the former giving less error in the linear case and the latter being more accurate in the quadratic and cubic cases (on non-uniform discretizations). In all cases, the bond-associated, higher-order results are superior to the base models. BA-RK-PD and BA-GMLS-PD solutions are very close in the quadratic case and uniform, linear case. In the other three cases, BA-RK-PD shows a better convergence rate. We speculate that this might be due to the specific choices of influence function made in this work. Further study is required to shed light on the role of the influence function.

It should be noted that the linear and quadratic results are similar in the uniform case (and to some extent in the non-uniform case). {\em Super-convergence} of the solution in some cases is observed, i.e., observing the rate of $n+1$ for an $n^\text{th}$-order method, which agrees with the results of \citep{leng2019super,trask2019asymptotically,hillman2020generalized,leng2020asymptotically}. In this example, not only for the linear case, but also for the quadratic case (due to super-convergence), the BA-PD and BA-RK-PD approaches produce the same results. Considering the non-uniform case, the bond-associated, higher-order results are considerably improved by increasing the order from quadratic to cubic (that does not occur for BA-PD). On the other hand, the RK-PD and GMLS-PD results are barely improved in that case. That is in fact due to the presence of instabilities in the solution, which is not observed for the uniform case in this problem. 

\begin{figure*}[!ht]
  \centering
  \subfloat[][RK-PD]{\includegraphics[height=0.21\textwidth]{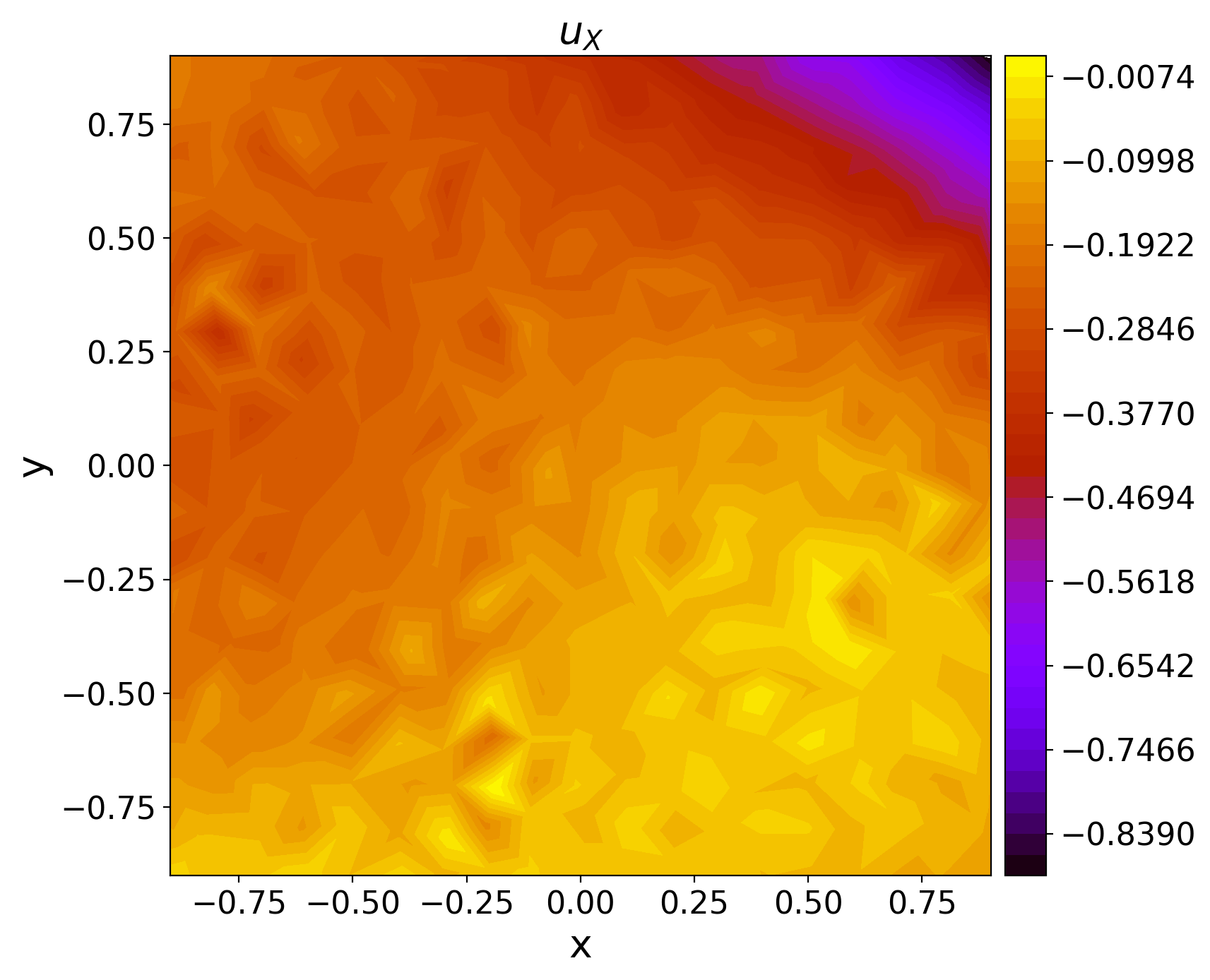}}
  \hspace*{0.2cm}
  \subfloat[][GMLS-PD]{\includegraphics[height=0.21\textwidth,trim={3cm 0 0 0},clip]{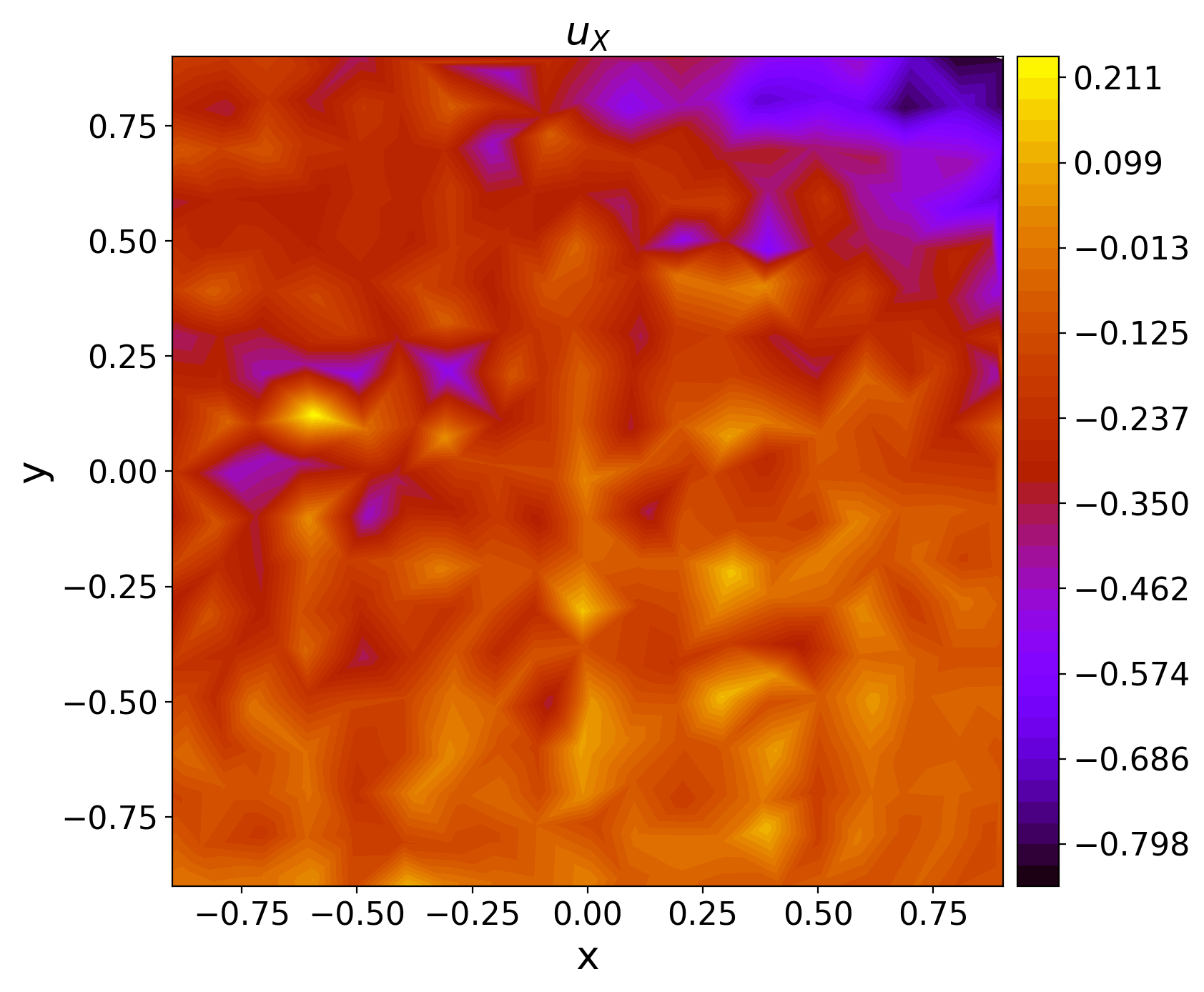}}
  \hspace*{0.2cm}
  \subfloat[][BA-RK-PD]{\includegraphics[height=0.21\textwidth,trim={3cm 0 0 0},clip]{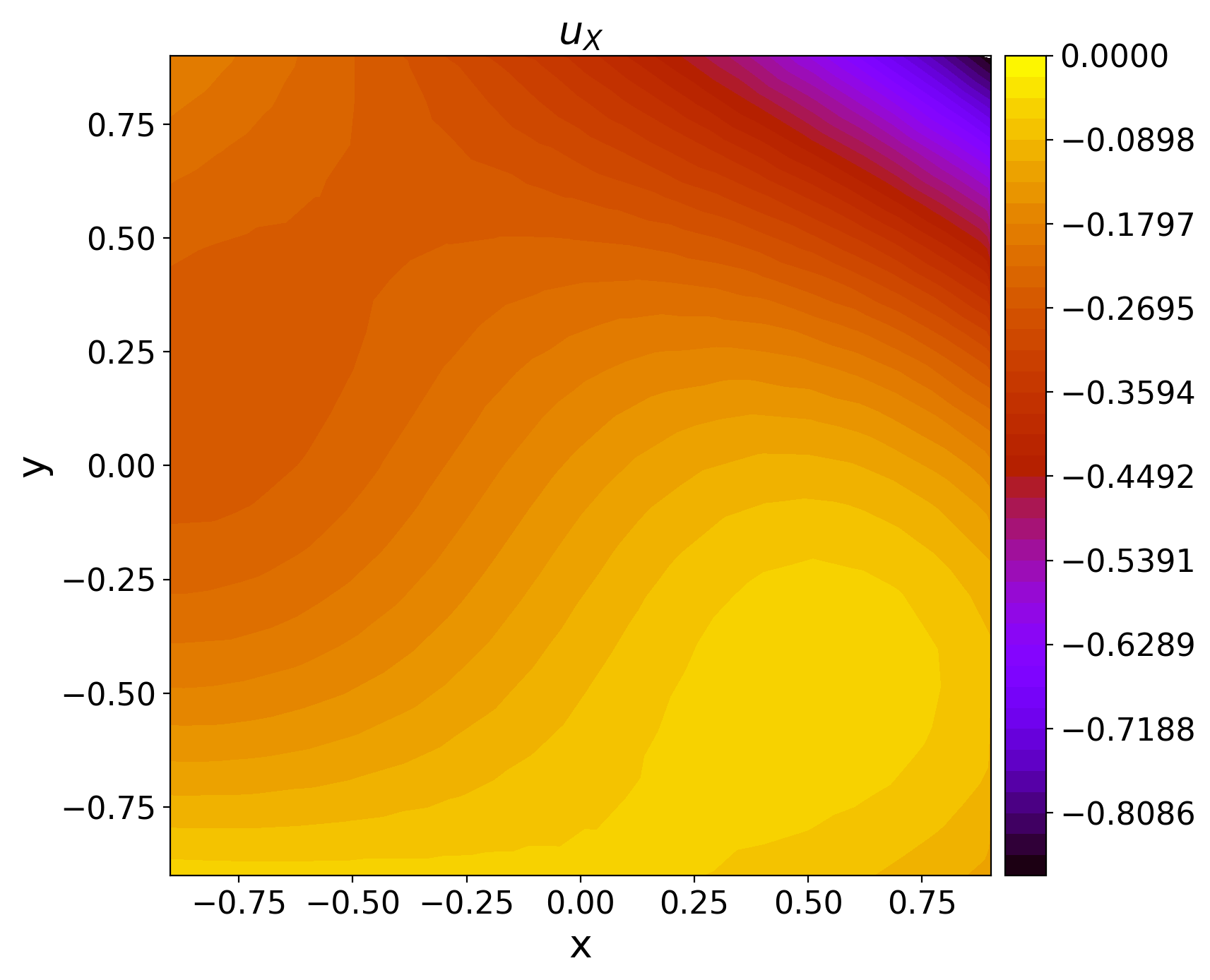}}
  \hspace*{0.2cm}
  \subfloat[][BA-GMLS-PD]{\includegraphics[height=0.21\textwidth,trim={3cm 0 0 0},clip]{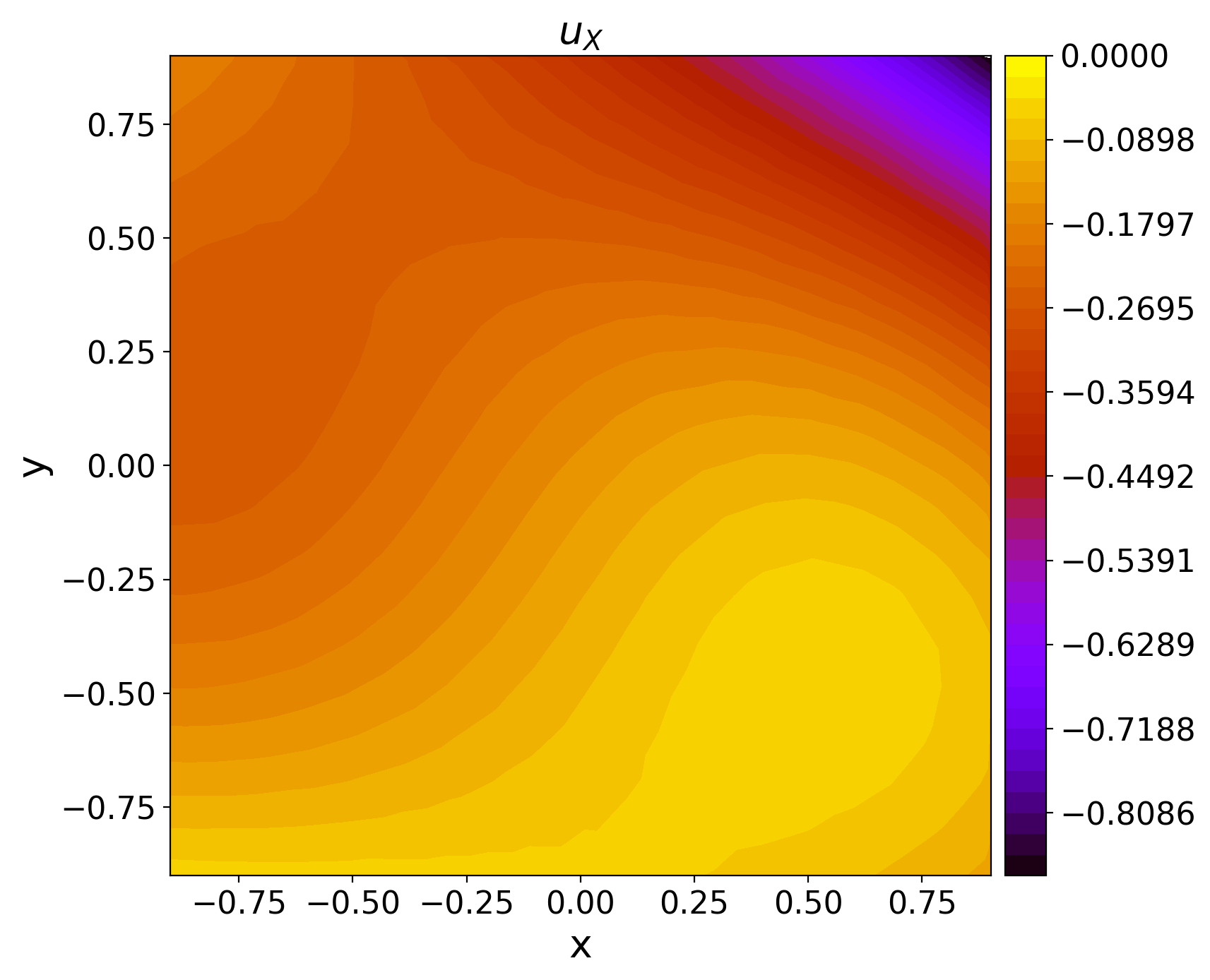}}

  \subfloat[][RK-PD]{\includegraphics[height=0.21\textwidth]{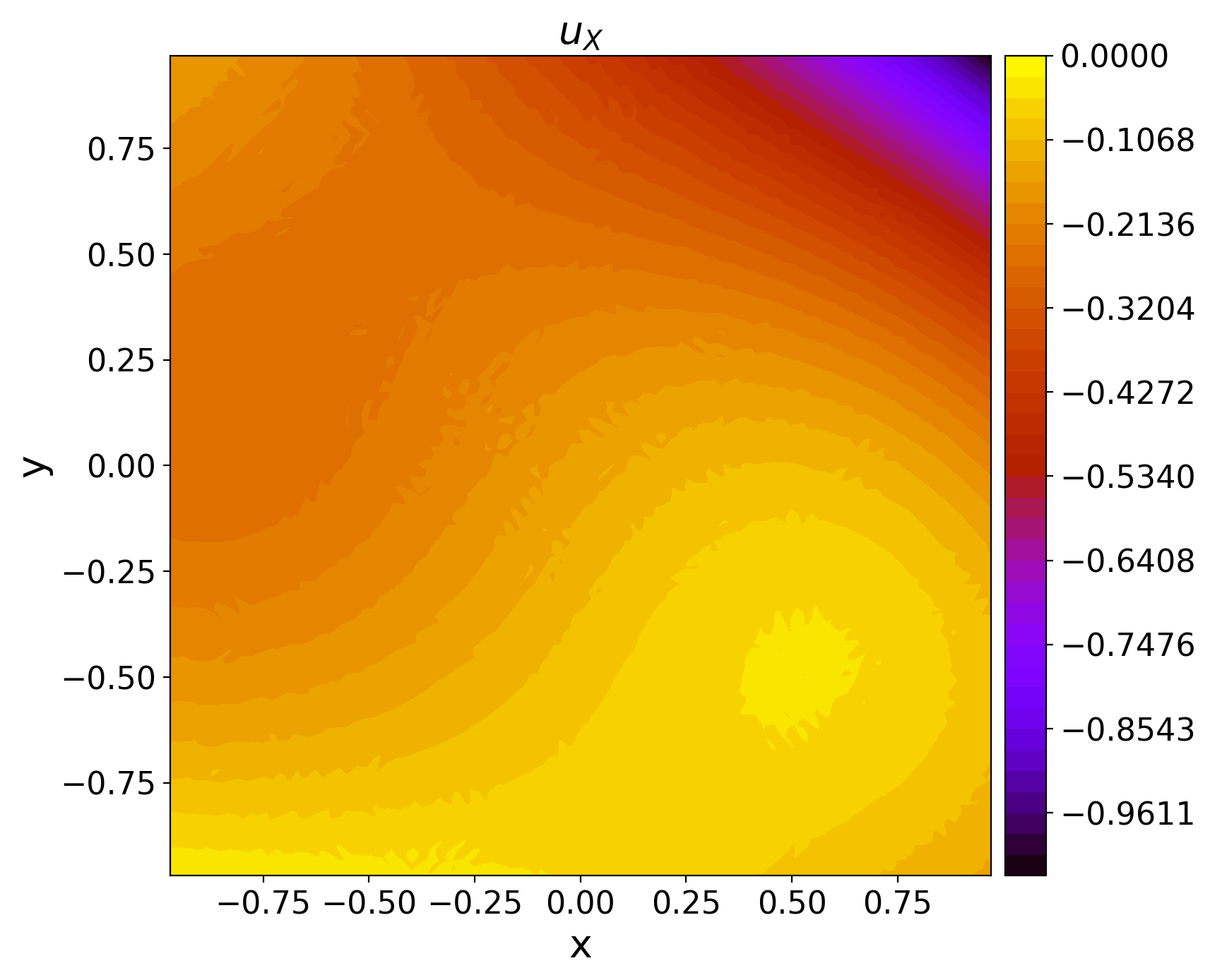}}
  \hspace*{0.2cm}
  \subfloat[][GMLS-PD]{\includegraphics[height=0.21\textwidth,trim={3cm 0 0 0},clip]{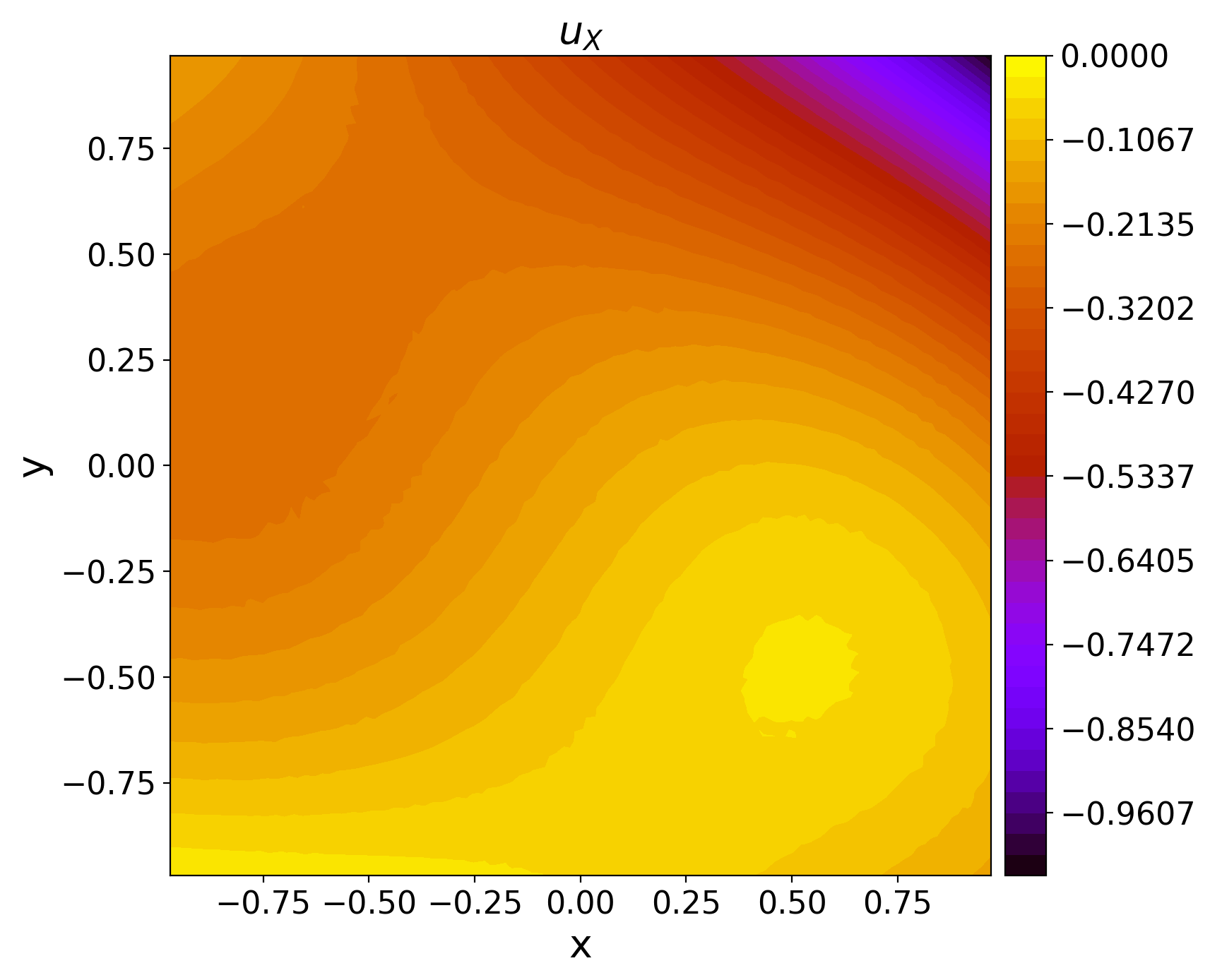}}
  \hspace*{0.2cm}
  \subfloat[][BA-RK-PD]{\includegraphics[height=0.21\textwidth,trim={3cm 0 0 0},clip]{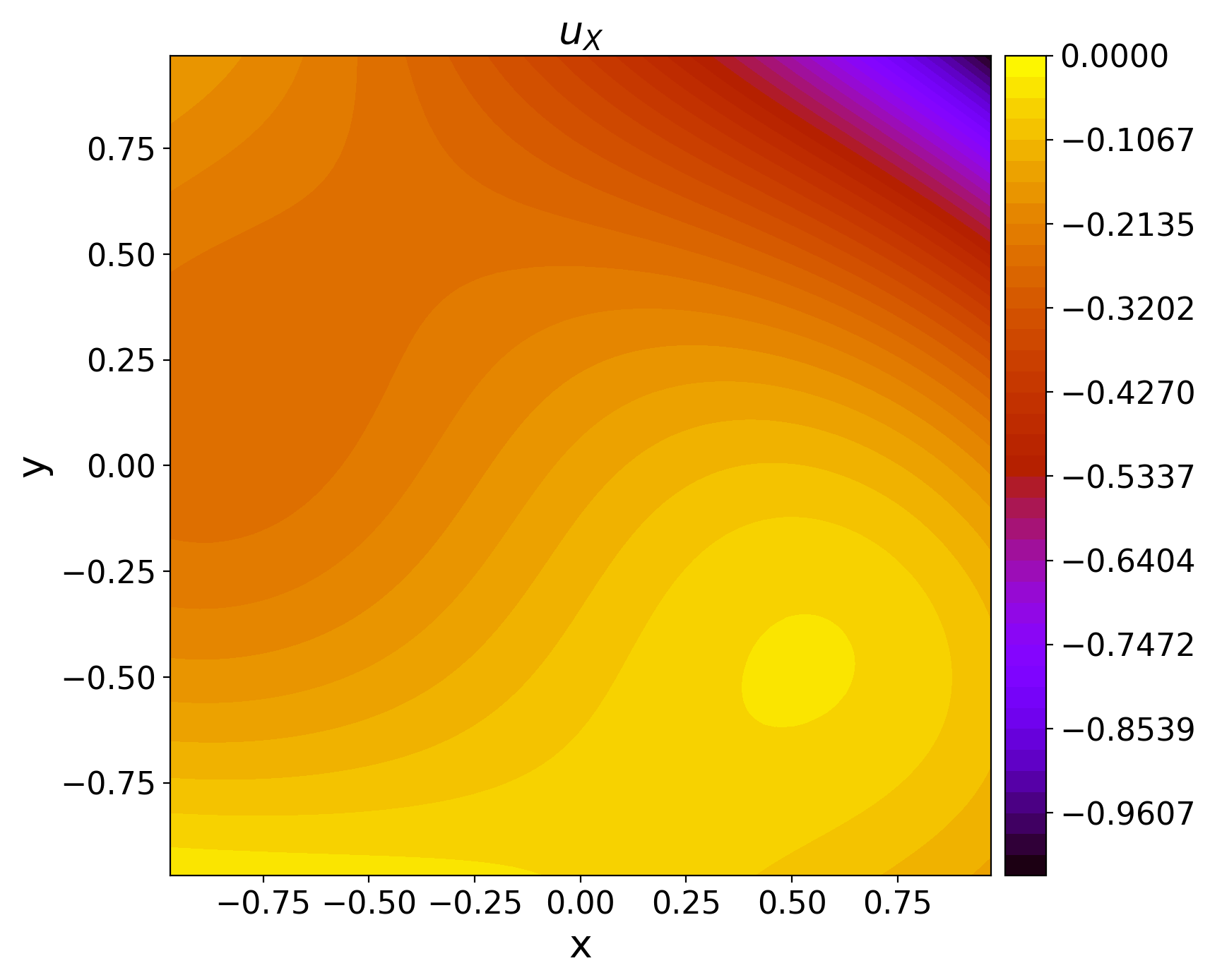}}
  \hspace*{0.2cm}
  \subfloat[][BA-GMLS-PD]{\includegraphics[height=0.21\textwidth,trim={3cm 0 0 0},clip]{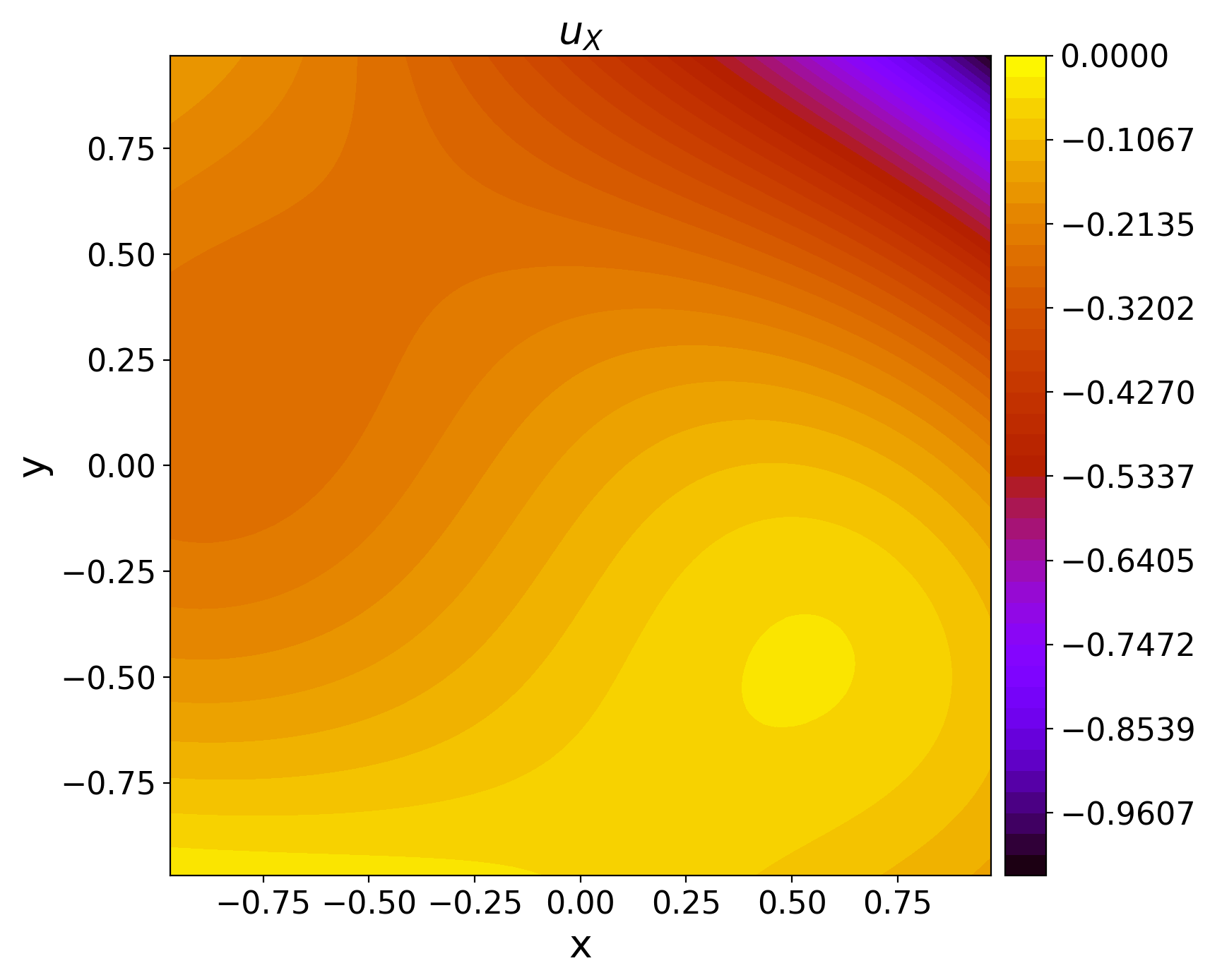}}
  \caption{Horizontal-displacement contours in the 2D manufactured problem using the quadratic formulations. Top and bottom rows correspond to the L1 and L3 non-uniform discretizations, respectively. While the oscillations in the RK-PD and GMLS-PD fields are reduced with nodal refinement, they are not completely eliminated. On the other hand, there is no oscillation in the bond-associated fields, even in the coarse grids. The refined bond-associated solutions (g-h) are in good agreement with the exact solution (not shown here).}
  \label{fig:square-contours}
\end{figure*}

The unstable behavior is demonstrated in \cref{fig:square-contours}, where the contour plots of the horizontal displacements, solved by the quadratic models, are provided for two levels of non-uniform grids (denoted as L1 and L3 in \cref{fig:square-mesh}). While the high degree of oscillations in the RK-PD and GMLS-PD results in the coarse level are reduced by refinement, the unstable behavior is still observed in the refined case. On the other hand, there is no instabilities in the bond-associated methods, and they provide meaningful solutions even on coarse meshes. This property is essential for practical applications of the methods.

\begin{figure*}[!ht]
  \centering
  \subfloat[][$\delta = 2.75 \, h$]{\includegraphics[height=0.28\textwidth]{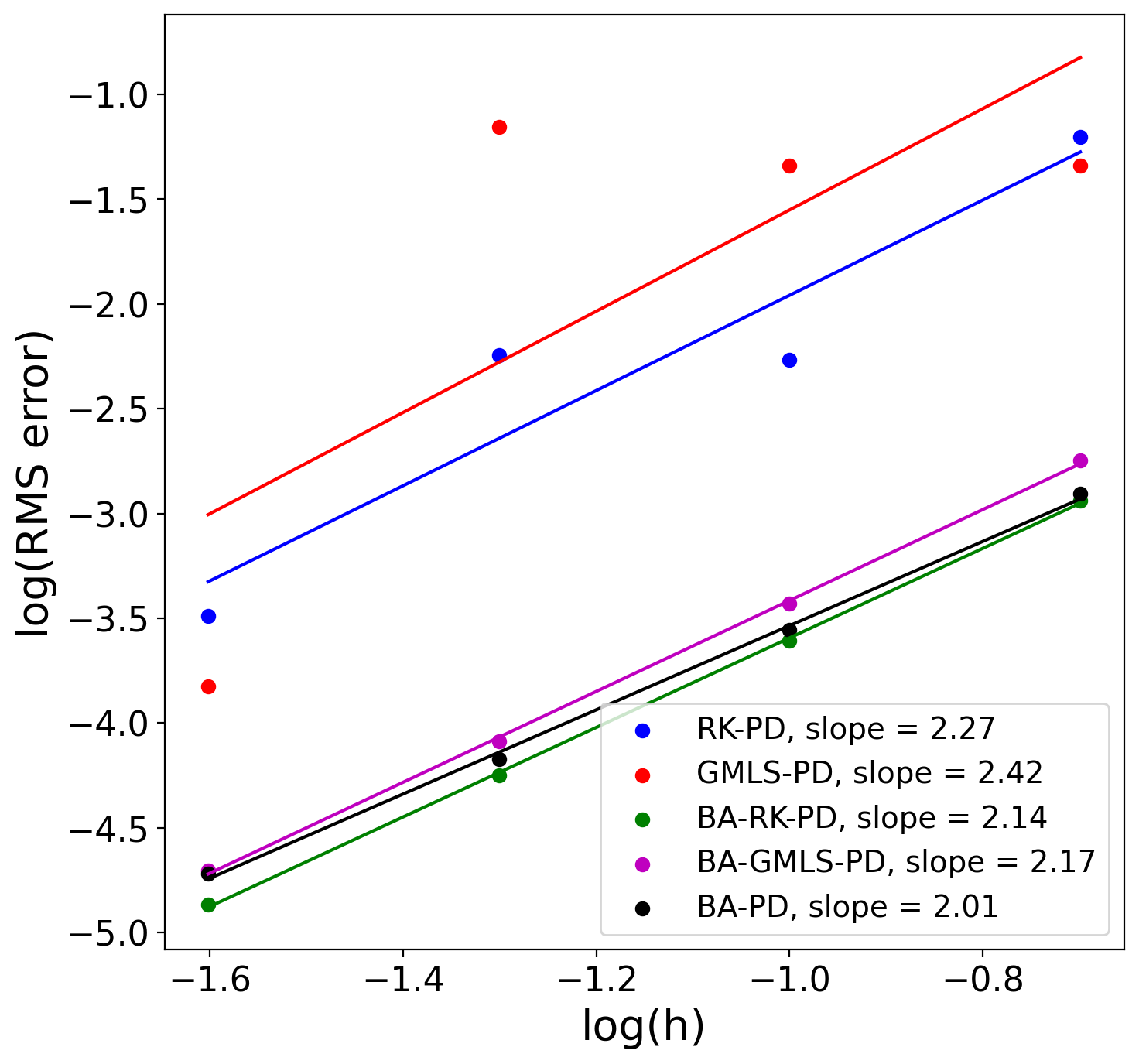}}
  \hspace*{0.3cm}
  \subfloat[][$\delta = 3.5 \, h$]{\includegraphics[height=0.28\textwidth]{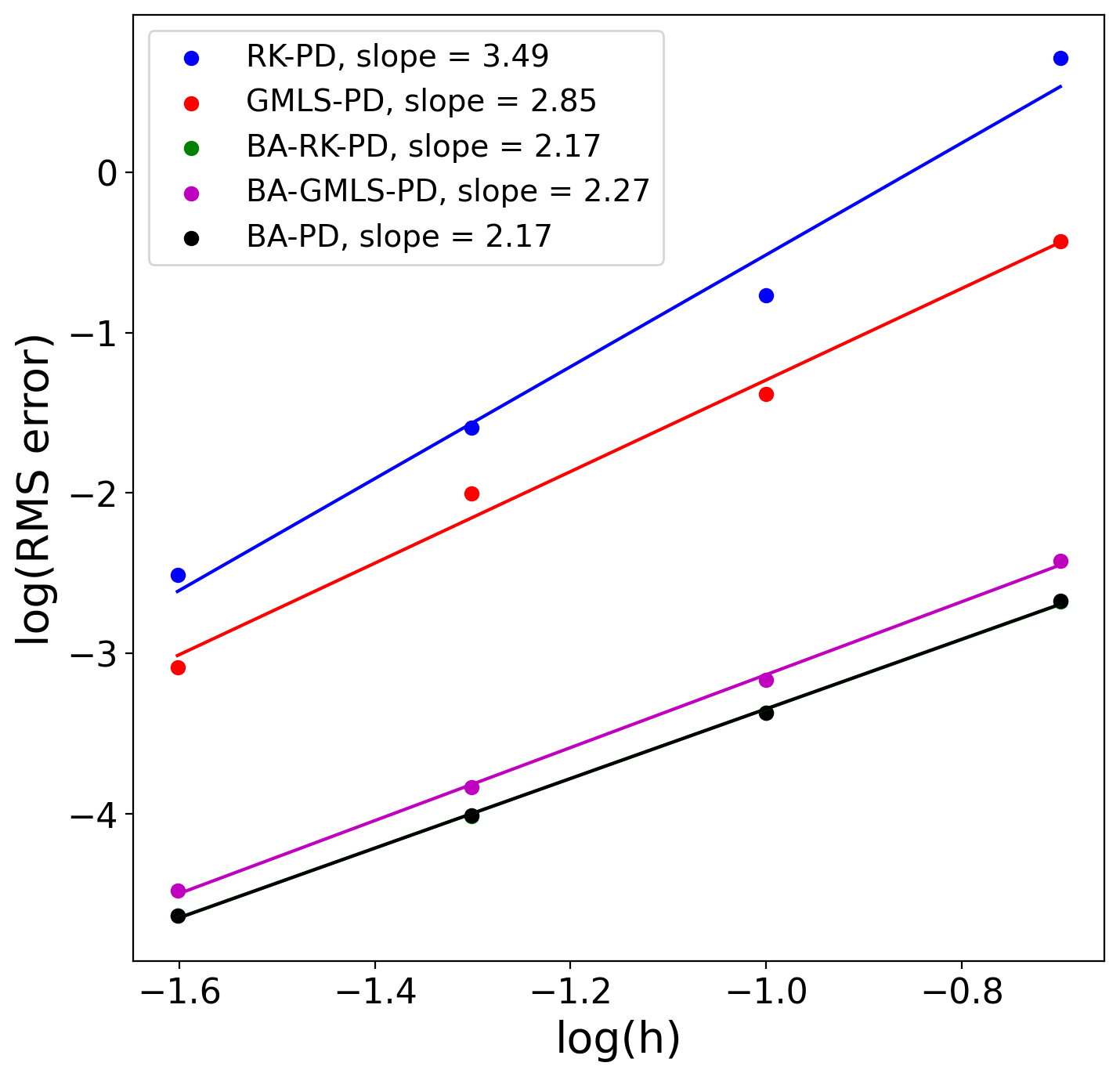}}
  \hspace*{0.3cm}
  \subfloat[][$\delta = 4.25 \, h$]{\includegraphics[height=0.28\textwidth]{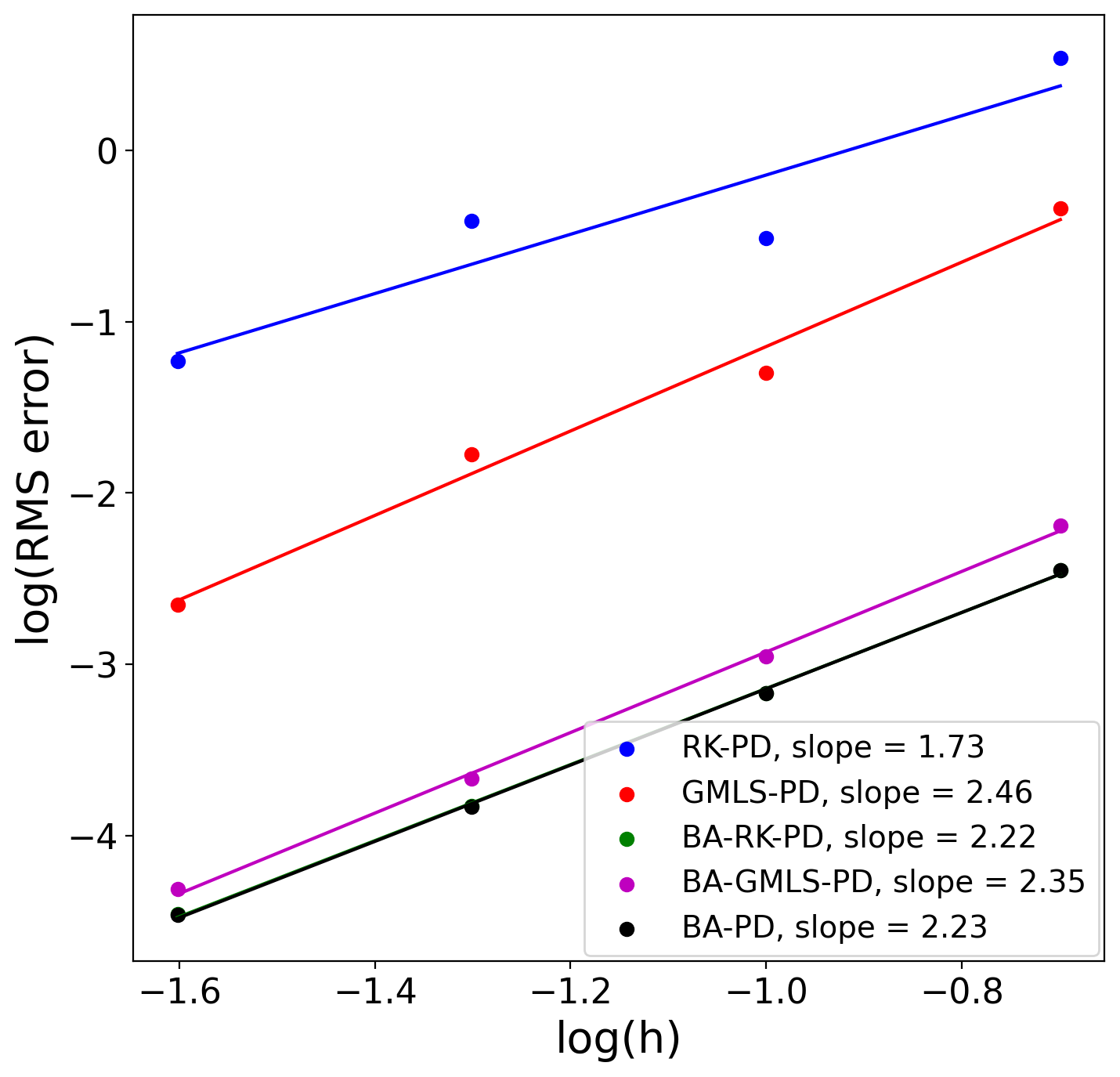}}
  \caption{Convergence of quadratic models in the 2D manufactured problem is compared between different horizon sizes. The horizon size has a large effect on RK-PD and GMLS-PD, while the bond-associated versions show robustness with respect to the size of horizon. Plots show the displacement RMS error values.}
  \label{fig:square-horizon}
\end{figure*}

The horizon-size effect is studied in \cref{fig:square-horizon}, where three different horizons are considered for the quadratic case with non-uniform discretization. That is, for a fixed nodal spacing $h$, the larger the horizon, the more neighbors each node has. RK-PD errors increase with bigger horizon in this case, while GMLS-PD has issues with the smallest horizon. The bond-associated models are stable across the board, thus demonstrating robustness with respect to the horizon size. This is also a very desirable property in real applications where a user cannot be expected to {\em tune} the horizon size to avoid instabilities in the solution. 

\subsection{Infinite plate with circular hole under far-field bi-axial tension}
\label{subsec:plate}

A classical problem is considered here, in which an infinite plate with a circular hole of radius $a$ is subjected to far-field bi-axial tension $T$, without any applied body force. Under plane strain condition, the following solution can be obtained using Airy stress functions \citep{michell1899direct}:
\begin{align*}
  & u_1(r,\theta) = \frac{T a}{2\mu} \Big[ \frac{1-\nu}{1+\nu} \frac{r}{a} + \frac{a}{r} \Big] \cos\theta , \notag \\
  & u_2(r,\theta) = \frac{T a}{2\mu} \Big[ \frac{1-\nu}{1+\nu} \frac{r}{a} + \frac{a}{r} \Big] \sin\theta .
\end{align*}

\begin{figure}[!htbp]
  \centering
  \includegraphics[width=0.45\textwidth]{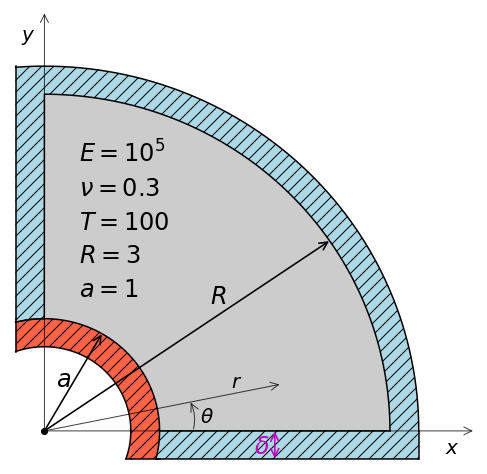}
  \caption{Elastic plate with a circular hole. Blue color indicates the Dirichlet boundary nodes, prescribed with the exact displacement values. Red color shows the free-surface region, where fictitious nodes are placed and given zero stress values.}
  \label{fig:geometry}
\end{figure}

This example serves as an important benchmark as it involves a curvilinear free surface. To set up this problem, the infinite problem is modeled as a finite quarter plate, with the exact displacements prescribed on the non-local collar surrounding the perimeter of the plate, except at the hole area where fictitious nodes are placed to model the free surface. That is, as mentioned previously, the free-surface nodes are supplied with zero stress values that directly contribute to the evaluation of $\nabla_h \cdot \mbt{P}$ for their material neighbors. The free-surface nodes do not contribute to the evaluation of the kinematic variable for their neighbors. \cref{fig:geometry} provides a schematic of the modeling approach for this problem.

\begin{figure*}[!ht]
  \centering
  \subfloat[][Polar - L0]{\includegraphics[height=0.25\textwidth]{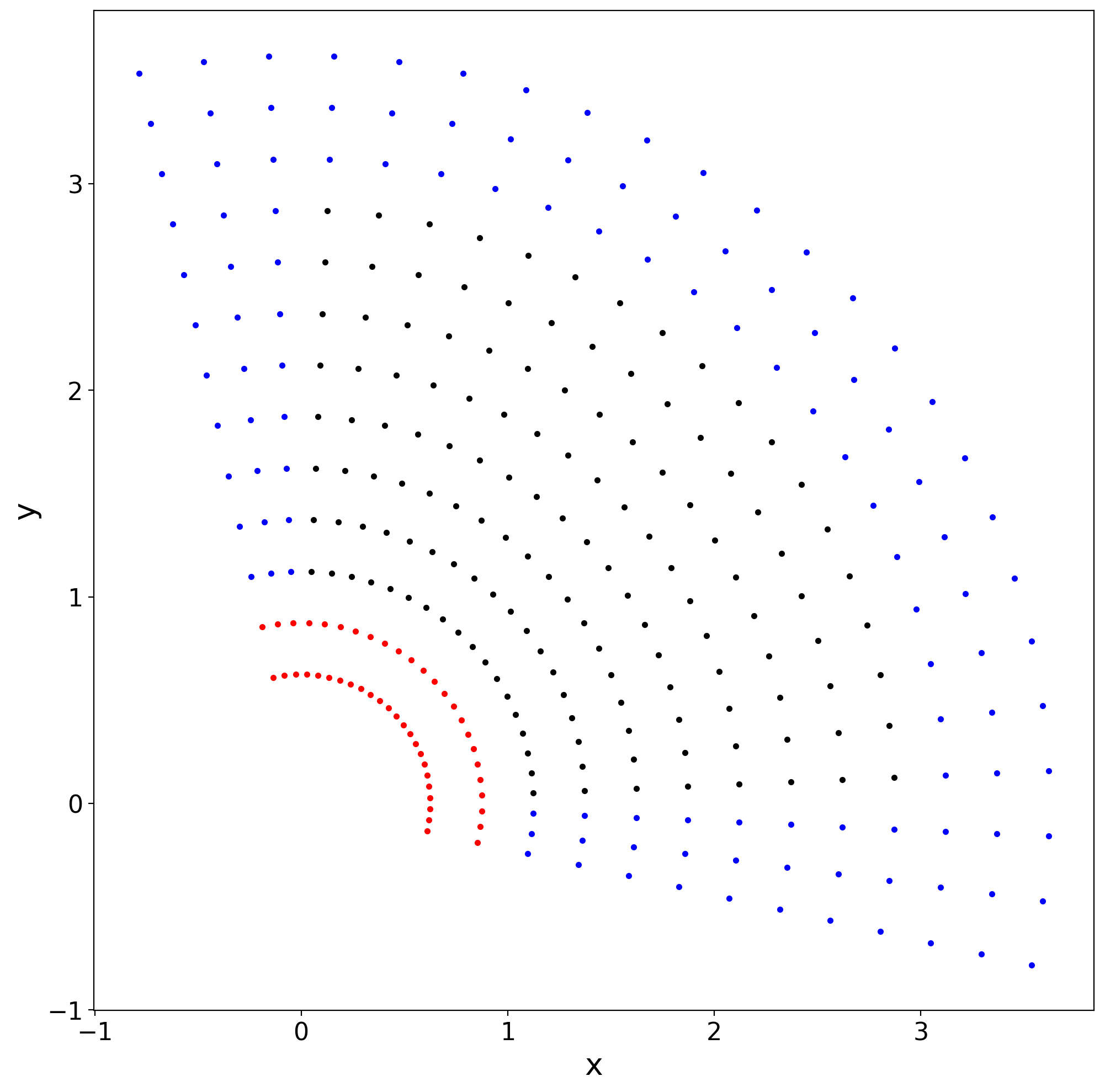}}
  \hspace*{0.1cm}
  \subfloat[][Polar - L1]{\includegraphics[height=0.25\textwidth,trim={2.5cm 0 0 0},clip]{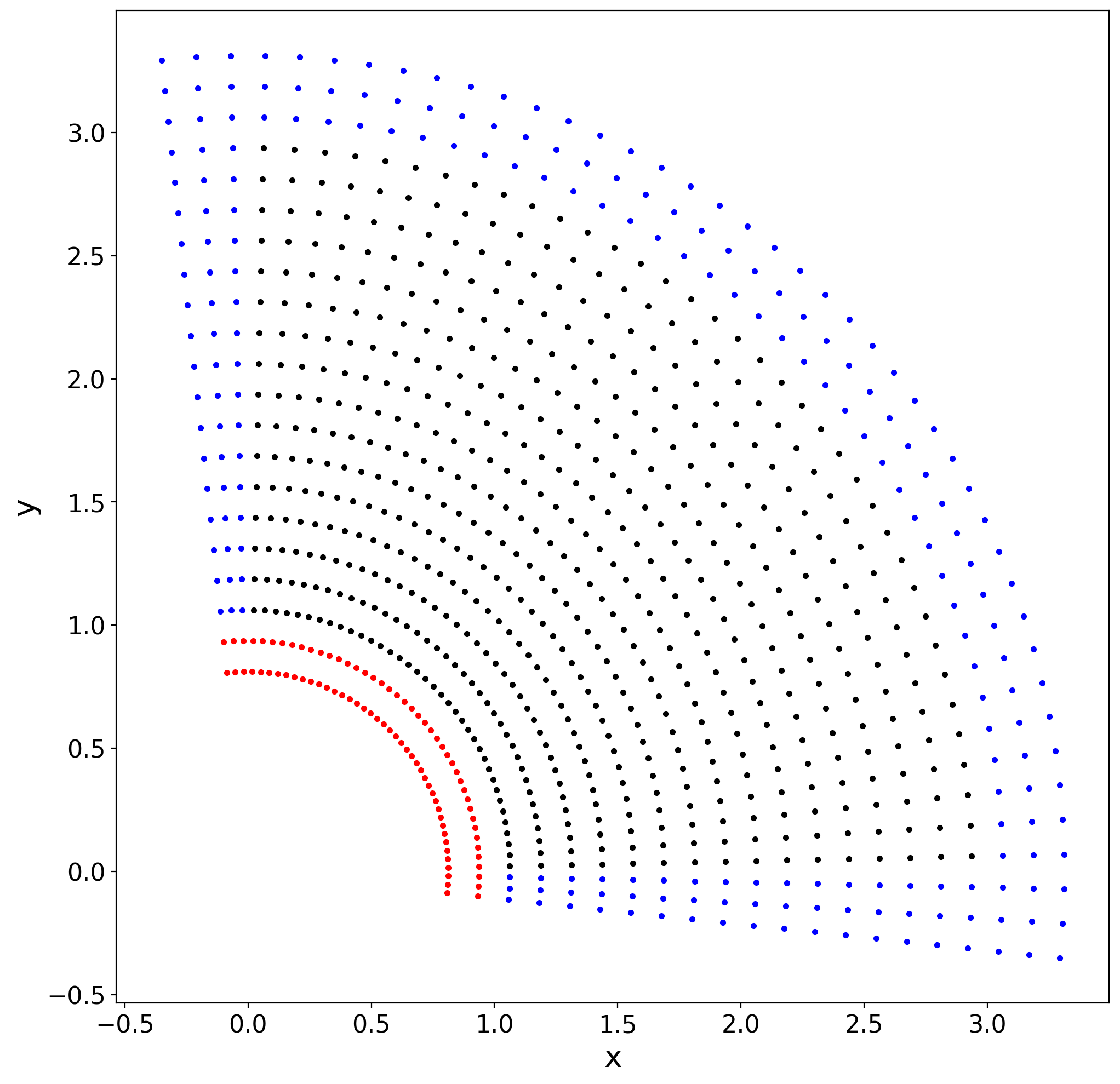}}
  \hspace*{0.1cm}
  \subfloat[][Polar - L2]{\includegraphics[height=0.25\textwidth,trim={2cm 0 0 0},clip]{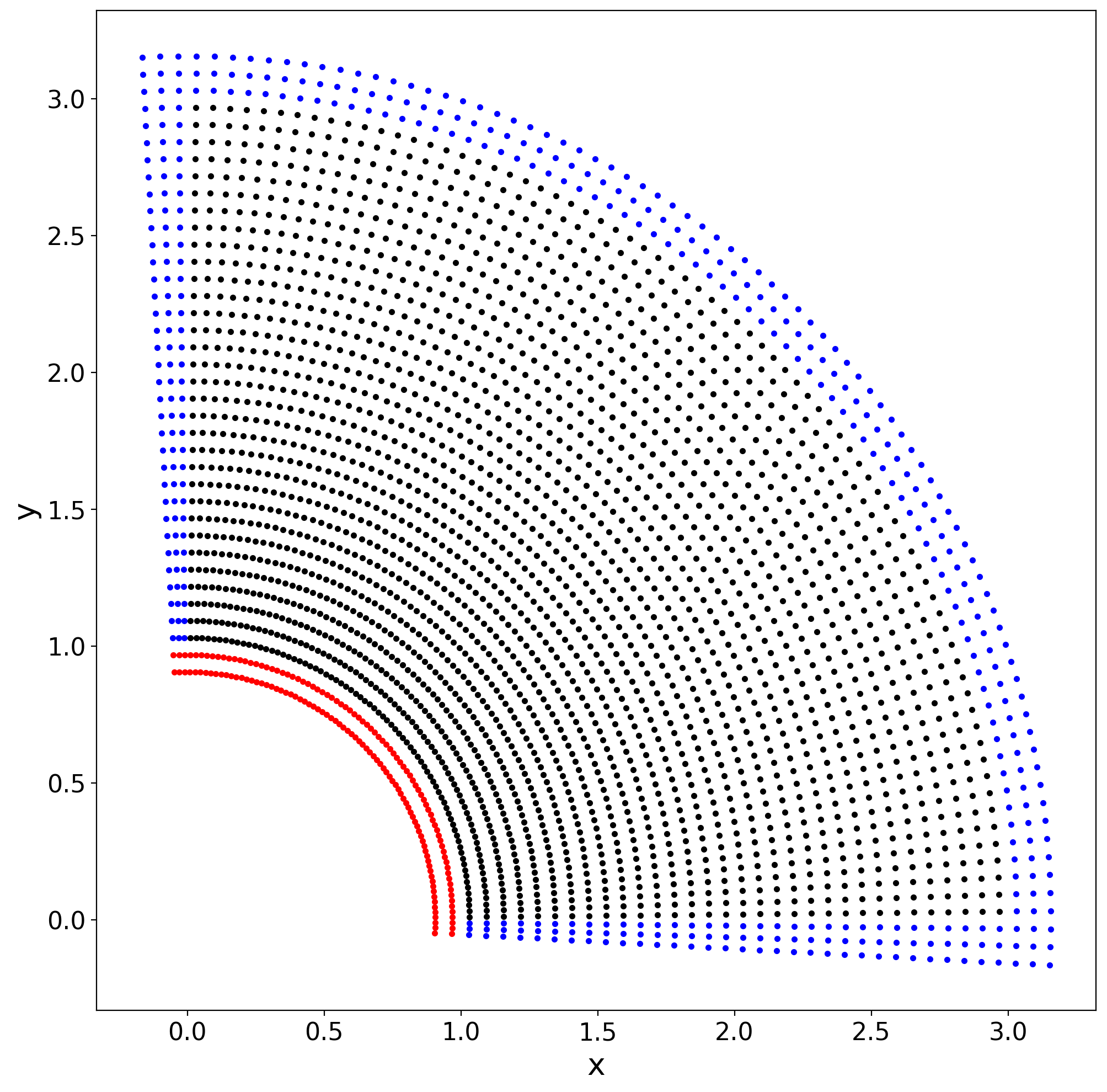}}
  \hspace*{0.1cm}
  \subfloat[][Polar - L3]{\includegraphics[height=0.25\textwidth,trim={2cm 0 0 0},clip]{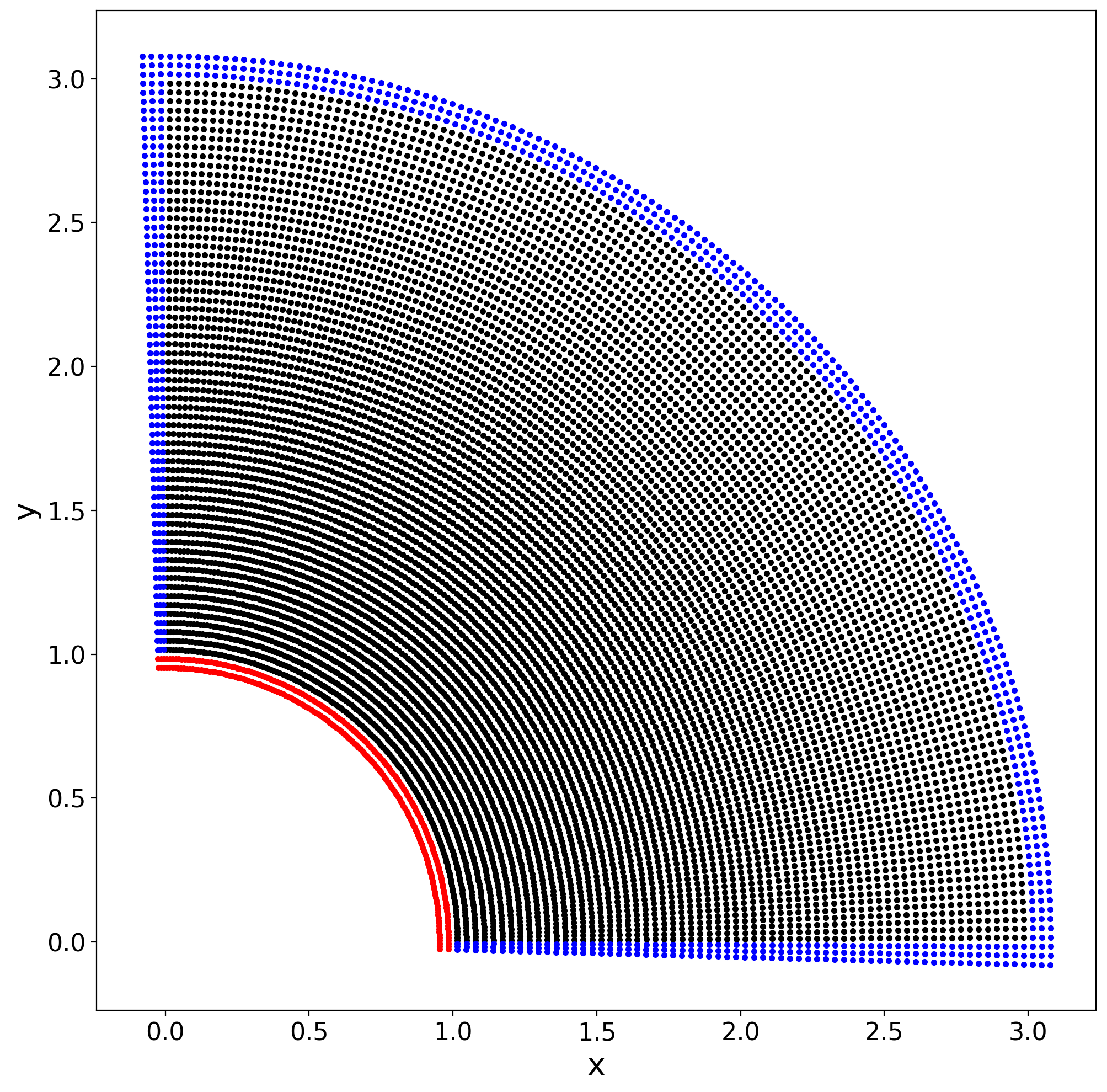}}

  \subfloat[][Triangular - L0]{\includegraphics[height=0.245\textwidth]{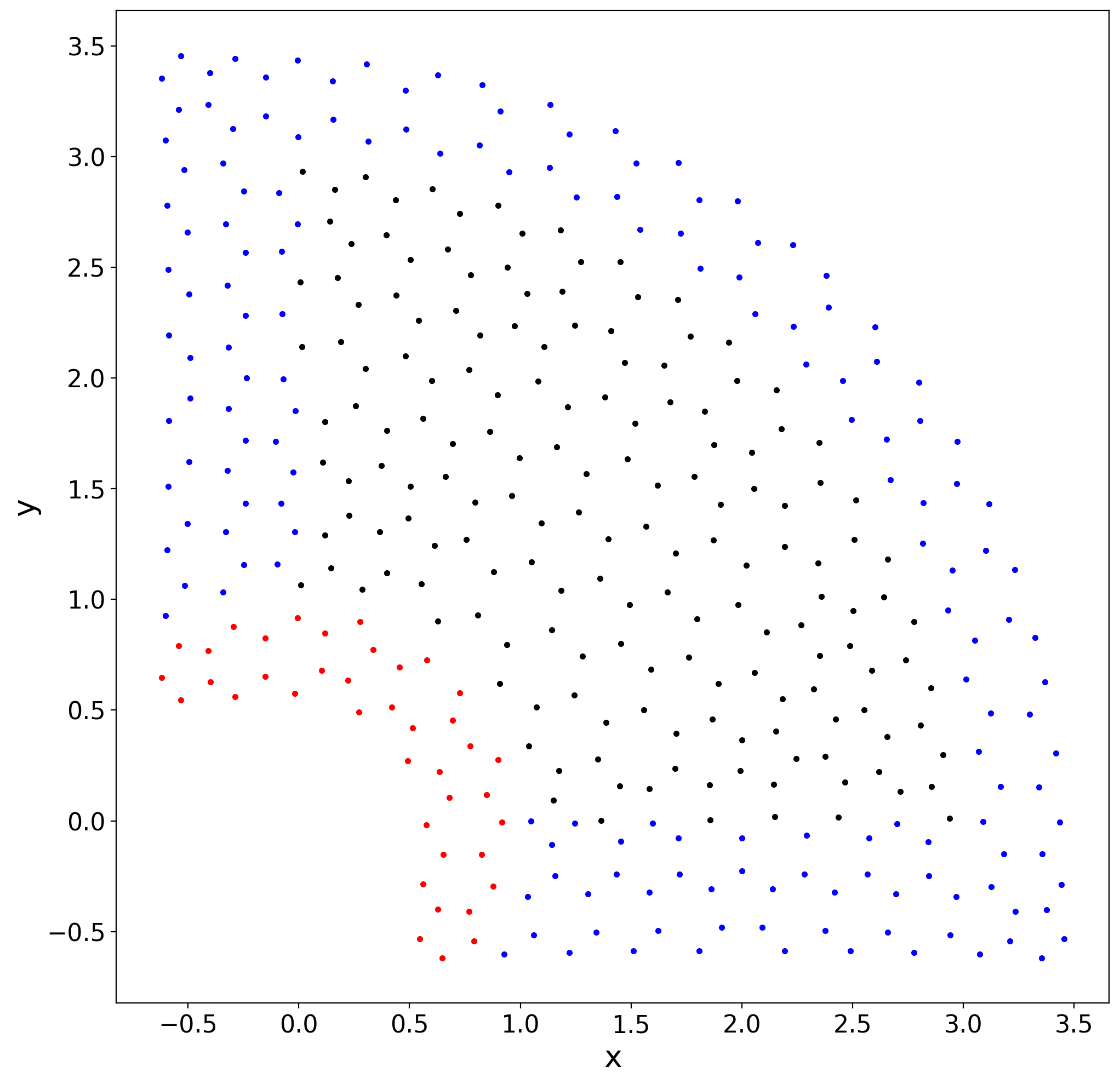}}
  \hspace*{0.1cm}
  \subfloat[][Triangular - L1]{\includegraphics[height=0.245\textwidth,trim={2cm 0 0 0},clip]{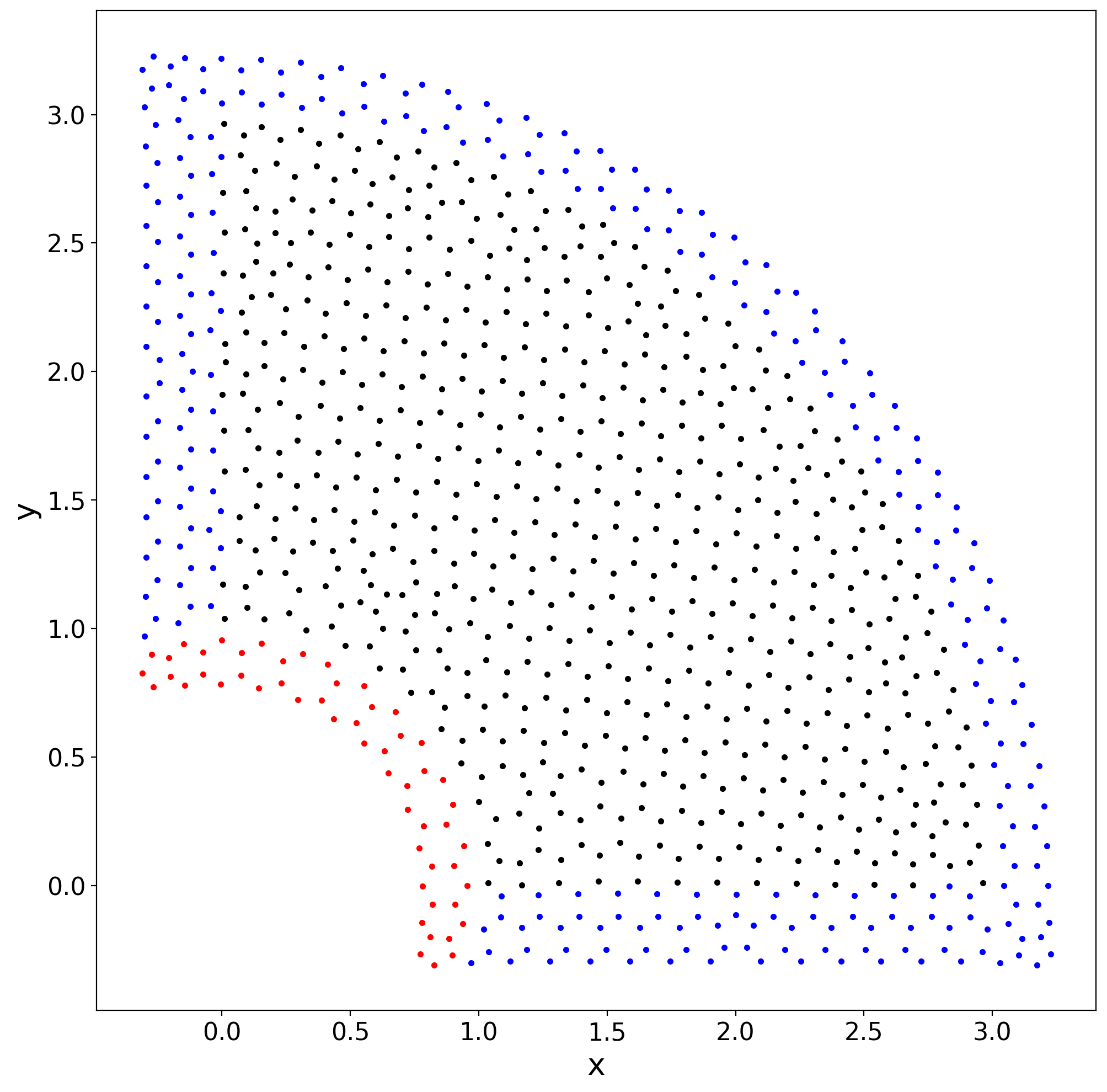}}
  \hspace*{0.1cm}
  \subfloat[][Triangular - L2]{\includegraphics[height=0.245\textwidth,trim={2cm 0 0 0},clip]{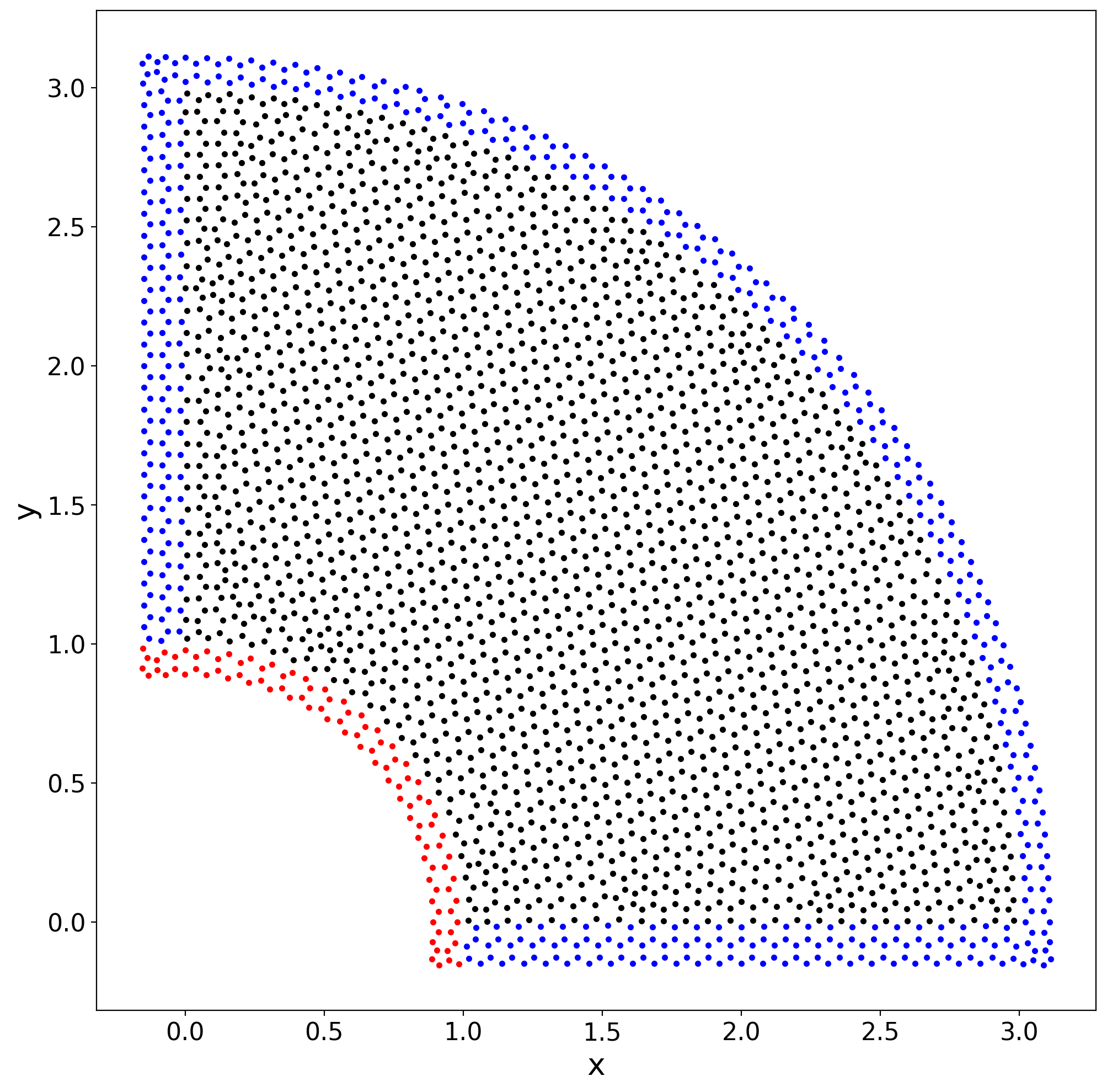}}
  \hspace*{0.1cm}
  \subfloat[][Triangular - L3]{\includegraphics[height=0.245\textwidth,trim={2cm 0 0 0},clip]{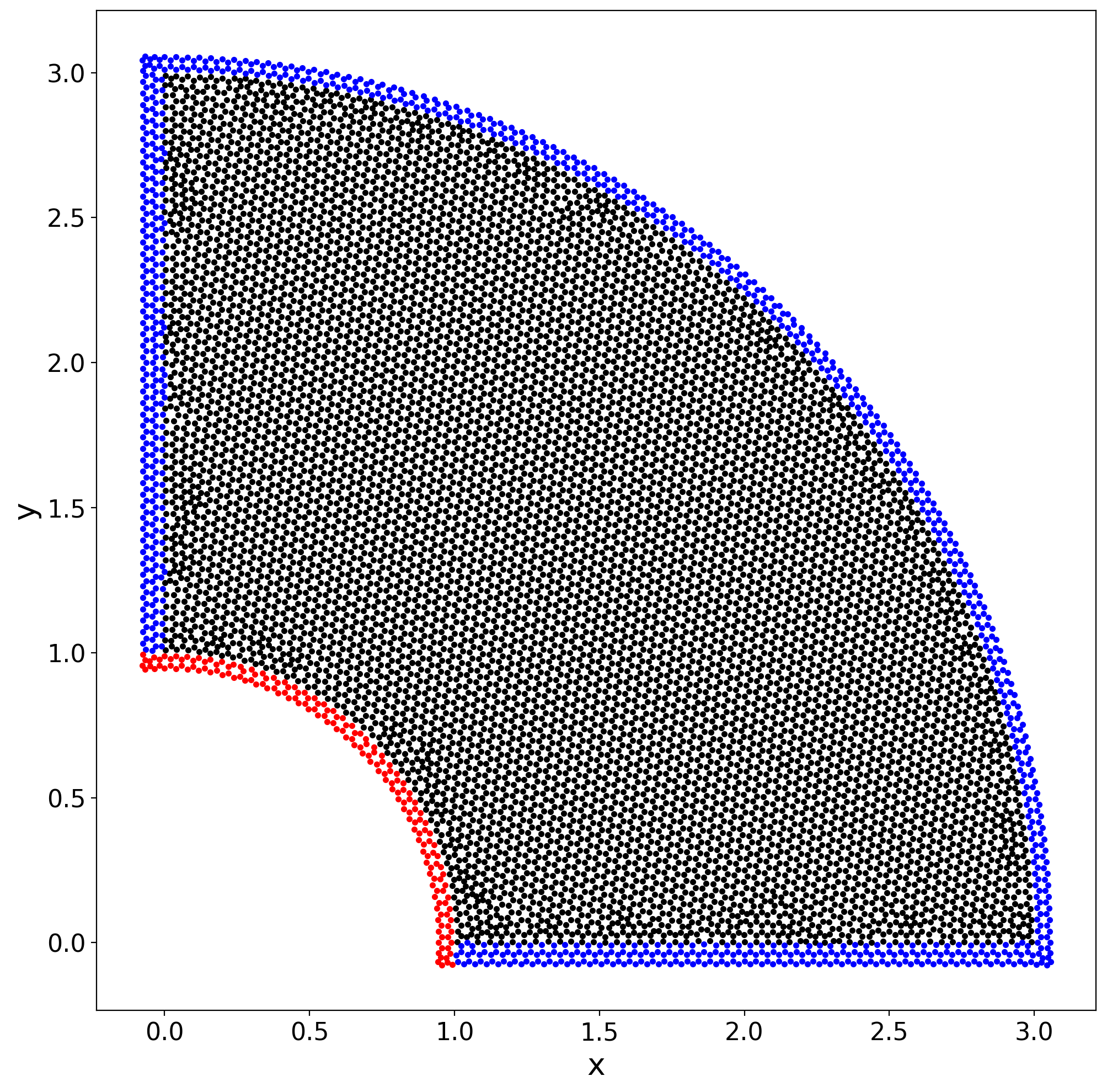}}
  \caption{Nodal discretizations and refinements for the plate with circular hole problem. The blue color denotes the nodes used to prescribe displacement-controlled boundary conditions. The red-colored nodes are used to model free surface.}
  \label{fig:plate-mesh}
\end{figure*}

To obtain nodal discretizations in this example, a finite element mesh is generated first, then each element is replaced by a meshfree node. The meshfree node is placed on the centroid of the element and is associated with its area (volume in 3D). Two strategies are employed to mesh the domain: (1) a structured discretization algorithm, where quadrilateral elements are obtained by having their nodes uniformly distributed along a polar coordinate system, and (2) a semi-unstructured algorithm using {\em pyGmsh}, a Python library based on Gmsh \citep{geuzaine2007gmsh}, where Delaunay triangular meshes are generated. In the remainder of this section, these two schemes are denoted as {\em polar} and {\em triangular} discretizations, respectively. Different levels of mesh refinement are utilized to perform a convergence study. For each case, average node spacing is defined by $$h = \sqrt{\frac{\sum_i^N A_i}{N}},$$ where $A_i$ is element area. The approximate average spacings of $h\approx[0.2, \, 0.1, \, 0.05, \, 0.025]$ are considered for both schemes. \cref{fig:plate-mesh} shows different levels of discretization in this approach.

\begin{figure*}[!htbp]
  \centering
  \includegraphics[width=0.65\textwidth]{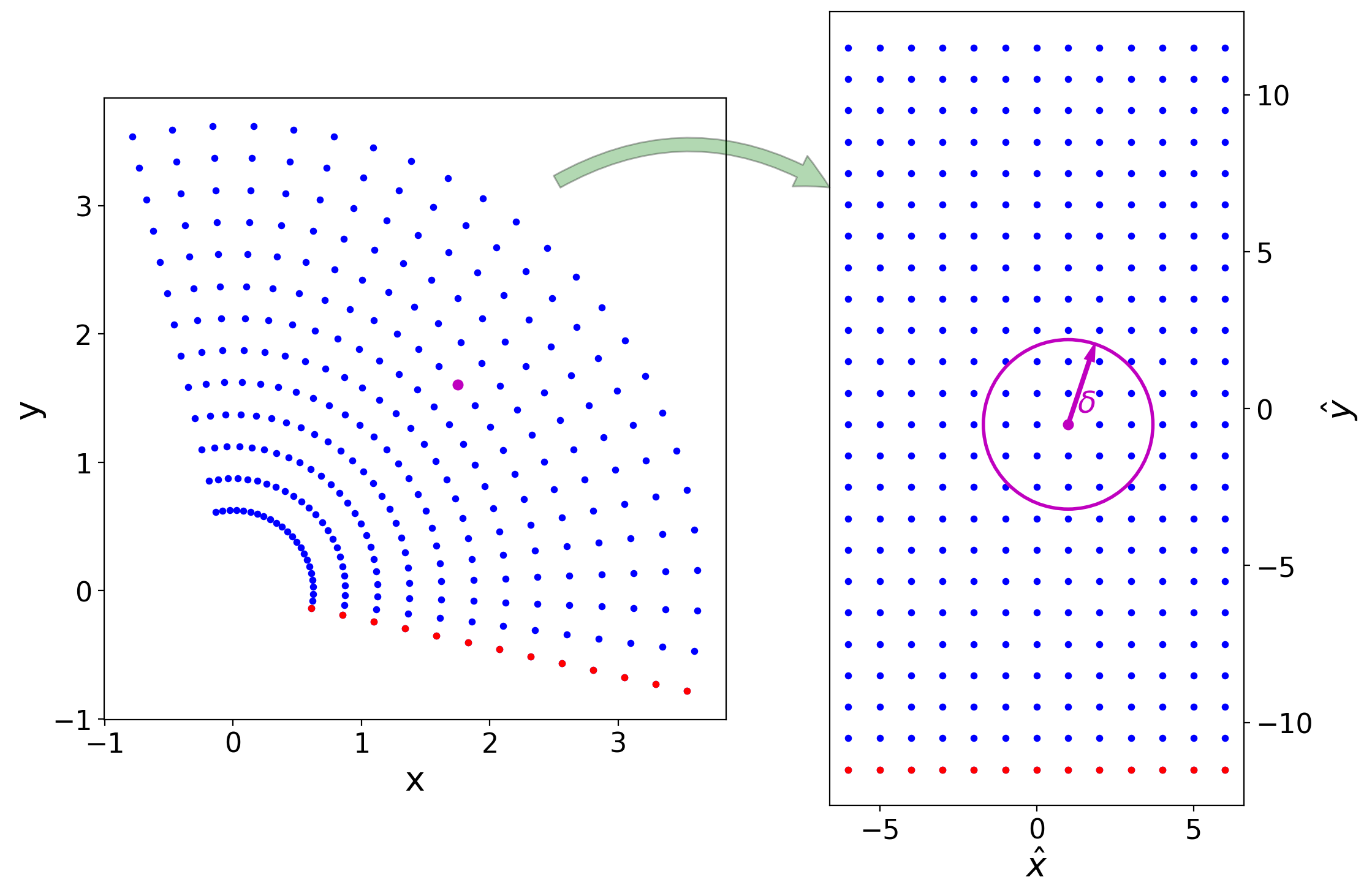}
  \caption{A mapping is used to define the neighbor sets in a parametric space, resulting in a uniform number of neighbors for each node (except for near the boundaries).}
  \label{fig:mapping}
\end{figure*}

The polar discretization scheme results in having more concentrated nodes in areas closer to the hole. Using a fixed horizon size in the {\em physical} domain to construct neighborhoods would result in a significantly larger neighbor sets near the hole, reducing efficiency, and potentially damaging the robustness of the RK and GMLS algorithms as they depend on the number of neighbors (cf. \cref{fig:square-horizon}). Using circular-shaped neighborhoods with varying horizons in the physical space, which results in cases where a node $I$ is a neighbor of $J$, but $J$ is not a neighbor of $I$, has been found problematic \citep{ren2016dual}. To obtain neighbor sets with the same number of neighbors (except near the boundaries) regardless of where the neighborhoods are formed, while keeping the consistency between neighbors (i.e., if a bond $IJ$ exists, $JI$ also exists), a mapping algorithm is used. That is, the nodes are transformed from the physical space to a {\em parametric} space, where a uniform spacing of 1 is obtained between the adjacent nodes. The neighborhoods are formed in the parametric space. \cref{fig:mapping} illustrates this concept. In the parametric space, the horizon size for the linear, quadratic, and cubic formulations is chosen as $\delta = 1.75$, $\delta = 2.75$, and $\delta = 3.75$, respectively. Note that while the influence functions are evaluated in the parametric space, the quadrature weights are computed directly in the physical space. 

For the triangular discretization scheme, the neighborhoods are formed in the physical space. In this case, the horizon size for the linear, quadratic, and cubic methods is chosen as $\delta = 2.25 h$, $\delta = 3.25 h$, and $\delta = 4.25 h$, respectively.

\begin{figure*}[!ht]
  \centering
  \subfloat[][Polar discretization - Linear]{\includegraphics[height=0.39\textwidth]{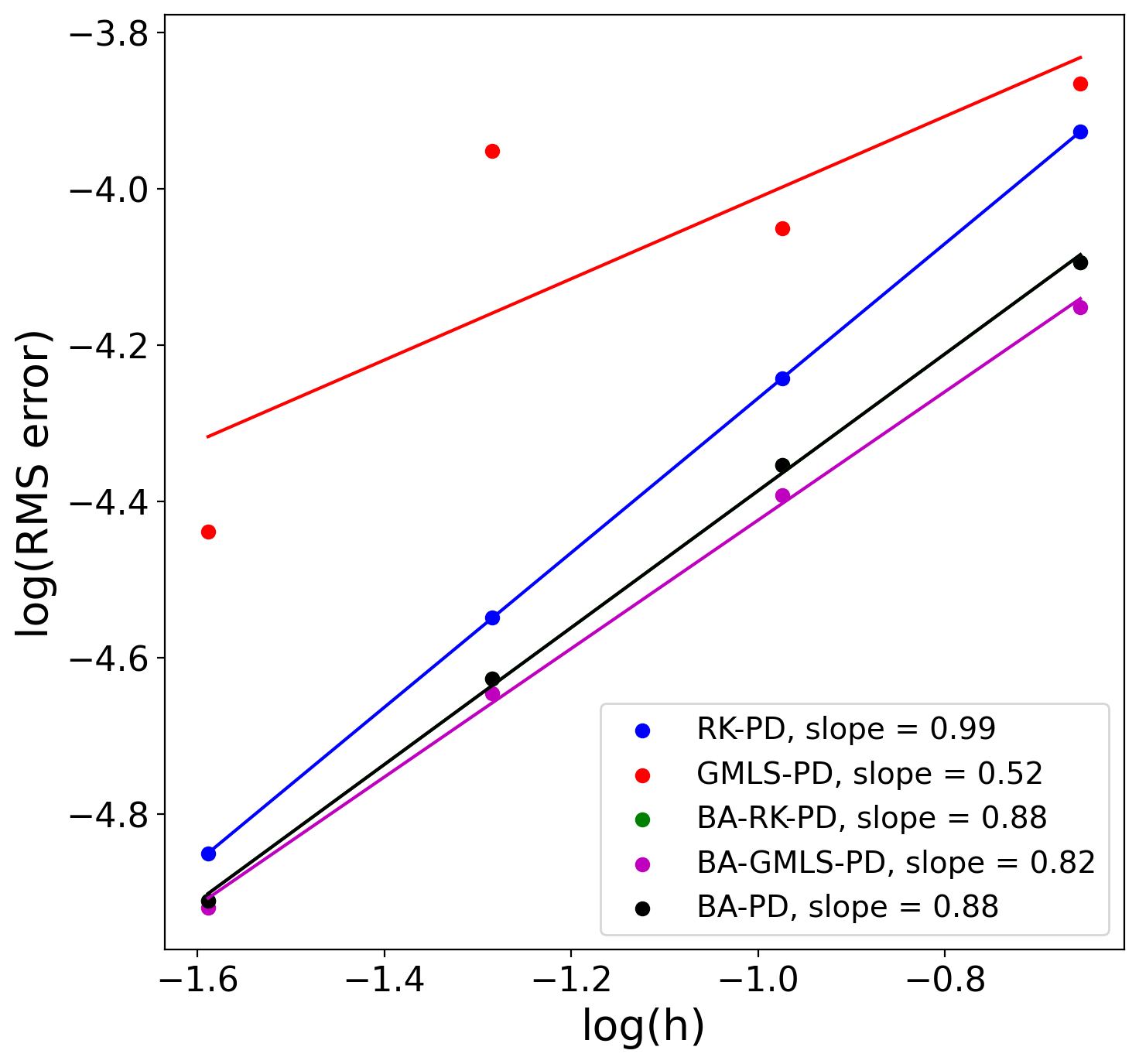}}
  \hspace*{1cm}
  \subfloat[][Triangular discretization - Linear]{\includegraphics[height=0.39\textwidth]{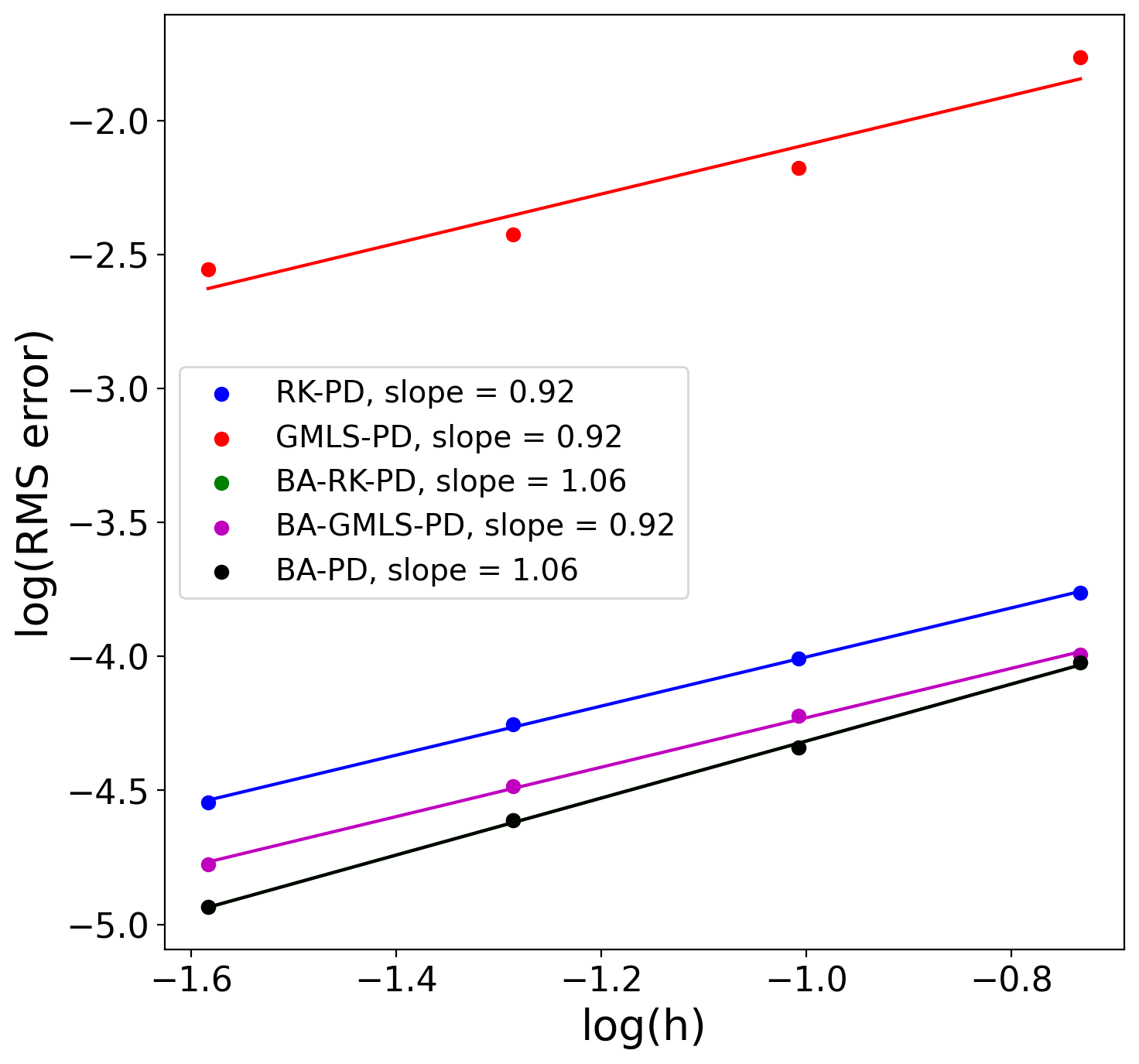}}

  \subfloat[][Polar discretization - Quadratic]{\includegraphics[height=0.39\textwidth]{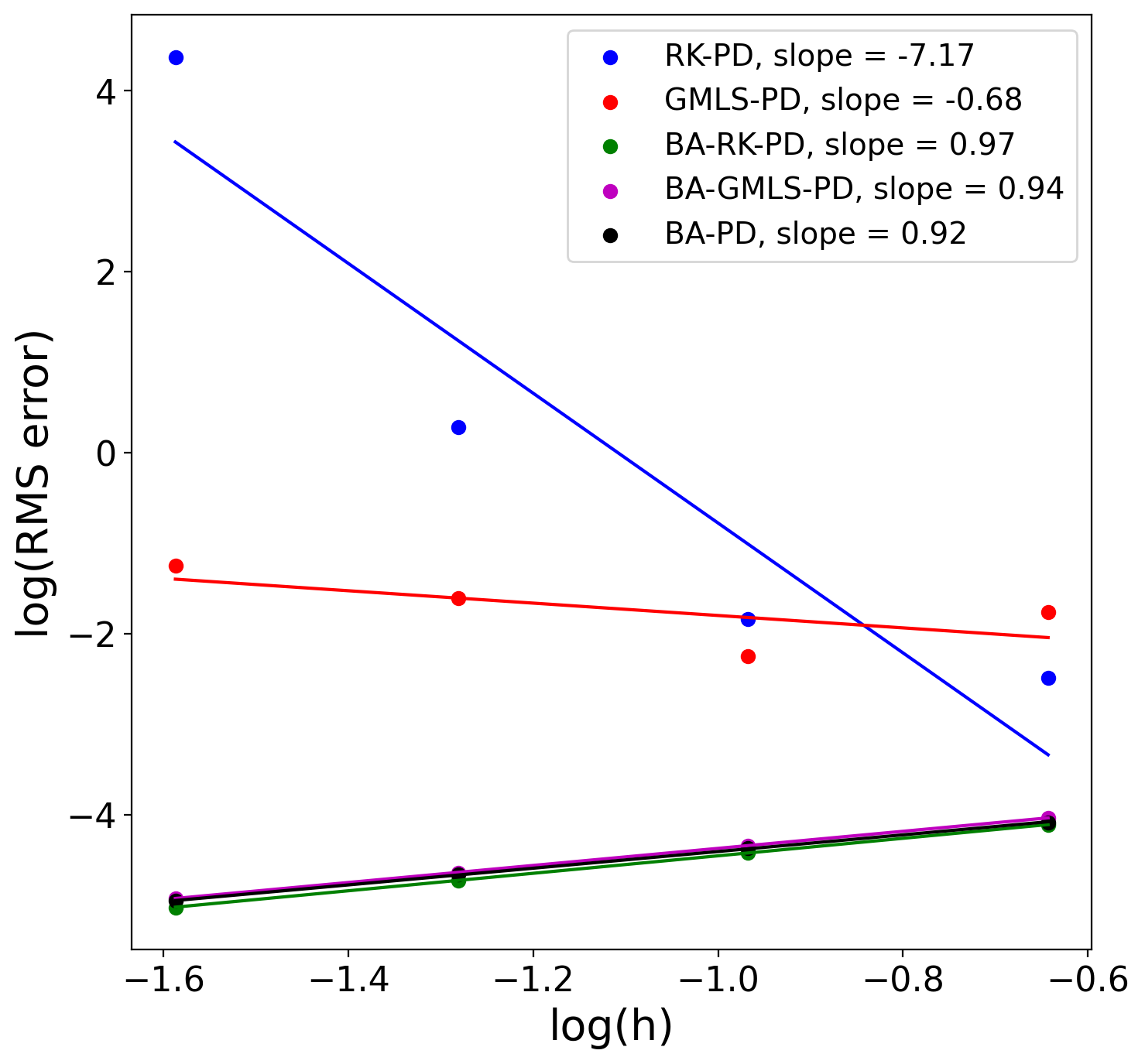}}
  \hspace*{1cm}
  \subfloat[][Triangular discretization - Quadratic]{\includegraphics[height=0.39\textwidth]{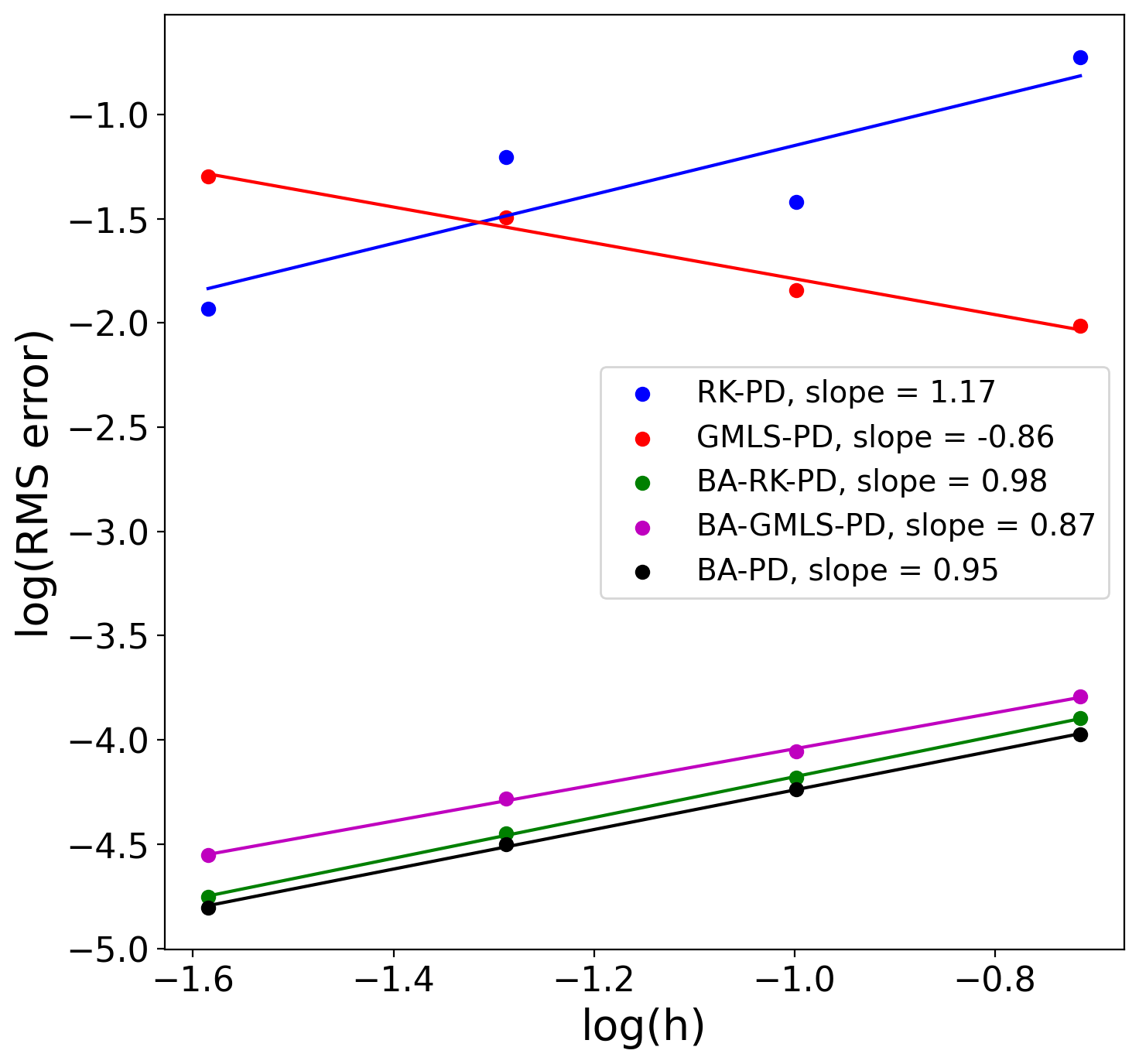}}

  \subfloat[][Polar discretization - Cubic]{\includegraphics[height=0.39\textwidth]{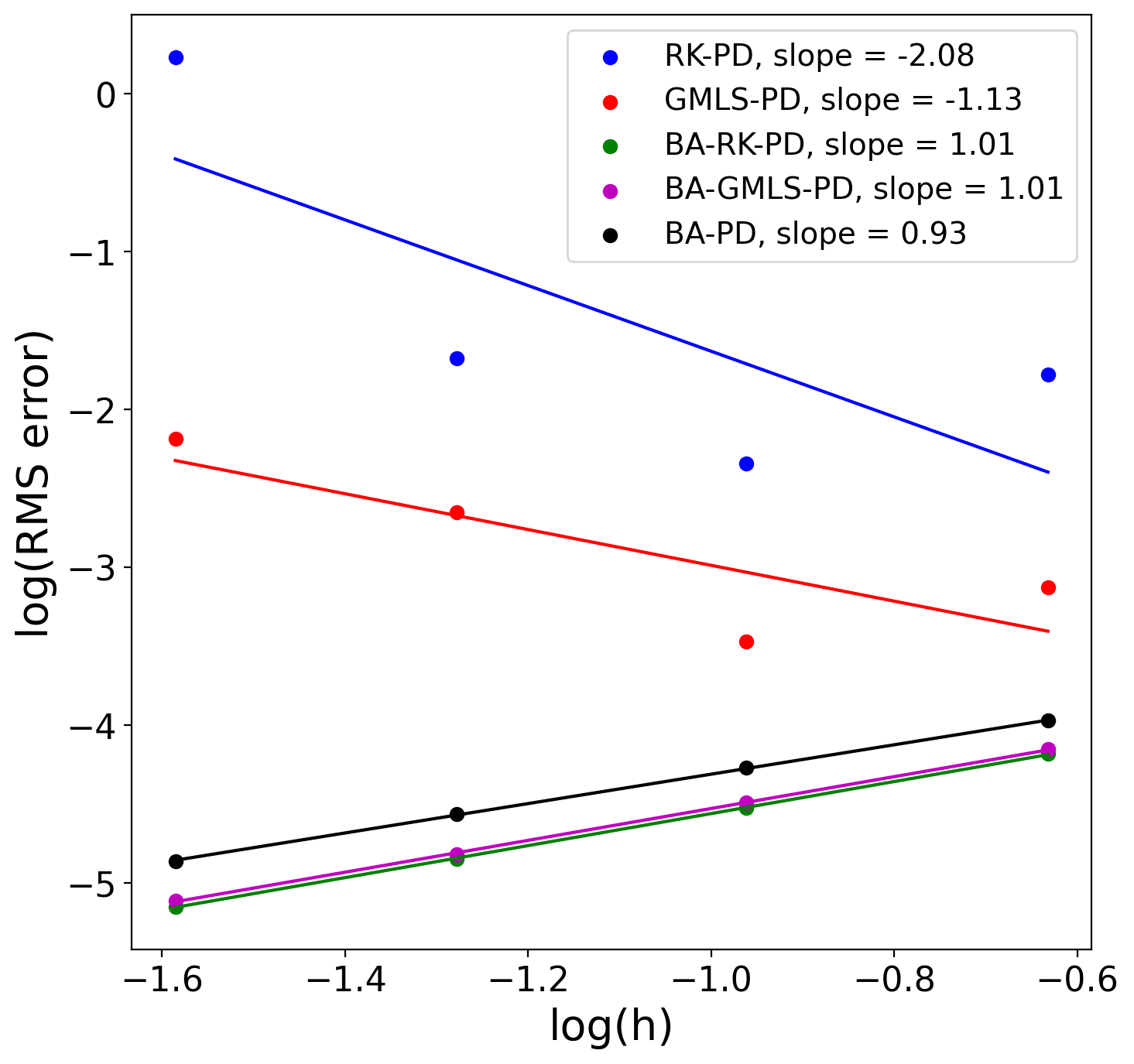}}
  \hspace*{1cm}
  \subfloat[][Triangular discretization - Cubic]{\includegraphics[height=0.39\textwidth]{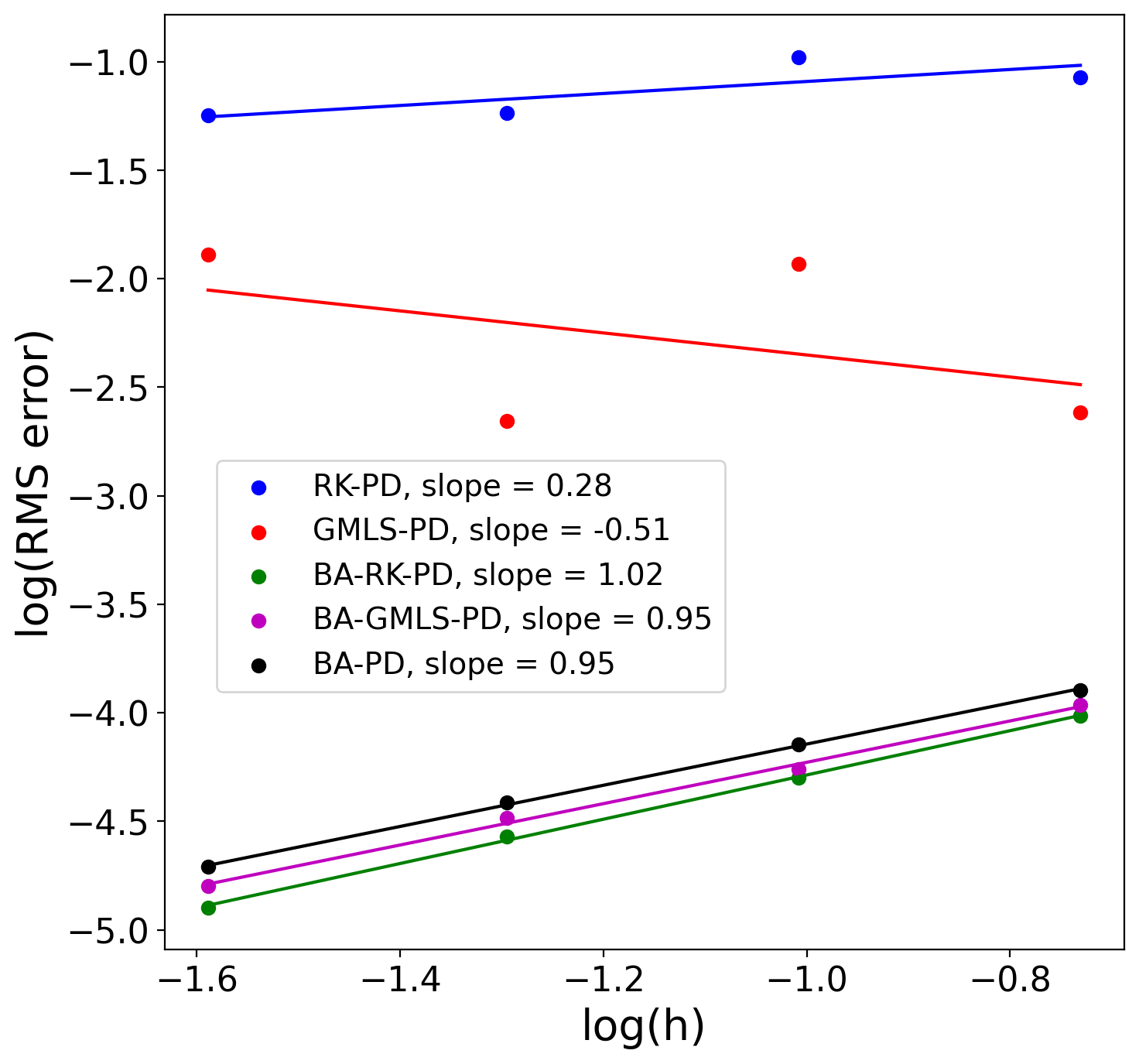}}
  \caption{Two-dimensional convergence test for RK-PD, GMLS-PD, BA-RK-PD, BA-GMLS-PD, and BA-PD in the plate with a circular hole problem. Linear, quadratic, and cubic formulations are tested on the polar and triangular discretizations. BA-PD and BA-RK-PD overlap in (a--b). Plots show the displacement RMS errors.}
  \label{fig:plate-convergence}
\end{figure*}

\cref{fig:plate-convergence} shows the displacement RMS error values for different permutations of the discretizations and models. The bond-associated models obtained a near-linear convergence rate for all cases, which demonstrates asymptotically compatible convergence to the local solution. Note that the solution in this problem is not smooth since there is a jump in stresses near the free surface; therefore, we cannot expect to obtain a higher-order convergence rate with the non-local approach taken. The linear RK-PD case also shows asymptotically compatibility; however, as soon as the number of neighbors increases (i.e., larger horizon), neither of the quadratic or cubic methods are able to solve the problem. The linear GMLS-PD convergence rate is closer to 1 than its quadratic and cubic versions. Except for the linear model on the polar mesh, high errors are seen in others, rooted to the noted instability issue. In this problem, some improvement is obtained by increasing the order of the bond-associated variants, by examining the cubic results. However, the improvements achieved through the stabilization corrections are significantly more tangible.

\begin{figure*}[!ht]
  \centering
  \subfloat[][RK-PD]{\includegraphics[height=0.4\textwidth]{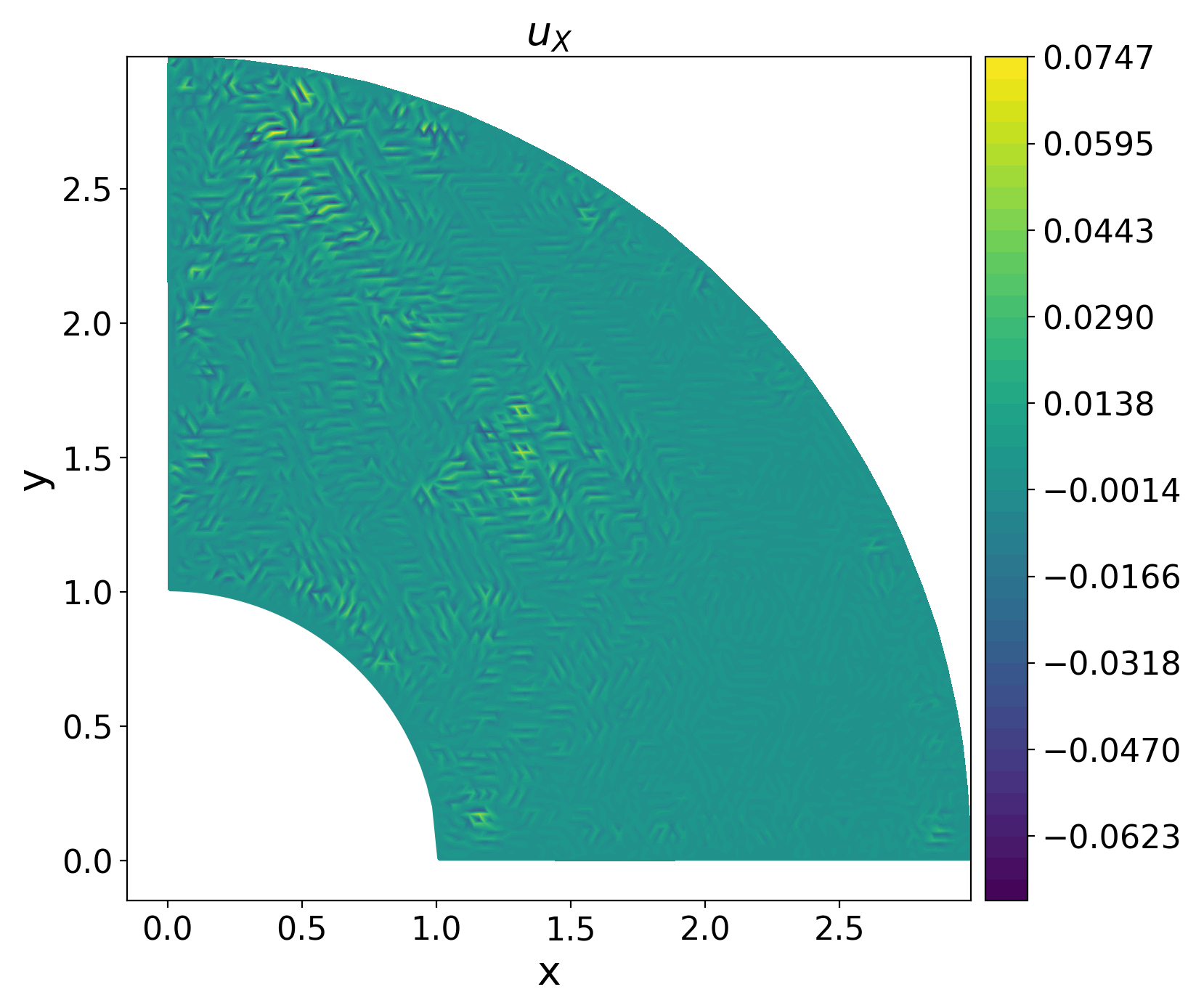}}
  \hspace*{0.5cm}
  \subfloat[][GMLS-PD]{\includegraphics[height=0.4\textwidth]{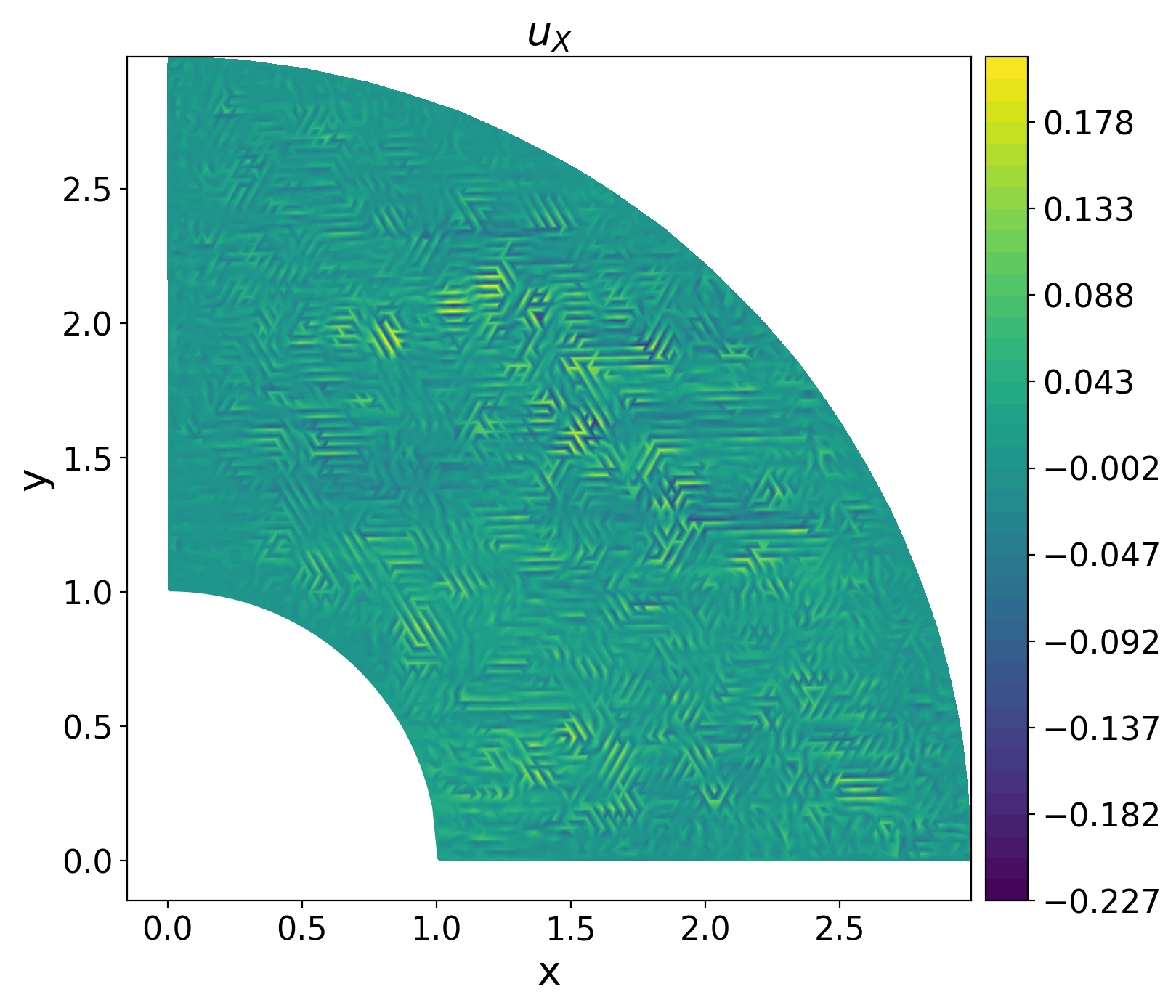}}

  \subfloat[][BA-RK-PD]{\includegraphics[height=0.4\textwidth]{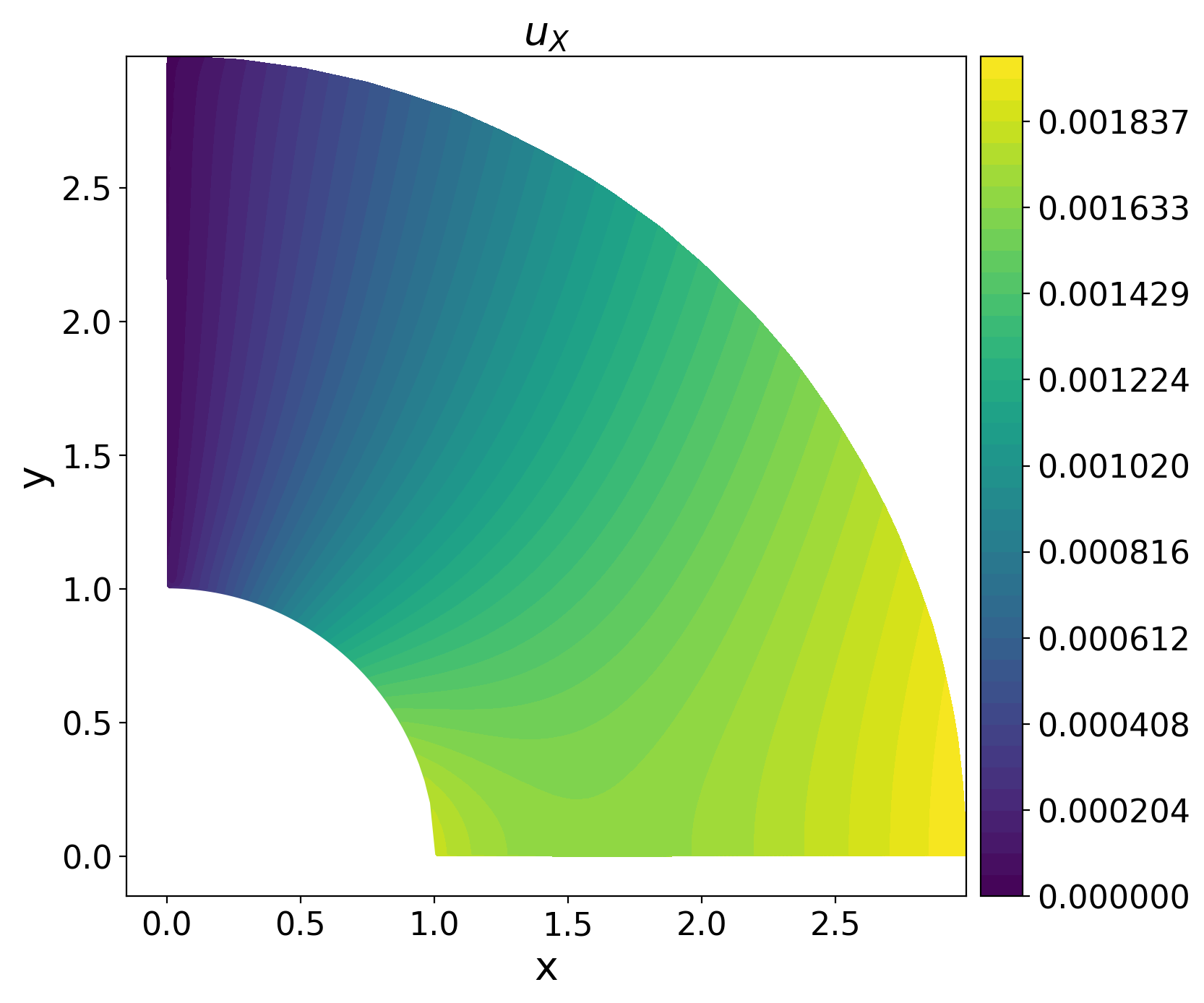}}
  \hspace*{0.5cm}
  \subfloat[][BA-GMLS-PD]{\includegraphics[height=0.4\textwidth]{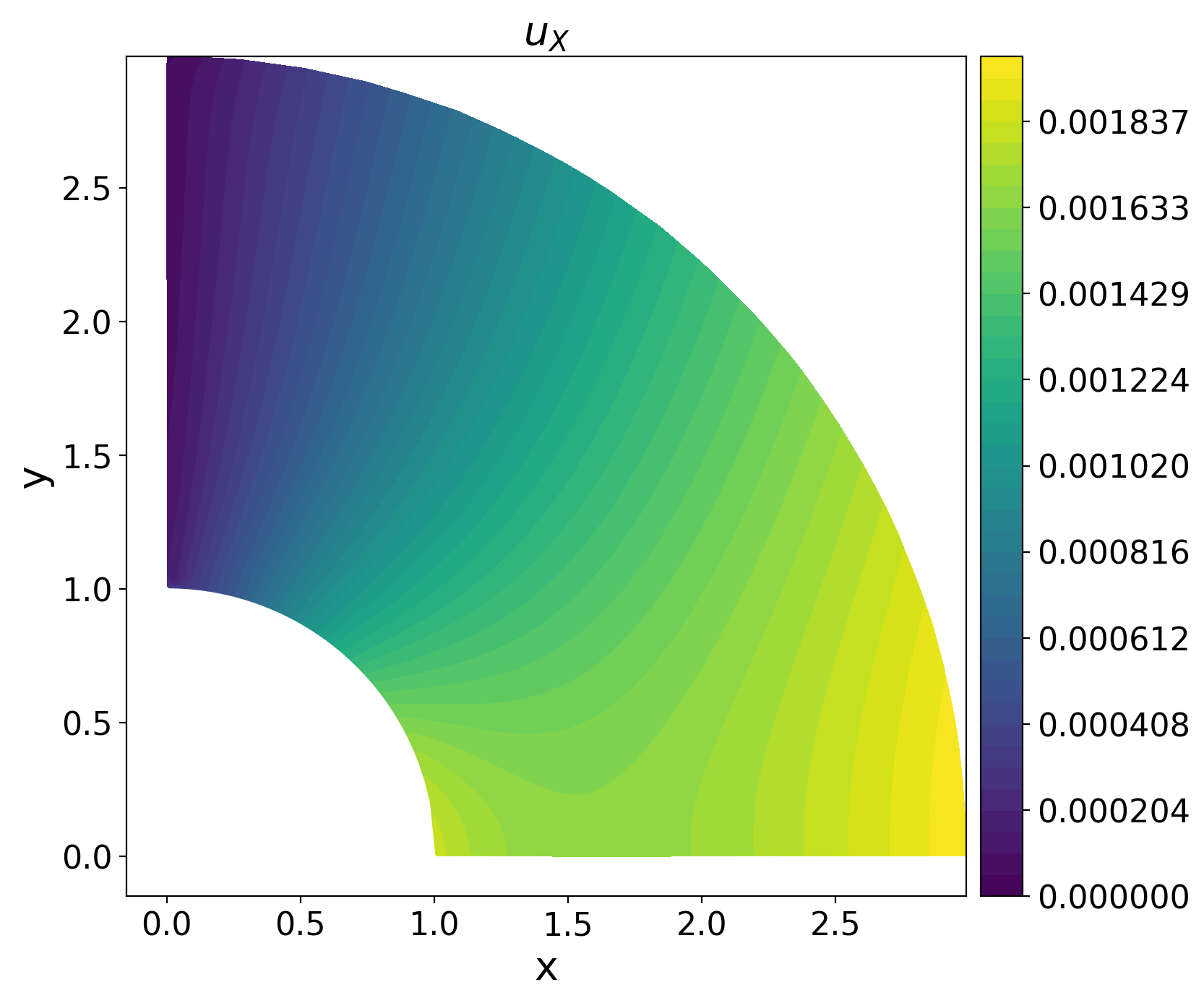}}
  \caption{Horizontal-displacement contours in the plate with hole problem using the quadratic formulations and the L3 triangular discretization. Only the bond-associated solutions (c-d) are in good agreement with the exact solution (not shown here).}
  \label{fig:plate-contours}
\end{figure*}

\cref{fig:plate-contours} illustrates the stability problem of the base models, by showing the solution for the quadratic models on the most refined triangular mesh. While the bond-associated methods capture the exact solution in this problem, the base models result in high errors.

\begin{figure*}[!ht]
  \centering
  \subfloat[][Polar discretization - Quadratic]{\includegraphics[height=0.4\textwidth]{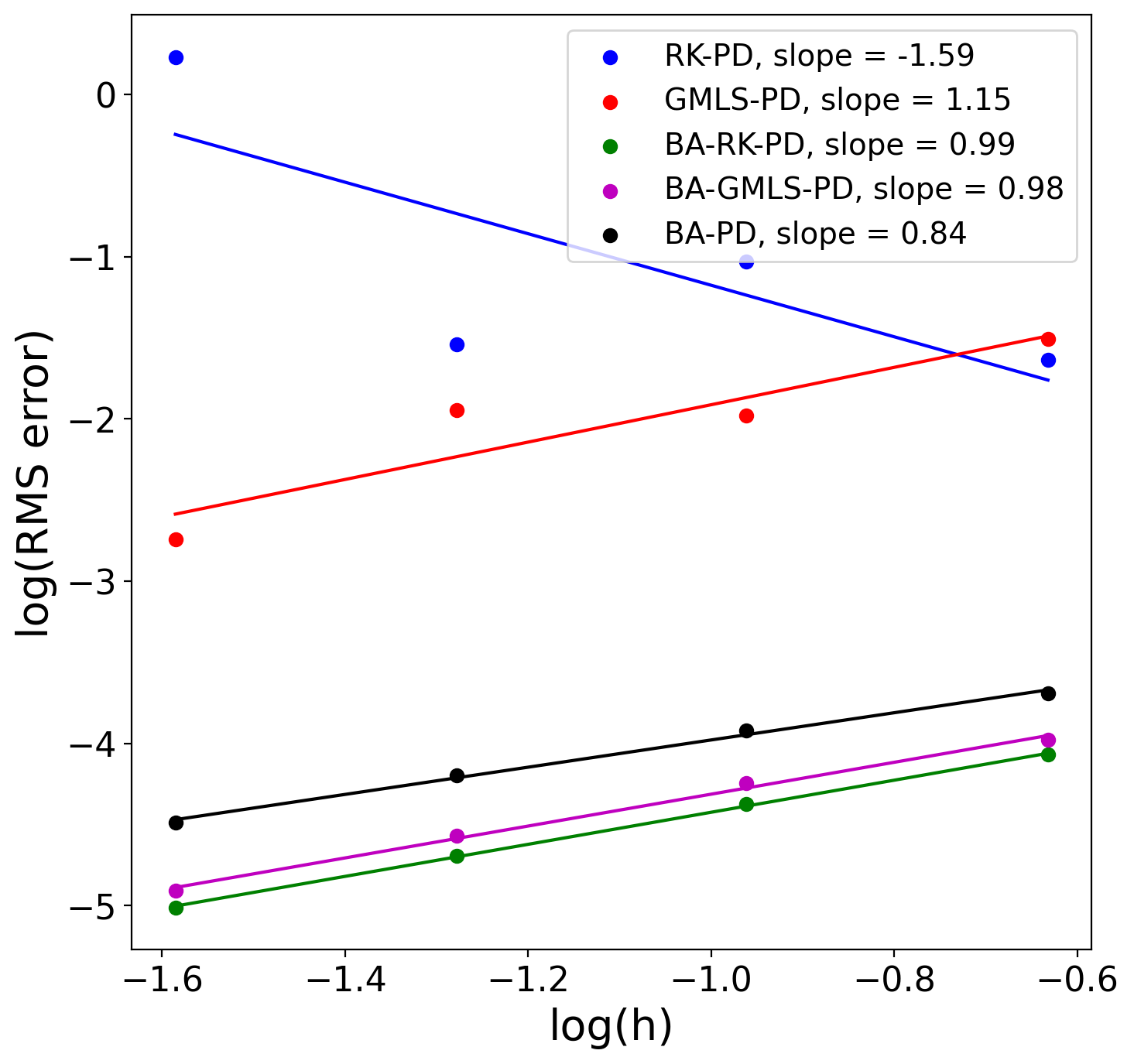}}
  \hspace*{1cm}
  \subfloat[][Triangular discretization - Quadratic]{\includegraphics[height=0.4\textwidth]{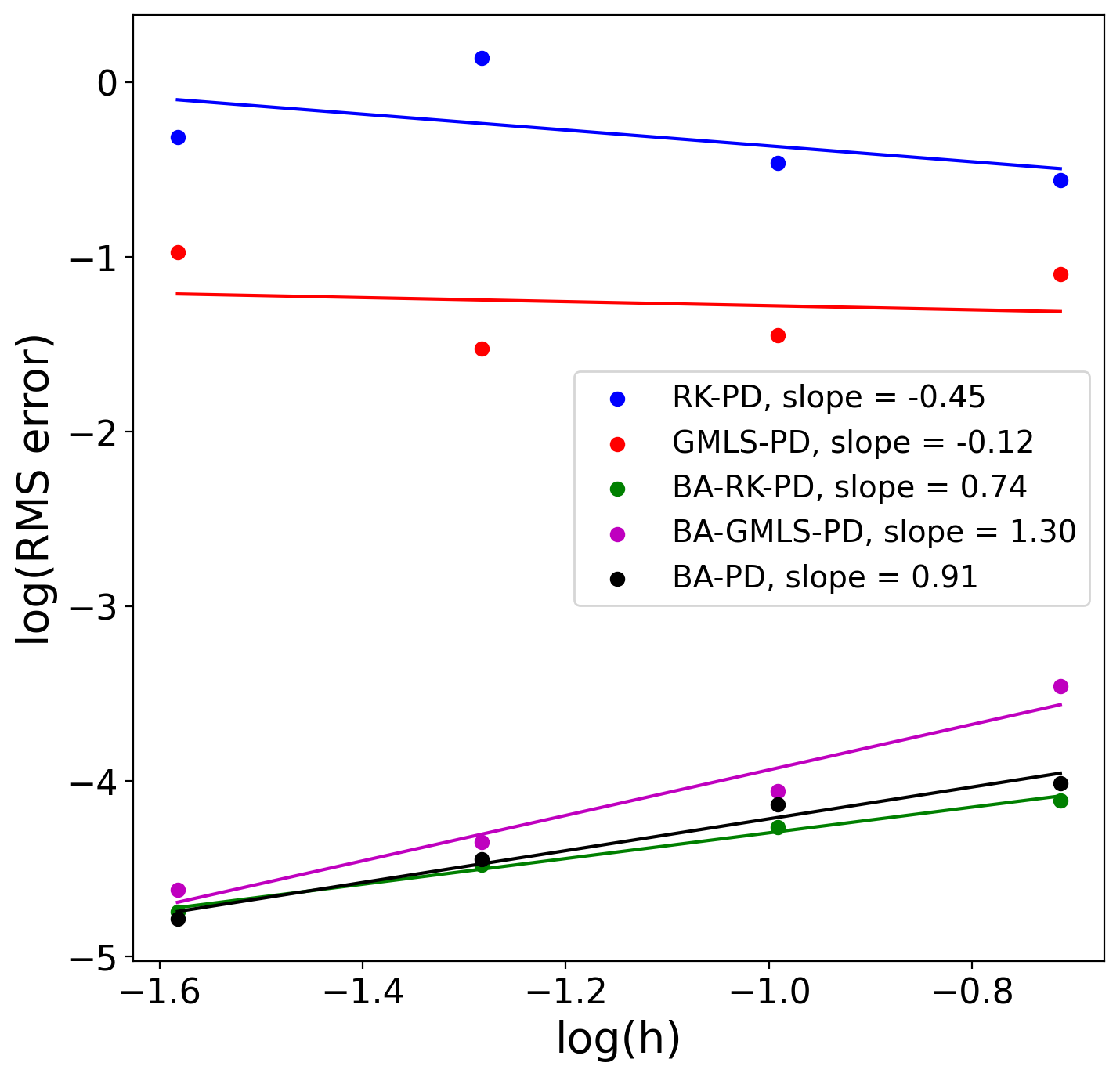}}
  \caption{Two-dimensional convergence test for the plate with a hole problem, in the mild near-incompressibility regime. Plots show the displacement RMS errors for the quadratic methods.}
  \label{fig:plate-incompressibility}
\end{figure*}

The models are also tested in the mild near-incompressibility limit ($\nu=0.495$), where similar results are obtained. \cref{fig:plate-incompressibility} shows the quadratic results where a near-linear convergence rate is observed for the bond-associated methods. For higher levels of near-incompressibility, special treatment will be required to avoid volumetric locking (see, e.g.~\citep{moutsanidis2020treatment}).

\section{Conclusions}
\label{sec:conclusions}


In this study, we developed a unified approach to the application of reproducing kernel and generalized moving least square methods to peridynamic correspondence models, under the strong-form framework. These methods were previously offered to employ higher-order corrections and improve the accuracy of the PD deformation gradient operator, which has well-known instability issues. We have investigated the robustness of these approaches by testing their linear, quadratic, and cubic versions in different problems. Our analysis demonstrates that while the RK and GMLS algorithms can be useful for obtaining higher-order convergence in smooth problems, they are insufficient for eliminating the zero-energy mode oscillations. The unstable behavior becomes more dominant under non-uniform discretizations. In addition, discrete solution instabilities are more likely to occur as the number of neighbors inside the horizon increases. 

We proposed a bond-associated formulation based on the higher-order framework, which naturally eliminates the instability issue, without introducing tunable parameters. The non-local divergence operator has been modified by utilizing a bond-level stress measure, which depends on the bond-level kinematic variable (deformation gradient). We demonstrated that the bond-associated correction is necessary to achieve a stable solution, and the improvements gained through the stabilization correction is more tangible than the benefits of the higher-order corrections. Our results indicate that the correspondence formulation is the most robust when the bond-associative approach is combined with the higher-order gradient operators. It is shown that the higher-order, bond-associated model can obtain second-order convergence for smooth problems and first-order convergence for problems involving discontinuities (e.g. curvilinear free surfaces). Additionally, it is demonstrated that the results hold for mild near-incompressible materials.

Future work would probe the application of the bond-associated models in solving dynamic problems involving crack growth. The presented approach may be readily implemented in other formulations such as the semi-Lagrangian peridynamics \citep{behzadinasab2020semi}. While our analysis has shown robustness of the bond-associated model with respect to the number of neighbors, further study could examine the role of horizon and influence function in the developed PD formulations. 

In Part II of this paper, we present a framework for incorporating traction (Neumann) boundary conditions in the strong-form PD. We also study wave propagation using the unified correspondence formulation. We demonstrate the applicability of the bond-associative modeling approach and the higher-order corrections to such problems, further proving the importance of both elements in achieving a robust methodology.

\begin{acknowledgements}
  Y.~Bazilevs was partially supported through the Sandia contract No. 2111577. N.~Trask acknowledges funding under the DOE ASCR PhILMS center (Grant number DE-SC001924) and the Laboratory Directed Research and Development program at Sandia National Laboratories. Sandia National Laboratories is a multi-mission laboratory managed and operated by National Technology and Engineering Solutions of Sandia, LLC., a wholly owned subsidiary of Honeywell International, Inc., for the U.S. Department of Energy’s National Nuclear Security Administration under contract DE-NA0003525.
\end{acknowledgements}

\bibliography{paper}

\end{document}